\documentclass{gtart}
\def\mailto#1{#1}
\def\MR#1{MR#1}
\let\gtemail\email
\def\gtref#1#2#3#4#5{Geom. Topol. #1 (#2) #4--#5}

\def\ifplaintex{\expandafter\ifx\csname documentclass\endcsname\relax}

\def\ifplaintex{\expandafter\ifx\csname documentclass\endcsname\relax}

\ifplaintex 
\else
\fi

\expandafter\ifx\csname epsfbox\endcsname\relax\input epsf\fi

\def\gt{{\mathsurround=0pt\it $\cal G\mskip-2mu$eometry \&\ 
$\cal T\!\!$opology}}        

\def\gtp{{\mathsurround=0pt\it $\cal G\mskip-2mu$eometry \&\ 
$\cal T\!\!$opology $\cal P\!$ublications}}  

\def\lognumber#1{\def\thelognumber{#1}}
\def\volumenumber#1{\def\thevolumenumber{#1}}
\def\papernumber#1{\def\thepapernumber{#1}}
\def\volumeyear#1{\def\thevolumeyear{#1}}

\def\pagenumbers#1#2{\def\startpage{#1}\def\finishpage{#2}}
\def\published#1{\def\publishdate{#1}}
\def\proposed#1{\def\theproposer{#1}}
\def\seconded#1{\def\theseconders{#1}}
\def\received#1{\def\receiveddate{#1}}

\def\accepted#1{\def\accepteddate{#1}}
\def\asciititle#1{\def\theasciititle{#1}}

\def\asciiaddress#1{\def\theasciiaddress{#1}}
\def\asciiemail#1{\def\theasciiemail{#1}}

\long\def\asciiabstract#1{\long\def\theasciiabstract{#1}}
\def\asciikeywords#1{\def\theasciikeywords{#1}}

\def\shorttitle#1{\def\theshorttitle{#1}}

\let\\\par\let\thelognumber\relax
\let\thevolumenumber\relax\let\thepapernumber\relax
\let\thevolumeyear\relax\let\thesamplenumber\relax\let\startpage\relax
\let\finishpage\relax\let\publishdate\relax\let\receiveddate\relax
\let\reviseddate\relax\let\accepteddate\relax\let\theasciititle\relax
\let\theasciiauthors\relax\let\theasciiaddress\relax
\let\theasciiabstract\relax\let\theasciikeywords\relax
\let\theasciiemail\relax\let\theshortauthors\relax\let\theshorttitle\relax

\long\def\maketitlep{   

\count0=\startpage

\gt\hfill      
\hbox to 77pt{\vbox to 0pt{\vglue -15pt\epsfbox{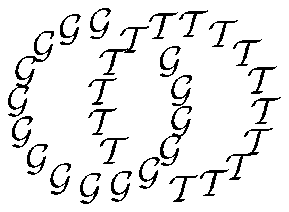}\vss}\hss}
\break
{\small\ifx\thesamplenumber\relax 
Volume \else Sample
\fi\thevolumenumber\ (\thevolumeyear)
\startpage--\finishpage\nl
Published: \publishdate}
\vglue 0.5truein plus 0.4fil minus 0.1truein

{\parskip=0pt\leftskip 0pt plus 1fil\def\\{\par\smallskip}{\ifplaintex\large
\else\Large\fi\bf\thetitle}\par\medskip}   

\vglue 0pt plus 0.1fil 

{\parskip=0pt\leftskip 0pt plus 1fil\def\\{\par}{\sc\theauthors}
\par\medskip}

\vglue 0pt plus 0.1fil 

{\small\parskip=0pt\let\newline\\
{\leftskip 0pt plus 1fil\def\\{\par}{\sl\theaddress}\par}
\expandafter\ifx\theemail\relax    
\relax\else\vglue 5pt plus 0.02fil minus 2pt\def\\{\stdspace{\rm 
and}\stdspace} 
\cl{Email:\stdspace\tt\theemail}\fi
\ifx\theurl\relax                  
\relax\else\vglue 5pt plus 0.02fil minus 2pt\def\\{\stdspace{\rm 
and}\stdspace}
\cl{URL:\stdspace\tt\theurl}\fi\par}

\vglue 7pt plus 0.3fil minus 3pt

{\bf Abstract}
\vglue 5pt plus 0.1fil minus 2pt

\theabstract

\vglue 7pt plus 0.3fil minus 3pt

{\bf AMS Classification numbers}\quad Primary:\quad \theprimaryclass

Secondary:\quad \thesecondaryclass

\vglue 5pt plus 0.3fil minus 2pt

{\bf Keywords:}\quad \thekeywords

\vglue 10pt plus 0.5fil minus 5pt

{\small  Proposed: \theproposer\hfill Received: \receiveddate\nl
Seconded: \theseconders\hfill 
\ifx\reviseddate\relax                         
Accepted: \accepteddate                        
\else
Revised: \reviseddate                          
\fi}
\eject
}       

\font\phead=cmsl9 scaled 950
\font=cmsl9 scaled 1050
\font\pnum=cmbx10 scaled 913
\font\lnum=cmbx10 
\font\pfoot=cmsl9 scaled 950
\font=cmsl9 scaled 1050
\ifplaintex
\headline{\vbox to 0pt{\vskip -4.5mm\line{\small\phead\ifnum
\count0=\startpage ISSN 1364-0380 (on line)
1465-3060 (printed) \hfill {\pnum\folio}\else\ifodd\count0\def\\{ }%
\ifx\theshorttitle\relax\thetitle\else\theshorttitle\fi\hfill{\pnum\folio}
\else\def\\{ and }{\pnum\folio}\hfill\ifx\theshortauthors\relax\theauthors
\else\theshortauthors\fi\fi\fi}\vss}}
\footline{\vbox to 0pt{\vglue 0mm\line{\small\pfoot\ifnum\count0=\startpage
\copyright\ \gtp\hfill\else
\gt, Volume \thevolumenumber\ (\thevolumeyear)\hfill\fi}\vss
}}
\else
\makeatletter
\def\@oddhead{{\small\lhead\ifnum\count0=\startpage ISSN 1364-0380 (on line)
1465-3060 (printed) \hfill {\lnum\number\count0}\else\ifodd\count0
\def\\{ }\ifx\theshorttitle\relax \thetitle \else\theshorttitle\fi\hfill
{\lnum\number\count0}\else\def\\{ and }{\lnum\number\count0}
\hfill\ifx\theshortauthors\relax 
\theauthors\else\theshortauthors\fi\fi\fi}}\def\@evenhead{\@oddhead}
\def\@oddfoot{\small\lfoot\ifnum\count0=\startpage\copyright\ \gtp\hfill\else
\gt, Volume \thevolumenumber\ (\thevolumeyear)\hfill\fi}
\def\@evenfoot{\@oddfoot}
\makeatother
\fi

\newwrite\gtoutfile
\long\gdef\makeheadfile{  
{\def\\{, }\def\s{ }
\immediate\openout\gtoutfile head.xxx
\immediate\write\gtoutfile{Proxy-for: \ifx\theasciiauthors\relax
\theauthors\else\theasciiauthors\fi\s<\ifx\theasciiemail\relax\theemail\else\theasciiemail\fi>}
\immediate\write\gtoutfile{\noexpand\\}
\immediate\write\gtoutfile{Authors: \ifx\theasciiauthors\relax
\theauthors\else\theasciiauthors\fi}
{\def\\{ }\immediate\write\gtoutfile{Title: \ifx\theasciititle\relax
\thetitle\else\theasciititle\fi}}
\immediate\write\gtoutfile{Subj-class: GT or SG or MG etc}
\immediate\write\gtoutfile{MSC-class: \theprimaryclass\ifx\thesecondaryclass\relax\else, \thesecondaryclass\fi}
\immediate\write\gtoutfile{Journal-ref: Geom. Topol. \thevolumenumber
(\thevolumeyear) \startpage-\finishpage}
\immediate\write\gtoutfile{Comments: Published by Geometry and Topology at}
\immediate\write\gtoutfile{\s\s http://www.maths.warwick.ac.uk/gt/GTVol\thevolumenumber/paper\thepapernumber.abs.html}
\immediate\write\gtoutfile{\noexpand\\}
\immediate\write\gtoutfile{}
\ifx\theasciiabstract\relax
\immediate\write\gtoutfile{\theabstract}\else
\immediate\write\gtoutfile{\theasciiabstract}\fi
\immediate\write\gtoutfile{}
\immediate\write\gtoutfile{\noexpand\\}
\immediate\write\gtoutfile{}
\immediate\closeout\gtoutfile}}  

\def\maketitlepage{\maketitlep\makeheadfile}
\let\maketitle\maketitlepage

\lognumber{289}

\volumenumber{6}
\papernumber{28} 
\volumeyear{2002}
\pagenumbers{917}{990} 
\received{22 October 2002}
\accepted{31 December 2002}
\published{31 December 2002}
\proposed{Haynes Miller}
\seconded{Ralph Cohen, Bill Dwyer}

\usepackage{amsmath,amssymb}

\def\S{section~\ignorespaces}

\usepackage[all]{xy}
\newcommand{\xymx}[2][]{\vcenter{\xymatrix#1{#2}}}

\usepackage{diagrams}

\diagramstyle[height=2em,width=2em,midshaft,labelstyle=\scriptstyle]
\newarrow{To}----{->}

\newcommand{\Z}{{\mathbb Z}}

\newcommand{\C}{{\mathbb C}}

\newcommand{\F}{{\mathbb F}}

\newcommand{\pcom}{{}_{p}^{\wedge}}
\newcommand{\doscom}{_{2}^{\wedge}}
\newcommand{\zploc}{\Z_{(p)}}

\DeclareMathAlphabet\EuR{U}{eur}{m}{n}
\SetMathAlphabet\EuR{bold}{U}{eur}{b}{n}

\newcommand{\SFL}[1][]{(S#1,\calf#1,\call#1)}

\newcommand{\Tor}{\operatorname{Tor}\nolimits}
\newcommand{\Inj}{\operatorname{Inj}\nolimits}
\newcommand{\diag}{\operatorname{diag}\nolimits}
\newcommand{\defeq}{\overset{\text{\textup{def}}}{=}}
\newcommand{\gen}[1]{{\langle}#1{\rangle}}
\newcommand{\pr}{\operatorname{pr}\nolimits}
\renewcommand{\:}{\colon\thinspace}

\newcommand{\calc}{\mathcal{C}}
\newcommand{\cale}{\mathcal{E}}
\newcommand{\calf}{\mathcal{F}}
\newcommand{\calj}{\mathcal{J}}
\newcommand{\call}{\mathcal{L}}
\newcommand{\calz}{\mathcal{Z}}

\newcommand{\Fr}{\textup{Fr}}

\newcommand{\orb}{\mathcal{O}}

\newcommand{\curs}{\EuR}
\newcommand{\Ab}{\curs{Ab}}

\newcommand{\widebar}[1]{\overset{\mskip3mu\hrulefill\mskip3mu}{#1}
                \vphantom{#1}}
\newcommand{\sminus}{\smallsetminus}
\newcommand{\nsg}{\vartriangleleft}
\newcommand{\Id}{\operatorname{Id}\nolimits}
\newcommand{\incl}{\operatorname{incl}\nolimits}
\newcommand{\Inn}{\operatorname{Inn}\nolimits}
\newcommand{\op}{^{\textup{op}}}

\let\oldcirc=\circ
\renewcommand{\circ}{\mathchoice
    {\mathbin{\scriptstyle\oldcirc}}{\mathbin{\scriptstyle\oldcirc}}
    {\mathbin{\scriptscriptstyle\oldcirc}}
    {\mathbin{\scriptscriptstyle\oldcirc}}}

\newcommand{\hclim}[1]{\setbox1=\hbox{\rm hocolim}
    \setbox2=\hbox to \wd1{\rightarrowfill} \ht2=0pt \dp2=-1pt
    \mathop{\vtop{\baselineskip=5pt\box1\box2}}
    _{#1}}

\newcommand{\invlim}[1]{\higherlim{#1}{}}

\newcommand{\map}{\operatorname{Map}\nolimits}

\newcommand{\Hom}{\operatorname{Hom}\nolimits}
\newcommand{\Rep}{\operatorname{Rep}\nolimits}
\newcommand{\Iso}{\operatorname{Iso}\nolimits}
\newcommand{\Aut}{\operatorname{Aut}\nolimits}
\newcommand{\Out}{\operatorname{Out}\nolimits}
\newcommand{\Mor}{\operatorname{Mor}\nolimits}
\newcommand{\Ker}{\operatorname{Ker}\nolimits}

\renewcommand{\Im}{\operatorname{Im}\nolimits}
\renewcommand{\rk}{\operatorname{rk}\nolimits}
\newcommand{\ztwo}{\widehat{\Z}_2}
\newcommand{\q}[1]{\widehat{#1}}

\newcommand{\longleft}[1]{\;{\leftarrow%
\count255=0 \loop \mathrel{\mkern-6mu}%
    \relbar\advance\count255 by1\ifnum\count255<#1\repeat}\;}
\newcommand{\longright}[1]{\;{\count255=0 \loop \relbar\mathrel{\mkern-6mu}%
    \advance\count255 by1\ifnum\count255<#1\repeat\rightarrow}\;}
\newcommand{\Right}[2]{\overset{#2}{\longright#1}}
\newcommand{\RIGHT}[3]{\mathrel{\mathop{\kern0pt\longright#1}
        \limits^{#2}_{#3}}}

\newcommand{\LEFT}[3]{\mathrel{\mathop{\kern0pt\longleft#1}\limits^{#2}_{#3}}
}
\newcommand{\dRIGHT}[3]{\mathrel{%
   \mathop{\vcenter{\baselineskip=0pt\hbox{$\kern0pt\longright#1$}%
   \hbox{$\kern0pt\longright#1$}}}\limits^{#2}_{#3}}}
\newcommand{\LRIGHT}[3]{\mathrel{%
   \mathop{\vcenter{\baselineskip=0pt\hbox{$\kern0pt\longleft#1$}%
   \hbox{$\kern0pt\longright#1$}}}\limits^{#2}_{#3}}}
\newcommand{\RLEFT}[3]{\mathrel{%
   \mathop{\vcenter{\baselineskip=0pt\hbox{$\kern0pt\longright#1$}%
   \hbox{$\kern0pt\longleft#1$}}}\limits^{#2}_{#3}}}
\newcommand{\onto}[1]{\;{\count255=0 \loop \relbar\joinrel
    \advance\count255 by1
    \ifnum\count255<#1 \repeat \twoheadrightarrow}\;}
\newcommand{\Onto}[2]{\overset{#2}{\onto#1}}

\newtheorem{Thm}{Theorem}[section]
\newtheorem{Prop}[Thm]{Proposition}

\newtheorem{Lem}[Thm]{Lemma}

\theoremstyle{definition}
\newtheorem{Defi}[Thm]{Definition}

\newenvironment{smallpmatrix}
                 {\bigl(\begin{smallmatrix}}{\end{smallmatrix}\bigr)}

\def\trp[#1,#2,#3]{[\hskip-1.5pt[#1,#2,#3]\hskip-1.5pt]}

\newcommand{\AAA}{\widehat{A}}
\newcommand{\BBB}{\widehat{B}}
\renewcommand{\gg}{\mathbb{G}}
\newcommand{\A}{\mathfrak{A}}
\newcommand{\B}{\mathfrak{B}}
\newcommand{\CC}{\mathfrak{C}}
\newcommand{\RR}{\mathfrak{R}}
\newcommand{\higherlim}[2]{\displaystyle\setbox1=\hbox{\rm lim}
        \setbox2=\hbox to \wd1{\leftarrowfill} \ht2=0pt \dp2=-1pt
        \setbox3=\hbox{$\scriptstyle{#1}$}
        \def\test{#1}\ifx\test\empty
        \mathop{\mathop{\vtop{\baselineskip=5pt\box1\box2}}}\nolimits^{#2}
        \else
        \ifdim\wd1<\wd3
        \mathop{\hphantom{^{#2}}\vtop{\baselineskip=5pt\box1\box2}^{#2}}_{#1}
        \else
        \mathop{\mathop{\vtop{\baselineskip=5pt\box1\box2}}_{#1}}%
        \nolimits^{#2}
        \fi\fi}
\newcommand{\hilim}[4][]{\def\test{#1}\def\tst{-}
	\ifx\test\tst{\setbox1=\hbox{\rm lim}
	\setbox2=\hbox to \wd1{\leftarrowfill} \ht2=0pt \dp2=-1pt
	\mathop{\vtop{\baselineskip=5pt\box1\box2}}\nolimits^{#3}(#4)}
	\else{\higherlim{#2}{#3}(#4)}\fi}

\newcommand{\X}{\mathfrak{X}}

\newcommand{\U}{\mathcal{U}}
\newcommand{\Spin}{\textup{Spin}}
\newcommand{\bb}{\mathfrak{b}}   
\newcommand{\Sol}{\textup{Sol}}
\def\fqbar{\widebar{\F}_q}
\def\Top{\curs{Top}}

\renewcommand{\labelenumi}{\textup{(\alph{enumi})}}%

\newcommand{\homf}{\Hom_{\calf}}
\newcommand{\autf}{\Aut_{\calf}}
\newcommand{\outf}{\Out_{\calf}}
\newcommand{\isof}{\Iso_{\calf}}

\newcommand{\sylp}[1]{\textup{Syl}_p(#1)}

\renewenvironment{enumerate}{\begin{list}%
{\labelenumi}
{\usecounter{enumi}%
\setlength{\itemindent}{0pt}%
\settowidth{\labelwidth}{\labelenumi}%
\addtolength{\labelwidth}{\labelsep}%
\setlength{\leftmargin}{\labelsep}%
\addtolength{\leftmargin}{\labelwidth}%
\setlength{\listparindent}{0pt}%
\setlength{\itemsep}{6pt}%
\setlength{\parsep}{0pt}%
\setlength{\topsep}{6pt}%
}}{\end{list}}

\def\Syl{\textup{Syl}}

\def\beq#1\eeq{\begin{equation*}#1\end{equation*}}

\newcommand{\spin}{_{\Spin}}
\newcommand{\sol}{_{\Sol}}
\newcommand{\fspin}{\calf\spin}
\newcommand{\fcspin}{\calf^c\spin}

\newcommand{\lcspin}{\call^c\spin}
\newcommand{\lccspin}{\call^{cc}\spin}
\newcommand{\ospin}{\orb\spin}
\newcommand{\ocspin}{\orb^c\spin}
\newcommand{\zspin}{\calz\spin}
\newcommand{\fsol}{\calf\sol}
\newcommand{\fcsol}{\calf^c\sol}
\newcommand{\fccsol}{\calf^{cc}\sol}
\newcommand{\lcsol}{\call^c\sol}
\newcommand{\lccsol}{\call^{cc}\sol}
\newcommand{\osol}{\orb\sol}
\newcommand{\ocsol}{\orb^c\sol}
\newcommand{\zsol}{\calz\sol}
\newcommand{\bspin}{B\Spin'}   
\newcommand{\bsol}{B\Sol'}   

\title{Construction of 2--local finite groups of a type\\studied by 
Solomon and Benson}
\asciititle{Construction of 2-local finite groups of a type studied by 
Solomon and Benson}
\shorttitle{Construction of 2--local finite groups} 

\authors{Ran Levi\\Bob Oliver}
\address{Department of Mathematical Sciences, University of 
Aberdeen\\Meston Building 339, Aberdeen AB24 3UE, UK}
\gtemail{\mailto{ran@maths.abdn.ac.uk}{\rm\qua and\qua}\mailto{bob@math.univ-paris13.fr}}
\secondaddress{LAGA -- UMR 7539 of the CNRS, Institut Galil\'ee\\Av 
J-B Cl\'ement, 93430 Villetaneuse, France}

\asciiaddress{Department of Mathematical Sciences, University of 
Aberdeen\\Meston Building 339, Aberdeen AB24 3UE, UK\\and\\LAGA - UMR 
7539 of the CNRS, Institut Galilee\\Av 
J-B Clement, 93430 Villetaneuse, France}
\asciiemail{ran@maths.abdn.ac.uk, bob@math.univ-paris13.fr} 

\asciikeywords{Classifying space, p-completion, finite groups, fusion.}

\begin{document}

\begin{abstract}  
A $p$--local finite group is an algebraic structure with a classifying 
space which has many of the properties of $p$--completed classifying spaces 
of finite groups. In this paper, we construct a family of 2--local finite 
groups, which are exotic in the following sense:  they are based on 
certain fusion systems over the Sylow 2--subgroup of $\Spin_7(q)$ ($q$ an 
odd prime power) shown by Solomon not to occur as the 2--fusion 
in any actual finite group. Thus, the resulting classifying spaces are not 
homotopy equivalent to the $2$--completed classifying space of any finite 
group. As predicted by Benson, these classifying spaces are also very 
closely related to the Dwyer--Wilkerson space $BDI(4)$.
\end{abstract}

\asciiabstract{A p-local finite group is an algebraic structure 
with a classifying space which has many of the properties of
p-completed classifying spaces of finite groups. In this paper, we
construct a family of 2-local finite groups, which are exotic in the
following sense: they are based on certain fusion systems over the
Sylow 2-subgroup of Spin_7(q) (q an odd prime power) shown by Solomon
not to occur as the 2-fusion in any actual finite group. Thus, the
resulting classifying spaces are not homotopy equivalent to the
2-completed classifying space of any finite group. As predicted by
Benson, these classifying spaces are also very closely related to the
Dwyer-Wilkerson space BDI(4).}

\maketitlepage


As one step in the classification of finite simple groups, Ron Solomon 
\cite{Solomon} considered the problem of classifying all finite simple 
groups whose Sylow 2--subgroups are isomorphic to those of the Conway
group $Co_3$. The end result of his paper was that $Co_3$ is the only
such group.  In the process of proving this, he needed to consider
groups $G$ in which all involutions are conjugate, and such that for
any involution $x\in{}G$, there are subgroups $K\nsg{}H\nsg{}C_G(x)$
such that $K$ and $C_G(x)/H$ have odd order and $H/K\cong\Spin_7(q)$
for some odd prime power $q$.  Solomon showed that such a group $G$
does not exist. The proof of this statement was also interesting, in
the sense that the 2--local structure of the group in question appeared
to be internally consistent, and it was only by analyzing its
interaction with the $p$--local structure (where $p$ is the prime of
which $q$ is a power) that he found a contradiction.

In a later paper \cite{Benson}, Dave Benson, inspired by Solomon's
work, constructed certain spaces 
which can be thought of as the 2--completed classifying spaces which the 
groups studied by Solomon would have if they existed.  He started with the 
spaces $BDI(4)$ constructed by Dwyer and Wilkerson  having the property 
that 
        $$ H^*(BDI(4);\F_2) \cong \F_2[x_1,x_2,x_3,x_4]^{GL_4(2)} $$ 
(the rank four Dickson algebra at the prime 2).  Benson then considered, for 
each odd prime power $q$, the homotopy fixed point set of the $\Z$--action 
on $BDI(4)$ generated by an ``Adams operation'' $\psi^q$ constructed by 
Dwyer and Wilkerson.  This homotopy fixed point set is denoted here 
$BDI_4(q)$.  

In this paper, we construct a family of \emph{2--local} finite groups, in 
the sense of \cite{BLO2}, which have the 2--local structure considered by 
Solomon, and whose classifying spaces are homotopy equivalent to Benson's 
spaces $BDI_4(q)$.  The results of \cite{BLO2} combined with those here 
allow us to make much more precise the statement that these spaces have 
many of the properties which the  2--completed classifying spaces of the 
groups studied by Solomon would have  had if they existed.  To explain 
what this means, we first recall some definitions.  

A \emph{fusion system} over a finite $p$--group $S$ is a category whose 
objects are the subgroups of $S$, and whose morphisms are monomorphisms of 
groups which include all those induced by conjugation by elements of $S$.  
A fusion system is \emph{saturated} if it satisfies certain axioms 
formulated by Puig \cite{Puig}, and also listed in \cite[Definition 
1.2]{BLO2} as well as at the beginning of Section 1 in this paper.  In
particular, for any  
finite group $G$ and any $S\in\sylp{G}$, the category $\calf_S(G)$ whose 
objects are the subgroups of $S$ and whose morphisms are those 
monomorphisms between subgroups induced by conjugation in $G$ is a 
saturated fusion system over $S$.

If $\calf$ is a saturated fusion system over $S$, then a subgroup 
$P\le{}S$ is called \emph{$\calf$--centric} if $C_S(P')=Z(P')$ for all $P'$ 
isomorphic to $P$ in the category $\calf$.  A \emph{centric linking system} 
associated to $\calf$ consists of a category $\call$ whose objects are the 
$\calf$--centric subgroups of $S$, together with a functor 
$\call\Right1{}\calf$ which is the inclusion on objects, is surjective on 
all morphism sets and which satisfies certain additional axioms (see 
\cite[Definition 1.7]{BLO2}).  These axioms suffice to ensure that the 
$p$--completed nerve $|\call|\pcom$ has all of the properties needed to 
regard it as a ``classifying space'' of the fusion system $\calf$.  A 
\emph{$p$--local finite group} consists of a triple $(S,\calf,\call)$, 
where $S$ is a finite $p$--group, $\calf$ is a saturated fusion system over 
$S$, and $\call$ is a linking system associated to $\calf$. The
classifying space of a $p$--local finite group $\SFL$ is the
$p$--completed nerve $|\call|\pcom$ (which is $p$--complete since
$|\call|$ is always $p$--good \cite[Proposition 1.12]{BLO2}).  For example,  
if $G$ is a finite group and $S\in\sylp{G}$, then there is an explicitly 
defined centric linking system $\call^c_S(G)$ associated to $\calf_S(G)$, 
and the  classifying space of the triple $(S,\calf_S(G),\call^c_S(G))$ is 
the space $|\call^c_S(G)|\pcom\simeq{}BG\pcom$.  

Exotic examples of $p$--local finite groups for odd primes $p$ --- 
ie, examples which do not represent actual groups --- have already been 
constructed in \cite{BLO2}, but using ad hoc methods which seemed to work 
only at odd primes.

In this paper, we first construct a fusion system $\calf_{\Sol}(q)$ (for 
any odd prime power $q$) over a 2--Sylow subgroup $S$ of $\Spin_7(q)$, with 
the properties that all elements of order 2 in $S$ are conjugate (ie, 
the subgroups they generated are all isomorphic in the category), and the 
``centralizer fusion system'' (see the beginning of Section 1) of each 
such element is isomorphic to the fusion system of $\Spin_7(q)$. We then 
show that $\calf_{\Sol}(q)$ is saturated, and has a  unique associated 
linking system $\call^c_{\Sol}(q)$.  We thus obtain a 2--local finite group 
$(S,\calf_{\Sol}(q),\call^c_{\Sol}(q))$ where by Solomon's theorem 
\cite{Solomon} (as explained in more detail in Proposition \ref{G<>Sol}), 
$\calf_{\Sol}(q)$ is not the fusion system of any finite group.  Let 
$B\Sol(q)\defeq|\call^c_{\Sol}(q)|\doscom$ denote the classifying space of 
$(S,\calf_{\Sol}(q),\call^c_{\Sol}(q))$.  Thus, $B\Sol(q)$  does not have 
the homotopy type of  $BG\doscom$ for any finite group $G$, but does have 
many of the nice properties of the 2--completed classifying space of a 
finite group (as described in \cite{BLO2}).

Relating $B\Sol(q)$ to $BDI_4(q)$ requires taking the ``union'' of the 
categories $\lcsol(q^n)$ for all  $n\ge1$. This however  is complicated by 
the fact that an inclusion of fields $\F_{p^m}\subseteq\F_{p^n}$ (ie, 
$m|n$) does not induce an inclusion of cenric linking systems.  Hence we 
have to replace the centric linking systems $\lcsol(q^n)$ by  the full 
subcategories $\lccsol(q^n)$ whose objects are those 2--subgroups which are 
centric in $\fcsol(q^\infty)=\bigcup_{n\ge1}\fcsol(q^n)$, and show that 
the inclusion induces a homotopy equivalence 
$\bsol(q^n)\defeq|\lccsol(q^n)|\doscom\simeq{}B\Sol(q^n)$.  Inclusions of 
fields do induce inclusions of these categories, so we can then define 
$\lcsol(q^\infty)\defeq\bigcup_{n\ge1}\lccsol(q^n)$, and spaces
	$$ B\Sol(q^\infty) = |\lcsol(q^\infty)|\doscom \simeq
	\Bigl(\bigcup_{n\ge1}\bsol(q^n)\Bigr){}\doscom. $$
The category $\lcsol(q^\infty)$ has an ``Adams map'' $\psi^q$ induced by 
the Frobenius automorphism $x\mapsto{}x^q$ of $\fqbar$.  We then show that 
$B\Sol(q^\infty)\simeq{}BDI(4)$, the space of Dwyer and Wilkerson 
mentioned above; and also that $B\Sol(q)$ is equivalent to the homotopy 
fixed point set of the $\Z$--action on $B\Sol(q^\infty)$ generated by 
$B\psi^q$. The space  $B\Sol(q)$ is thus equivalent to Benson's spaces
$BDI_4(q)$ for any odd prime power $q$. 

The paper is organized as follows.  Two propositions used for constructing 
saturated fusion systems, one very general and one more specialized, are 
proven in Section 1.  These are then applied in Section 2 to construct the 
fusion systems $\calf_{\Sol}(q)$, and to prove that they are
saturated. In Section 3 we prove the existence and uniqueness of a
centric linking systems associated to $\calf_{\Sol}(q)$ and study their
automorphisms. Also in Section 3 is the proof that $\calf_{\Sol}(q)$
is not the fusion system of any finite group.  The connections with the 
space $BDI(4)$ of Dwyer and Wilkerson is shown in Section 4.  Some background 
material on the spinor groups $\Spin(V,\bb)$ over fields of characteristic 
$\ne2$ is collected in an appendix.  

We would like to thank Dave Benson, Ron Solomon, and Carles Broto for their 
help while working on this paper.


\section{Constructing saturated fusion systems}

In this section, we first prove a general result which is useful for 
constructing saturated fusion systems.  This is then followed by a more 
technical result, which is designed to handle the specific 
construction in Section 2.  

We first recall some definitions from \cite{BLO2}.  A \emph{fusion system} 
over a $p$--group $S$ is a category $\calf$ whose objects are the subgroups 
of $\calf$, such that 
        $$ \Hom_S(P,Q) \subseteq \Mor_\calf(P,Q) \subseteq \Inj(P,Q) $$
for all $P,Q\le{}S$, and such that each morphism in $\calf$ factors as the 
composite of an $\calf$--isomorphism followed by an inclusion.  We write 
$\homf(P,Q)=\Mor_\calf(P,Q)$ to emphasize that the morphisms are all 
homomorphisms of groups.  We say that two subgroups $P,Q\le{}S$ are 
\emph{$\calf$--conjugate} if they are isomorphic in $\calf$.  A subgroup 
$P\le{}S$ is \emph{fully centralized} (\emph{fully normalized}) in 
$\calf$ if $|C_S(P)|\ge|C_S(P')|$ ($|N_S(P)|\ge|N_S(P')|$) for all 
$P'\le{}S$ which is $\calf$--conjugate to $P$.  A \emph{saturated fusion 
system} is a fusion system $\calf$ over $S$ which satisfies the following 
two additional conditions:
\begin{enumerate}\renewcommand{\labelenumi}{\textup{(\Roman{enumi})}}%
\item For each fully normalized subgroup $P\le{}S$, $P$ is fully 
centralized and $\Aut_S(P)\in\sylp{\Aut_\calf(P)}$.
\item For each $P\le{}S$ and each $\varphi\in\homf(P,S)$ such that 
$\varphi(P)$ is fully centralized in $\calf$, if we set
        $$ N_\varphi = \bigl\{g\in{}N_S(P) \,\big|\, 
        \varphi c_g\varphi^{-1} \in \Aut_S(\varphi(P)) \bigr\}, $$
then $\varphi$ extends to a homomorphism 
$\widebar{\varphi}\in\homf(N_\varphi,S)$.  
\end{enumerate}

For example, if $G$ is a finite group and $S\in\sylp{G}$, then the 
category $\calf_S(G)$ whose objects are the subgroups of $S$ and whose 
morphisms are the homomorphisms induced by conjugation in $G$ is a 
saturated fusion system over $S$.  A subgroup $P\le{}S$ is fully 
centralized in $\calf_S(G)$ if and only if $C_S(P)\in\sylp{C_G(P)}$, and 
$P$ is fully normalized in $\calf_S(G)$ if and only if 
$N_S(P)\in\sylp{N_G(P)}$.  

For any fusion system $\calf$ over a $p$--group $S$, and any subgroup 
$P\le{}S$, the ``centralizer fusion system'' $C_{\calf}(P)$ over $C_S(P)$ 
is defined by setting 
        $$ \Hom_{C_{\calf}(P)}(Q,Q') = \bigl\{ (\varphi|_Q) \,\big|\,
        \varphi\in\homf(PQ,PQ'),\ \varphi(Q)\le{}Q',\ \varphi|_P=\Id_P 
        \bigr\} $$
for all $Q,Q'\le{}C_S(P)$ (see \cite[Definition A.3]{BLO2} or 
\cite{Puig} for more detail).  We also write 
$C_{\calf}(g)=C_{\calf}(\gen{g})$ for $g\in{}S$.  If $\calf$ is a 
saturated fusion system and $P$ is fully centralized in $\calf$, then 
$C_{\calf}(P)$ is saturated by \cite[Proposition A.6]{BLO2} (or 
\cite{Puig}).

\begin{Prop}  \label{sat.crit.}
Let $\calf$ be any fusion system over a $p$--group $S$.  Then $\calf$ is 
saturated if and only if there is a set $\X$ of elements of order $p$ 
in $S$ such that the following conditions hold:
\begin{enumerate}  
\item Each $x\in{}S$ of order $p$ is $\calf$--conjugate to some element of 
$\X$.

\item If $x$ and $y$ are $\calf$--conjugate and $y\in\X$, then there is 
some morphism $\psi\in\homf(C_S(x),C_S(y))$ such that $\psi(x)=y$.

\item For each $x\in\X$, $C_{\calf}(x)$ is a saturated fusion system 
over $C_S(x)$.
\end{enumerate}
\end{Prop}

\begin{proof}  Throughout the proof, conditions (I) and (II) always refer 
to the conditions in the definition of a saturated fusion system, as 
stated above or in \cite[Definition 1.2]{BLO2}.

Assume first that $\calf$ is saturated, and let $\X$ be the set of all 
$x\in{}S$ of order $p$ such that $\gen{x}$ is fully centralized.  Then 
condition (a) holds by definition, (b) follows from condition (II), and 
(c) holds by \cite[Proposition A.6]{BLO2} or \cite{Puig}.

Assume conversely that $\X$ is chosen such that conditions (a--c) hold for 
$\calf$.  Define
        \begin{small}  
        $$ \U = \bigl\{(P,x)\,\big|\, P\le{}S,\ |x|=p, 
        \textup{ $x\in{}Z(P)^T$,
        some $T\in\sylp{\Aut_{\calf}(P)}$, $T\ge\Aut_S(P)$} \bigr\}, 
        $$
        \end{small}
where $Z(P)^T$ is the subgroup of elements of $Z(P)$ fixed by the 
action of $T$.  Let $\U_0\subseteq\U$ be the set of pairs 
$(P,x)$ such that $x\in\X$.  For each $1\ne{}P\le{}S$, there is some $x$ 
such that $(P,x)\in\U$ (since every action of a $p$--group on $Z(P)$ has 
nontrivial fixed set); but $x$ need not be unique.  

We first check that
        \beq (P,x)\in\U_0, \textup{ $P$ fully centralized in $C_{\calf}(x)$}
        \ \Longrightarrow\ \textup{$P$ fully centralized in $\calf$.} 
        \tag{1} \eeq
Assume otherwise:  that $(P,x)\in\U_0$ and $P$ is fully centralized in 
$C_{\calf}(x)$, but $P$ is not fully centralized in $\calf$.  Let 
$P'\le{}S$ and  $\varphi\in\isof(P,P')$ be such that $|C_S(P)|<|C_S(P')|$. 
Set $x'=\varphi(x)\le{}Z(P')$.  By (b), there exists 
$\psi\in\homf(C_S(x'),C_S(x))$ such that $\psi(x')=x$.  Set 
$P''=\psi(P')$.  Then $\psi\circ\varphi\in\Iso_{C_{\calf}(x)}(P,P'')$, and 
in particular $P''$ is $C_{\calf}(x)$--conjugate to $P$.  Also, since 
$C_S(P')\le{}C_S(x')$, $\psi$ sends $C_S(P')$ injectively into $C_S(P'')$, 
and $|C_S(P)|<|C_S(P')|\le|C_S(P'')|$.  Since $C_S(P)=C_{C_S(x)}(P)$ 
and $C_S(P'')=C_{C_S(x)}(P'')$, this contradicts the original assumption 
that $P$ is fully centralized in $C_{\calf}(x)$.  

By definition, for each $(P,x)\in\U$, $N_S(P)\le{}C_S(x)$ and hence 
$\Aut_{C_S(x)}(P)=\Aut_S(P)$.  By assumption, there is 
$T\in\sylp{\Aut_{\calf}(P)}$ such that $\tau(x)=x$ for all $\tau\in{}T$; 
ie, such that $T\le\Aut_{C_{\calf}(x)}(P)$. In particular, it follows that
	\begin{small}  
        \beq \forall (P,x)\in\U:\ \ \Aut_S(P)\in\sylp{\Aut_{\calf}(P)}
        \Longleftrightarrow
        \Aut_{C_S(x)}(P)\in\sylp{\Aut_{C_{\calf}(x)}(P)}. \tag{2} \eeq
	\end{small}%

We are now ready to prove condition (I) for $\calf$; namely, to show for 
each $P\le{}S$ fully normalized in $\calf$ that $P$ is fully centralized 
and $\Aut_S(P)\in\sylp{\autf(P)}$.  By definition, $|N_S(P)|\ge|N_S(P')|$ 
for all $P'$ $\calf$--conjugate to $P$.  Choose $x\in{}Z(P)$ such that 
$(P,x)\in\U$; and let $T\in\sylp{\autf(P)}$ be such that $T\ge\Aut_S(P)$ 
and $x\in{}Z(P)^T$.  By (a) and (b), there is an element $y\in\X$ and a 
homomorphism $\psi\in\homf(C_S(x),C_S(y))$ such that $\psi(x)=y$.  Set 
$P'=\psi{}P$, and set $T'=\psi{}T\psi^{-1}\in\sylp{\autf(T')}$.  Since 
$T\ge\Aut_S(P)$ by definition of $\U$, and $\psi(N_S(P))=N_S(P')$ by the 
maximality assumption, we see that $T'\ge\Aut_S(P')$.  Also, 
$y\in{}Z(P')^{T'}$ ($T'y=y$ since $Tx=x$), and this shows that 
$(P',y)\in\U_0$.  The maximality of $|N_S(P')|=|N_{C_S(y)}(P')|$ implies 
that $P'$ is fully normalized in $C_{\calf}(y)$.  Hence by condition (I) 
for the saturated fusion system $C_\calf(y)$, together with (1) and (2), 
$P$ fully centralized in $\calf$ and $\Aut_S(P)\in\sylp{\Aut_{\calf}(P)}$. 

It remains to prove condition (II) for $\calf$.  Fix $1\ne{}P\le{}S$ and 
$\varphi\in\homf(P,S)$ such that $P'\defeq\varphi{}P$ is fully centralized 
in $\calf$, and set
        $$ N_\varphi = \bigl\{g\in{}N_S(P)\,\big|\, 
        \varphi{}c_g\varphi^{-1}\in\Aut_S(P') \bigr\}. $$ 
We must show that $\varphi$ extends to some 
$\widebar{\varphi}\in\homf(N_\varphi,S)$.  Choose some $x'\in{}Z(P')$ of 
order $p$ which is fixed under the action of $\Aut_S(P')$, and set 
$x=\varphi^{-1}(x')\in{}Z(P)$.  For all $g\in{}N_\varphi$, 
$\varphi{}c_g\varphi^{-1}\in\Aut_S(P')$ fixes $x'$, and hence $c_g(x)=x$.  
Thus 
        \beq x\in{}Z(N_\varphi) \quad\textup{and hence}\quad 
        N_\varphi\le{}C_S(x);
        \qquad\textup{and}\qquad N_S(P')\le C_S(x'). \tag{3} \eeq
Fix $y\in\X$ which is $\calf$--conjugate to $x$ and $x'$, and choose 
        $$ \psi\in\homf(C_S(x),C_S(y))
        \qquad\textup{and}\qquad
        \psi'\in\homf(C_S(x'),C_S(y)) $$ 
such that $\psi(x)=\psi'(x')=y$.  Set $Q=\psi(P)$ and $Q'=\psi'(P')$.  
Since $P'$ is fully centralized in $\calf$, $\psi'(P')=Q'$, and 
$C_S(P')\le{}C_S(x')$, we have 
        \beq \psi'(C_{C_S(x')}(P')) = \psi'(C_S(P')) = C_S(Q')
        = C_{C_S(y)}(Q'). \tag{4} \eeq
Set $\tau=\psi'\varphi\psi^{-1}\in\isof(Q,Q')$.  By construction, 
$\tau(y)=y$, and thus $\tau\in\Iso_{C_{\calf}(y)}(Q,Q')$.  Since $P'$ is 
fully centralized in $\calf$, (4) implies that $Q'$ is fully centralized 
in $C_{\calf}(y)$.  Hence condition (II), when applied to the saturated 
fusion system $C_{\calf}(y)$, shows that $\tau$ extends to a homomorphism 
$\widebar{\tau}\in\Hom_{C_{\calf}(y)}(N_\tau,C_S(y))$, where 
        \[ N_\tau = \bigl\{g\in{}N_{C_S(y)}(Q) \,\big|\,
        \tau c_g\tau^{-1}\in\Aut_{C_S(y)}(Q') \bigr\}. \]
Also, for all $g\in{}N_\varphi\le{}C_S(x)$ (see (3)), 
        \[ c_{\widebar{\tau}(\psi(g))} = \tau{}c_{\psi(g)}\tau^{-1} = 
        (\tau\psi)c_g(\tau\psi)^{-1} = 
        (\psi'\varphi)c_g(\psi'\varphi)^{-1}=c_{\psi'(h)} \in 
        \Aut_{C_S(y)}(Q') \]
for some $h\in{}N_S(P')$ such that $\varphi{}c_g\varphi^{-1}=c_h$.  This 
shows that $\psi(N_\varphi)\le{}N_\tau$; and also (since 
$C_S(Q')=\psi'(C_S(P'))$ by (4)) that
        \[ \widebar{\tau}(\psi(N_\varphi))\le \psi'(N_{C_S(x')}(P')). \] 
We can now define
        \[ \widebar{\varphi}\defeq
        (\psi')^{-1}\circ(\widebar{\tau}\circ\psi)|_{N_\varphi}
        \in\homf(N_\varphi,S), \]
and $\widebar{\varphi}|_P=\varphi$.
\end{proof}

Proposition \ref{sat.crit.} will also be applied in a 
separate paper of Carles Broto and Jesper M\o{}ller \cite{BM} to give a 
construction of some ``exotic'' $p$--local finite groups at certain odd 
primes.


Our goal now is to construct certain saturated fusion systems, by starting 
with the fusion system
of $\Spin_7(q)$ for some odd prime power $q$, and then adding to that the 
automorphisms of some subgroup of $\Spin_7(q)$.  This is a special case
of the general problem of studying fusion systems generated by fusion
subsystems, and then showing that they are saturated.  We first fix some
notation.  If $\calf_1$ and $\calf_2$ are two fusion systems over the 
same $p$--group $S$, then $\gen{\calf_1,\calf_2}$ denotes the fusion system 
over $S$ generated by $\calf_1$ and $\calf_2$:  the smallest fusion system 
over $S$ which contains both $\calf_1$ and $\calf_2$.  More generally, if 
$\calf$ is a fusion system over $S$, and $\calf_0$ is a fusion system over 
a subgroup $S_0\le{}S$, then $\gen{\calf;\calf_0}$ denotes the fusion 
system over $S$ generated by the morphisms in $\calf$ between subgroups of 
$S$, together with morphisms in $\calf_0$ between subgroups of $S_0$ only. 
In other words, a morphism in $\gen{\calf;\calf_0}$ is a composite
        $$ P_0 \Right2{\varphi_1} P_1 \Right2{\varphi_2} P_2 \Right2{}
        \cdots \Right2{} P_{k-1} \Right2{\varphi_k} P_k, $$
where for each $i$, either $\varphi_i\in\Hom_{\calf}(P_{i-1},P_i)$, or 
$\varphi_i\in\Hom_{\calf_0}(P_{i-1},P_i)$ (and $P_{i-1},P_i\le{}S_0$).

As usual, when $G$ is a finite group and $S\in\sylp{G}$, then $\calf_S(G)$ 
denotes the fusion system of $G$ over $S$.  If $\Gamma\le\Aut(G)$ is a 
group of automorphisms which contains $\Inn(G)$, then $\calf_S(\Gamma)$ 
will denote the fusion system over $S$ whose morphisms consist of all 
restrictions of automorphisms in $\Gamma$ to monomorphisms between 
subgroups of $S$.  

The next proposition provides some fairly specialized conditions which
imply that the fusion system generated by the fusion system of a
group $G$ together with certain automorphisms of a subgroup of $G$ is 
saturated.

\begin{Prop}  \label{construct-sFs}
Fix a finite group $G$, a prime $p$ dividing $|G|$, and a Sylow 
$p$--subgroup $S\in\sylp{G}$.  Fix a normal subgroup $Z\nsg{}G$ of order 
$p$, an elementary abelian subgroup $U\nsg{}S$ of rank two containing $Z$ 
such that $C_S(U)\in\sylp{C_G(U)}$, and a subgroup $\Gamma\le\Aut(C_G(U))$ 
containing $\Inn(C_G(U))$ such that $\gamma(U)=U$ for all 
$\gamma\in\Gamma$.  Set 
        $$ S_0=C_S(U) \qquad\textup{and}\qquad 
        \calf\defeq\gen{\calf_S(G);\calf_{S_0}(\Gamma)}, $$  
and assume the following hold.
\begin{enumerate}  
\item All subgroups of order $p$ in $S$ different from $Z$ are $G$--conjugate.

\item $\Gamma$ permutes transitively the subgroups of order $p$ in $U$.

\item $\{\varphi\in\Gamma\,|\,\varphi(Z)=Z\}=\Aut_{N_G(U)}(C_G(U))$.

\item For each $E\le{}S$ which is elementary abelian of rank 
three, contains $U$, and is fully centralized in $\calf_S(G)$, 
$$ \{\alpha\in\autf(C_S(E))\,|\,\alpha(Z)=Z\} = \Aut_G(C_S(E)). $$

\item For all $E,E'\le{}S$ which are elementary abelian of rank three and 
contain $U$, if $E$ and $E'$ are $\Gamma$--conjugate, then they are 
$G$--conjugate.
\end{enumerate}
Then $\calf$ is a saturated fusion system over $S$.  Also, for 
any $P\le{}S$ such that $Z\le{}P$,
\beq \{\varphi\in\homf(P,S)\,|\,\varphi(Z)=Z\} = \Hom_G(P,S). 
\tag{1} \eeq
\end{Prop}

Proposition \ref{construct-sFs} follows from the following three lemmas.  
Throughout the proofs of these lemmas, references to points (a--e) mean 
to those points in the hypotheses of the proposition, unless otherwise
stated.

\begin{Lem}\label{9.2.1}
Under the hypotheses of Proposition \ref{construct-sFs}, 
for any $P\le{}S$ and any central subgroup $Z'\le{}Z(P)$ of order $p$, 
        \beq \textup{$Z\ne Z'\le U$ $\Longrightarrow$ 
        $\exists\,\varphi\in\Hom_\Gamma(P,S_0)$ such that $\varphi(Z')=Z$} 
        \tag{1} \eeq
and
        \beq 
        \textup{$Z'\nleq U$ $\Longrightarrow$ 
        $\exists\,\psi\in\Hom_G(P,S_0)$ such that $\psi(Z')\le U$.} 
        \tag{2} \eeq
\end{Lem}

\begin{proof}  
Note first that $Z\le{}Z(S)$, since it is a normal subgroup of order $p$ 
in a $p$--group.  

Assume $Z\ne{}Z'\le{}U$. Then $U=ZZ'$, and 
        $$ P\le{}C_S(Z')=C_S(ZZ')=C_S(U)=S_0 $$
since $Z'\le{}Z(P)$ by assumption.  By (b), there is $\alpha\in\Gamma$ 
such  that $\alpha(Z')=Z$.  Since $S_0\in\sylp{C_G(U)}$, there is 
$h\in{}C_G(U)$ such that $h{\cdot}\alpha(P){\cdot}h^{-1}\le{}S_0$; and 
since 
        $$ c_h \in \Aut_{N_G(U)}(C_G(U)) \le \Gamma $$
by (c), $\varphi\defeq c_h\circ\alpha\in\Hom_\Gamma(P,S_0)$ and sends 
$Z'$ to $Z$.

If $Z'\nleq{}U$, then by (a), there is $g\in{}G$ such that 
$gZ'g^{-1}\le{}U{\sminus}Z$.  Since $Z$ is central in $S$, $gZ'g^{-1}$ is 
central in $gPg^{-1}$, and $U$ is generated by $Z$ and $gZ'g^{-1}$, it 
follows that  $gPg^{-1}\le{}C_G(U)$.  Since $S_0\in\sylp{C_G(U)}$, there 
is $h\in{}C_G(U)$  such that $h(gPg^{-1})h^{-1}\le{}S_0$; and we can take 
$\psi=c_{hg}\in\Hom_{G}(P,S_0)$.  
\end{proof}

We are now ready to prove point (1) in Proposition \ref{construct-sFs}.

\begin{Lem}\label{9.2.2}
Assume the hypotheses of Proposition \ref{construct-sFs}, and let 
	$$ \calf=\gen{\calf_S(G);\calf_{S_0}(\Gamma)} $$ 
be the fusion system generated by $G$ and $\Gamma$.  Then for all 
$P,P'\le{}S$ which contain $Z$,
        $$ \{\varphi\in\homf(P,P') \,|\, \varphi(Z)=Z \} = \Hom_G(P,P'). $$
\end{Lem}

\begin{proof} Upon replacing $P'$ by $\varphi(P)\le{}P'$, we can assume 
that $\varphi$ is an isomorphism, and thus that it factors as a composite 
of isomorphisms 
        $$ P=P_0 \RIGHT2{\varphi_1}{\cong} P_1 \RIGHT2{\varphi_2}{\cong} 
        P_2 \RIGHT2{\varphi_3}{\cong} \cdots \RIGHT2{\varphi_{k-1}}{\cong} 
        P_{k-1} \RIGHT2{\varphi_k}{\cong} P_k=P', $$
where for each $i$, $\varphi_i\in\Hom_G(P_{i-1},P_i)$ or 
$\varphi_i\in\Hom_\Gamma(P_{i-1},P_i)$.  Let $Z_i\le{}Z(P_i)$ be the 
subgroups of order $p$ such that $Z_0=Z_k=Z$ and $Z_i=\varphi_i(Z_{i-1})$. 

To simplify the discussion, we say that a morphism in $\calf$ is of type 
$(G)$ if it is given by conjugation by an element of $G$, and of type 
$(\Gamma)$ if it is the restriction of an automorphism in $\Gamma$.  More 
generally, we say that a morphism is of type $(G,\Gamma)$ if it is the 
composite of a morphism of type $(G)$ followed by one of type $(\Gamma)$, 
etc.  We regard $\Id_P$, for all $P\le{}S$, to be of both types, even
if $P\nleq{}S_0$.  By definition, if any nonidentity isomorphism is of 
type $(\Gamma)$, then its source and image are both contained in 
$S_0=C_S(U)$.  

For each $i$, using Lemma \ref{9.2.1}, choose some 
$\psi_i\in\Hom_{\calf}(P_iU,S)$ such that $\psi_i(Z_i)=Z$.  More 
precisely, using points (1) and (2) in Lemma \ref{9.2.1}, we can choose 
$\psi_i$ to be of type $(\Gamma)$ if $Z_i\le{}U$ (the inclusion if 
$Z_i=Z$), and to be of type $(G,\Gamma)$ if $Z\nleq{}U$.  Set 
$P'_i=\psi_i(P_i)$. To keep track of the effect of morphisms on the 
subgroups $Z_i$, we write them as morphisms between pairs, as shown below. 
Thus, $\varphi$ factors as a composite of isomorphisms
        $$ (P'_{i-1},Z) \Right4{\psi_{i-1}^{-1}} (P_{i-1},Z_{i-1}) 
        \Right4{\varphi_i} (P_i,Z_i) \Right4{\psi_i} (P'_i,Z).  $$
If $\varphi_i$ is of type $(G)$, then this composite (after replacing 
adjacent morphisms of the same type by their composite) is of type 
$(\Gamma,G,\Gamma)$.  If $\varphi_i$ is of type $(\Gamma)$, then the 
composite is again of type $(\Gamma,G,\Gamma)$ if either $Z_{i-1}\le{}U$ 
or $Z_i\le{}U$, and is of type $(\Gamma,G,\Gamma,G,\Gamma)$ if neither 
$Z_{i-1}$ nor $Z_i$ is contained in $U$.  So we are reduced to assuming
that $\varphi$ is of one of these two forms.

\smallskip

\noindent\textbf{Case 1}\qua  Assume first that $\varphi$ is of type 
$(\Gamma,G,\Gamma)$; ie, a composite of isomorphisms of the form
        $$ (P_0,Z) \RIGHT3{\varphi_1}{(\Gamma)} (P_1,Z_1) 
        \RIGHT3{\varphi_2}{(G)} (P_2,Z_2) 
        \RIGHT3{\varphi_3}{(\Gamma)} (P_3,Z). $$
Then $Z_1=Z$ if and only if $Z_2=Z$ because $\varphi_2$ is of type $(G)$. 
If $Z_1=Z_2=Z$, then $\varphi_1$ and $\varphi_3$ are of type $(G)$ by (c), 
and the result follows.  

If $Z_1\ne{}Z\ne{}Z_2$, then $U=ZZ_1=ZZ_2$, and thus $\varphi_2(U)=U$.  
Neither $\varphi_1$ nor $\varphi_3$ can be the identity, so 
$P_i\le{}S_0=C_S(U)$ for all $i$ by definition of $\Hom_\Gamma(-,-)$, 
and hence $\varphi_2$ is of type $(\Gamma)$ by (c).  It follows that 
$\varphi\in\Iso_\Gamma(P_0,P_3)$ sends $Z$ to itself, and is of type $(G)$ 
by (c) again.

\smallskip

\noindent\textbf{Case 2}\qua  Assume now that $\varphi$ is of type 
$(\Gamma,G,\Gamma,G,\Gamma)$; more precisely, that it is a composite of the 
form
	\begin{small}  
        $$ { \def\2#1#2{\!\RIGHT2{#1}{#2}}
        (P_0,Z) \2{\varphi_1}{(\Gamma)} (P_1,Z_1) 
        \2{\varphi_2}{(G)} (P_2,Z_2) 
        \2{\varphi_3}{(\Gamma)} (P_3,Z_3) 
        \2{\varphi_4}{(G)} (P_4,Z_4) 
        \2{\varphi_5}{(\Gamma)} (P_5,Z), } $$
	\end{small}%
where $Z_2,Z_3\nleq{}U$.  Then $Z_1,Z_4\le{}U$ and are distinct from $Z$, 
and the groups $P_0,P_1,P_4,P_5$ all contain $U$ since $\varphi_1$ and 
$\varphi_5$ (being of type $(\Gamma)$) leave $U$ invariant.  In 
particular, $P_2$ and $P_3$ contain $Z$, since $P_1$ and $P_4$ do and 
$\varphi_2,\varphi_4$ are of type $(G)$.  We can also assume that 
$U\le{}P_2,P_3$, since otherwise $P_2\cap{}U=Z$ or $P_3\cap{}U=Z$, 
$\varphi_3(Z)=Z$, and hence $\varphi_3$ is of type $(G)$ by (c) again.  
Finally, we assume that $P_2,P_3\le{}S_0=C_S(U)$, since otherwise 
$\varphi_3=\Id$.  

Let $E_i\le{}P_i$ be the rank three elementary abelian subgroups defined 
by the requirements that $E_2=UZ_2$, $E_3=UZ_3$, and 
$\varphi_i(E_{i-1})=E_i$.  In particular, $E_i\le{}Z(P_i)$ for $i=2,3$ 
(since $Z_i\le{}Z(P_i)$, and $U\le{}Z(P_i)$ by the above remarks); and 
hence $E_i\le{}Z(P_i)$ for all $i$.  Also, $U=ZZ_4\le\varphi_4(E_3)=E_4$ 
since $\varphi_4(Z)=Z$, and thus $U=\varphi_5(U)\le{}E_5$.  Via similar 
considerations for $E_0$ and $E_1$, we see that $U\le{}E_i$ for all $i$.

Set $H=C_G(U)$ for short.  Let $\cale_3$ be the set of all elementary 
abelian subgroups $E\le{}S$ of rank three which contain $U$, and with the 
property that $C_S(E)\in\sylp{C_H(E)}$. Since $C_S(E)\leq 
C_S(U)=S_0\in\sylp{H}$, the last condition implies that $E$ is fully 
centralized in the fusion system $\calf_{S_0}(H)$. If $E\le{}S$ is any 
rank three elementary  abelian subgroup which contains $U$, then there is 
some $a\in{}H$ such that $E'=aEa^{-1}\in\cale_3$, since $\calf_{S_0}(H)$ 
is saturated and $U\nsg H$. Then 
$c_a\in\Iso_{G}(E,E')\cap\Iso_\Gamma(E,E')$ by (c).  So upon composing 
with such isomorphisms, we can assume that $E_i\in\cale_3$ for all $i$, 
and also that $\varphi_i(C_S(E_{i-1}))=C_S(E_i)$ for each $i$.  

In this way, $\varphi$ can be assumed to extend to an $\calf$--isomorphism 
$\widebar{\varphi}$ from $C_S(E_0)$ to $C_S(E_5)$ which sends $Z$ to 
itself.  By (e), the rank three subgroups $E_i$ are all $G$--conjugate to 
each other.  Choose $g\in{}G$ such that $gE_5g^{-1}=E_0$.  Then 
$g{\cdot}C_S(E_5){\cdot}g^{-1}$ and $C_S(E_0)$ are both Sylow 
$p$--subgroups of $C_G(E_0)$, so there is $h\in{}C_G(E_0)$ such that 
$(hg)C_S(E_5)(hg)^{-1}=C_S(E_0)$.  By (d), 
$c_{hg}\circ\widebar{\varphi}\in\Aut_\calf(C_S(E_0))$ is of type $(G)$;
and thus $\varphi\in\Iso_{G}(P_0,P_5)$.
\end{proof}

To finish the proof of Proposition \ref{construct-sFs}, it remains only 
to show:

\begin{Lem}\label{9.2.3}
Under the hypotheses of Proposition \ref{construct-sFs}, the fusion system 
$\calf$ generated by $\calf_S(G)$ and $\calf_{S_0}(\Gamma)$ is saturated. 
\end{Lem}

\begin{proof}  We apply Proposition \ref{sat.crit.}, by letting $\X$ be 
the set of generators of $Z$.  Condition (a) of the proposition (every 
$x\in{}S$ of order $p$ is $\calf$--conjugate to an element of $\X$) holds 
by Lemma \ref{9.2.1}.  Condition (c) holds since $C_{\calf}(Z)$ is the 
fusion system of the group $C_G(Z)$ by Lemma \ref{9.2.2}, and hence 
is saturated by \cite[Proposition 1.3]{BLO2}.

It remains to prove condition (b) of Proposition \ref{sat.crit.}.  We must 
show that if $y,z\in{}S$ are $\calf$--conjugate and $\gen{z}=Z$, then there 
is $\psi\in\homf(C_S(y),C_S(z))$ such that $\psi(y)=z$.  If $y\notin{}U$, 
then by Lemma \ref{9.2.1}(2), there is $\varphi\in\homf(C_S(y),S_0)$ such 
that $\varphi(y)\in{}U$.  If $y\in{}U{\sminus}Z$, then by Lemma 
\ref{9.2.1}(1), there is $\varphi\in\homf(C_S(y),S_0)$ such that 
$\varphi(y)\in{}Z$.  We are thus reduced to the case where $y,z\in{}Z$ 
(and are $\calf$--conjugate).  

In this case, then by Lemma \ref{9.2.2}, there is $g\in{}G$ such that 
$z=gyg^{-1}$.  Since $Z\nsg{}G$, $[G{:}C_G(Z)]$ is prime to $p$, so $S$ 
and $gSg^{-1}$ are both Sylow $p$--subgroups of $C_G(Z)$, and hence are 
$C_G(Z)$--conjugate.  We can thus choose $g$ such that $z=gyg^{-1}$ and 
$gSg^{-1}=S$.  Since $C_S(y)=C_S(z)=S$ ($Z\le{}Z(S)$ since it is a normal 
subgroup of order $p$), this shows that $c_g\in\Iso_G(C_S(y),C_S(z))$, and 
finishes the proof of (b) in Proposition \ref{sat.crit.}.
\end{proof}



\section{A fusion system of a type considered by Solomon}
\label{sect:Sol}

The main result of this section and the next is the following theorem:

\begin{Thm}  \label{Sol(q)}
Let $q$ be an odd prime power, and fix 
$S\in\Syl_2(\Spin_7(q))$.  Let $z\in{}Z(\Spin_7(q))$ be the central element 
of order 2.  Then there is a saturated fusion system 
$\calf=\calf_{\Sol}(q)$ which satisfies the following conditions:
\begin{enumerate}  
\item $C_{\calf}(z)=\calf_S(\Spin_7(q))$ as fusion systems over $S$.

\item All involutions of $S$ are $\calf$--conjugate.
\end{enumerate}
Furthermore, there is a unique centric linking system 
$\call=\call^c_{\Sol}(q)$ associated to $\calf$.
\end{Thm}

Theorem \ref{Sol(q)} will be proven in Propositions 
\ref{F_Sol(fqbar)} and \ref{Sol(q)+}.  Later, at the end of Section 
\ref{sect:Sol2}, we explain why Solomon's theorem \cite{Solomon} implies 
that these fusion systems are not the fusion systems of any finite groups, 
and hence that the spaces $B\Sol(q)$ are not homotopy equivalent to the 
2--completed classifying spaces of any finite groups. 

Background results needed for computations in $\Spin(V,\bb)$ have been 
collected in Appendix \ref{Appx:Spin}.  We focus attention here on 
$SO_7(q)$ and $\Spin_7(q)$.  In fact, since we want to compare the 
constructions over $\F_q$ with those over its field extensions, most of 
the constructions will first be made in the groups $SO_7(\fqbar)$ and 
$\Spin_7(\fqbar)$. 

We now fix, for the rest of the section, an odd prime power $q$.  It will 
be convenient to write $\Spin_7(q^\infty)\defeq\Spin_7(\fqbar)$, etc.  In 
order to make certain computations more explicit, we set 
        $$ V_\infty= M_2(\fqbar) \oplus M_2^0(\fqbar) \cong 
        (\fqbar)^7
        \qquad\textup{and}\qquad
        \bb(A,B)=\det(A)+\det(B) $$
(where $M_2^0(-)$ is the group of $(2\times2)$ matrices of trace zero),
and for each $n\ge1$ set 
$V_n=M_2(\F_{q^n})\oplus{}M_2^0(\F_{q^n})\subseteq{}V_\infty$.  Then 
$\bb$ is a nonsingular quadratic form on $V_\infty$ and on $V_n$.  
Identify $SO_7(q^\infty)=SO(V_\infty,\bb)$ and $SO_7(q^n)=SO(V_n,\bb)$, 
and similarly for $\Spin_7(q^n)\le\Spin_7(q^\infty)$.  For all 
$\alpha\in\Spin(M_2(\fqbar),\det)$ and 
$\beta\in\Spin(M_2^0(\fqbar),\det)$, we write $\alpha\oplus\beta$ for 
their image in $\Spin_7(q^\infty)$ under the natural homomorphism
        $$ \iota_{4,3} \: \Spin_4(q^\infty)\times\Spin_3(q^\infty) 
        \Right4{} \Spin_7(q^\infty). $$
There are isomorphisms
        $$ \widetilde{\rho}_4 \: SL_2(q^\infty)\times SL_2(q^\infty) 
        \Right1{\cong} \Spin_4(q^\infty)
        \quad\text{and}\quad
        \widetilde{\rho}_3 \: SL_2(q^\infty) \Right1{\cong} 
        \Spin_3(q^\infty) $$
which are defined explicitly in Proposition \ref{SO3-4}, and which restrict 
to isomorphisms 
	$$ SL_2(q^n)\times{}SL_2(q^n)\cong\Spin_4(q^n)
	\qquad\textup{and}\qquad
	SL_2(q^n)\cong\Spin_3(q^n) $$ 
for each $n$.  Let 
        $$ z = \widetilde{\rho}_4(-I,-I)\oplus1 =
        1 \oplus \widetilde{\rho}_3(-I) \in Z(\Spin_7(q)) $$
denote the central element of order two, and set 
        $$ z_1=\widetilde{\rho}_4(-I,I)\oplus 1 \in \Spin_7(q). $$
Here, $1\in\Spin_k(q)$ ($k=3,4$) denotes the identity element.  Define
$U=\gen{z,z_1}$.  

\begin{Defi}  \label{Gamma}
Define
        $$ \omega \: SL_2(q^\infty)^3 \Right4{} \Spin_7(q^\infty) $$
by setting
        $$ \omega(A_1,A_2,A_3)=\widetilde{\rho}_4(A_1,A_2) \oplus
        \widetilde{\rho}_3(A_3) $$
for $A_1,A_2,A_3\in{}SL_2(q^\infty)$.  Set
        $$ H(q^\infty)=\omega(SL_2(q^\infty)^3)
        \qquad\textup{and}\qquad
        \trp[A_1,A_2,A_3]=\omega(A_1,A_2,A_3) \,. $$
\end{Defi}

Since $\widetilde{\rho}_3$ and $\widetilde{\rho}_4$ are isomorphisms, 
$\Ker(\omega)=\Ker(\iota_{4,3})$, and thus
	$$ \Ker(\omega) = \gen{(-I,-I,-I)}. $$
In particular, $H(q^\infty)\cong(SL_2(q^\infty)^3)/\{\pm(I,I,I)\}$.  Also, 
	$$ z =\trp[I,I,-I] \qquad\textup{and}\qquad z_1=\trp[-I,I,I], $$
and thus 
	$$ U=\bigl\{\trp[\pm{}I,\pm{}I,\pm{}I]\bigr\} $$ 
(with all combinations of signs).  

For each $1\le{}n<\infty$, the natural homomorphism 
        $$ \Spin_7(q^n)\Right5{}SO_7(q^n) $$ 
has kernel and cokernel both of order 2.  The image of this homomorphism 
is the commutator subgroup $\Omega_7(q^n)\nsg{}SO_7(q^n)$, which is partly 
described by Lemma \ref{SO-invol}(a).  In contrast, since all elements of 
$\fqbar$ are squares, the natural homomorphism from $\Spin_7(q^\infty)$ to 
$SO_7(q^\infty)$ is surjective.

\begin{Lem} \label{tau}
There is an element $\tau\in{}N_{\Spin_7(q)}(U)$ of order 2 such that
        \beq \tau{\cdot}\trp[A_1,A_2,A_3]{\cdot}\tau^{-1} = \trp[A_2,A_1,A_3] 
        \tag{1} \eeq
for all $A_1,A_2,A_3\in{}SL_2(q^\infty)$.  
\end{Lem}

\begin{proof}  Let $\widebar{\tau}\in{}SO_7(q)$ be the involution 
defined by setting
        $$ \widebar{\tau}(X,Y)=(-\theta(X),-Y) $$
for $(X,Y)\in{}V_\infty=M_2(\fqbar)\oplus{}M_2^0(\fqbar)$, where
        $$ \theta\bigl(\begin{smallmatrix}a&b\\c&d\end{smallmatrix}\bigr)
        =\bigl(\begin{smallmatrix}d&-b\\-c&a\end{smallmatrix}\bigr). $$  
Let $\tau\in\Spin_7(q^\infty)$ be a lifting of $\widebar{\tau}$.  The 
$(-1)$--eigenspace of $\widebar{\tau}$ on $V_\infty$ has orthogonal basis 
	$$ \bigl\{(I,0)\,,
\bigl(0,\bigl(\begin{smallmatrix}1&0\\0&-1\end{smallmatrix}\bigr)\bigr)\,,
\bigl(0,\bigl(\begin{smallmatrix}0&1\\1&0\end{smallmatrix}\bigr)\bigr)\,,
\bigl(0,\bigl(\begin{smallmatrix}0&1\\-1&0\end{smallmatrix}\bigr)\bigr)
	\bigr\}, $$
and in particular has discriminant 1 with respect to this basis.  Hence by 
Lemma \ref{SO-invol}(a), $\widebar{\tau}\in\Omega_7(q)$, and so 
$\tau\in\Spin_7(q)$.  Since in addition, the $(-1)$--eigenspace of 
$\widebar{\tau}$ is 4--dimensional, Lemma \ref{SO-invol}(b) applies to show 
that $\tau^2=1$.

By definition of the isomorphisms $\widetilde{\rho}_3$ and 
$\widetilde{\rho}_4$, for all $A_i\in{}SL_2(q^\infty)$ ($i=1,2,3$) and all 
$(X,Y)\in{}V_\infty$,
	$$ \trp[A_1,A_2,A_3](X,Y)= (A_1XA_2^{-1},A_3YA_3^{-1}). $$
Here, $\Spin_7(q^\infty)$ acts on $V_\infty$ via its projection to 
$SO_7(q^\infty)$.  Also, for all $X,Y\in{}M_2(\fqbar)$,
        \beq \theta(X) = 
        \begin{smallpmatrix}0&1\\-1&0\end{smallpmatrix}{\cdot}X^t{\cdot}
        \begin{smallpmatrix}0&1\\-1&0\end{smallpmatrix}^{-1}
        \qquad\textup{and in particular}\qquad
        \theta(XY)=\theta(Y){\cdot}\theta(X);  \eeq
and $\theta(X)=X^{-1}$ if $\det(X)=1$.  Hence for all 
$A_1,A_2,A_3\in{}SL_2(q^\infty)$ and all $(X,Y)\in{}V_\infty$,
        \begin{align*} 
        \bigl(\tau{\cdot}\trp[A_1,A_2,A_3]{\cdot}\tau^{-1}\bigr)(X,Y) 
        & = \tau(-A_1{\cdot}\theta(X){\cdot}A_2^{-1},-A_3YA_3^{-1}) \\
        & = (A_2XA_1^{-1},A_3YA_3^{-1}) = \trp[A_2,A_1,A_3](X,Y).
        \end{align*}
This shows that (1) holds modulo $\gen{z}=Z(\Spin_7(q^\infty))$.  We thus 
have two automorphisms of 
$H(q^\infty)\cong(SL_2(q^\infty)^3)/\{\pm(I,I,I)\}$ --- conjugation by 
$\tau$ and the permutation automorphism --- which are liftings of the same 
automorphism of $H(q^\infty)/\gen{z}$.  Since $H(q^\infty)$ is perfect, 
each automorphism of $H(q^\infty)/\gen{z}$ has at most one lifting to an 
automorphism of $H(q^\infty)$, and thus (1) holds.  Also, since $U$ is the 
subgroup of all elements $\trp[\pm{}I,\pm{}I,\pm{}I]$ with all 
combinations of signs, formula (1) shows that $\tau\in{}N_{\Spin_7(q)}(U)$.
\end{proof}

\newcommand{\Sthree}{\widehat{\Sigma}_3}

\begin{Defi} \label{Gamma-def}
For each $n\ge1$, set
        $$ H(q^n)=H(q^\infty)\cap\Spin_7(q^n) 
        \qquad\textup{and}\qquad
        H_0(q^n)=\omega(SL_2(q^n)^3)\le{}H(q^n). $$
Define
	$$ \Gamma_n = \Inn(H(q^n)) \rtimes \Sthree \le\Aut(H(q^n)), $$ 
where $\Sthree$ denotes the group of permutation automorphisms
	$$ \Sthree = \bigl\{
	\trp[A_1,A_2,A_3]\mapsto \trp[A_{\sigma1},A_{\sigma2},A_{\sigma3}]
	\,\big|\, \sigma\in\Sigma_3 \bigr\} \le \Aut(H(q^n)) \,. $$
\end{Defi}

For each $n$, let $\psi^{q^n}$ be the automorphism of $\Spin_7(q^\infty)$ 
induced by the field isomorphism $(q\mapsto{}q^{p^n})$.  By Lemma 
\ref{Galois-Spin}, $\Spin_7(q^n)$ is the fixed subgroup of $\psi^{q^n}$.  
Hence each element of $H(q^n)$ is of the form $\trp[A_1,A_2,A_3]$, where 
either $A_i\in{}SL_2(q^n)$ for each $i$ (and the element lies in 
$H_0(q^n)$), or $\psi^{q^n}(A_i)=-A_i$ for each $i$.  This shows that 
$H_0(q^n)$ has index 2 in $H(q^n)$.  

The goal is now to choose compatible Sylow subgroups 
$S(q^n)\in\Syl_2(\Spin_7(q^n))$ (all $n\ge1$) contained in 
$N(H(q^n))$, and let $\calf_{\Sol}(q^n)$ be the fusion system 
over $S(q^n)$ generated by conjugation in $\Spin_7(q^n)$ and by restrictions 
of $\Gamma_n$.

\begin{Prop}  \label{omega}
The following hold for each $n\ge1$.
\begin{enumerate}  
\item  $H(q^n)=C_{\Spin_7(q^n)}(U)$.

\item $N_{\Spin_7(q^n)}(U)=N_{\Spin_7(q^n)}(H(q^n))=H(q^n){\cdot}\gen{\tau}$, 
and 
contains a Sylow 2--subgroup of $\Spin_7(q^n)$.
\end{enumerate}
\end{Prop}

\begin{proof}  Let $\widebar{z}_1\in{}SO_7(q)$ be the image of 
$z_1\in\Spin_7(q)$.  Set $V_-=M_2(\fqbar)$ and $V_+=M_2^0(\fqbar)$:  the 
eigenspaces of $\widebar{z}_1$ acting on $V$.  By Lemma \ref{SO-invol}(c), 
	$$ C_{\Spin_7(q^\infty)}(U) = C_{\Spin_7(q^\infty)}(z_1) $$
is the group of all elements $\alpha\in\Spin_7(q^\infty)$ whose image 
$\widebar{\alpha}\in{}SO_7(q^\infty)$ has the form 
	$$ \widebar{\alpha} = \alpha_-\oplus\alpha_+
	\qquad\textup{where}\qquad
	\alpha_{\pm}\in{}SO(V_{\pm}). $$
In other words,
	$$ C_{\Spin_7(q^\infty)}(U) = 
	\iota_{4,3}\bigl(\Spin_4(q^\infty) \times \Spin_3(q^\infty)\bigr) 
	= \omega(SL_2(q^\infty)^3) = H(q^\infty). $$
Furthermore, since 
	$$ \tau{}z_1\tau^{-1}= \tau\trp[-I,I,I]\tau^{-1} = \trp[I,-I,I] = 
	zz_1 $$
by Lemma \ref{tau}, and since any element of $N_{\Spin_7(q^\infty)}(U)$ 
centralizes $z$, conjugation by $\tau$ generates 
$\Out_{\Spin_7(q^\infty)}(U)$.  Hence
	$$ N_{\Spin_7(q^\infty)}(U) = H(q^\infty){\cdot}\gen{\tau}. $$
Point (a), and the first part of point (b), now follow upon taking 
intersections with $\Spin_7(q^n)$.  

If $N_{\Spin_7(q^n)}(U)$ did not contain a Sylow 2--subgroup of 
$\Spin_7(q^n)$, then since every noncentral involution of $\Spin_7(q^n)$ 
is conjugate to $z_1$ (Proposition \ref{elem.abel.}), the Sylow 
2--subgroups of $\Spin_7(q)$ would have no normal subgroup isomorphic to 
$C_2^2$.  By a theorem of Hall (cf \cite[Theorem 5.4.10]{Gorenstein}), 
this would imply that they are cyclic, dihedral, quaternion, or 
semidihedral.  This is clearly not the case, so $N_{\Spin_7(q^n)}(U)$ must 
contain a Sylow 2--subgroup of $\Spin_7(q)$, and this finishes the proof of 
point (b).

Alternatively, point (b) follows from the standard formulas for the orders of 
these groups (cf \cite[pages 19,140]{Taylor}), which show that
        $$ \frac{|\Spin_7(q^n)|}{|H(q^n){\cdot}\gen{\tau}|}=
        \frac{q^{9n}(q^{6n}-1)(q^{4n}-1)(q^{2n}-1)}
        {2{\cdot}[q^n(q^{2n}-1)]^3} 
        = q^{6n}(q^{4n}+q^{2n}+1)\Bigl(\frac{q^{2n}+1}2\Bigr) $$
is odd.
\end{proof}

We next fix, for each $n$, a Sylow 2--subgroup of $\Spin_7(q^n)$ which is 
contained in $H(q^n){\cdot}\gen{\tau}=N_{\Spin_7(q^n)}(U)$.  

\begin{Defi} \label{S-def}
Fix elements $A,B\in{}SL_2(q)$ such that $\gen{A,B}\cong{}Q_8$ (a 
quaternion group of order 8), and set 
$\AAA=\trp[A,A,A]$ and $\BBB=\trp[B,B,B]$.  Let 
$C(q^\infty)\le{}C_{SL_2(q^\infty)}(A)$ be the subgroup of elements of 
2--power order in the centralizer (which is abelian), and set
$Q(q^\infty)=\gen{C(q^\infty),B}$.  Define
	$$ S_0(q^\infty) = \omega(Q(q^\infty)^3) \le H_0(q^\infty) $$
and
        $$ S(q^\infty)= S_0(q^\infty){\cdot}\gen{\tau} \le H(q^\infty)
        \le \Spin_7(q^\infty). $$
Here, $\tau\in\Spin_7(q)$ is the element of Lemma \ref{tau}.  Finally, 
for each $n\ge1$, define
	\begin{multline*}  
	C(q^n)=C(q^\infty)\cap{}SL_2(q^n),\qquad
	Q(q^n)=Q(q^\infty)\cap{}SL_2(q^n), \\
	S_0(q^n) = S_0(q^\infty)\cap \Spin_7(q^n),
	\qquad\textup{and}\qquad
        S(q^n) = S(q^\infty)\cap \Spin_7(q^n). 
	\end{multline*}
\end{Defi}

Since the two eigenvalues of $A$ are distinct, its centralizer in 
$SL_2(q^\infty)$ is conjugate to the subgroup of diagonal matrices, which 
is abelian.  Thus $C(q^\infty)$ is conjugate to the subgroup of diagonal 
matrices of 2--power order.  This shows that each finite subgroup of 
$C(q^\infty)$ is cyclic, and that each finite subgroup of $Q(q^\infty)$ is 
cyclic or quaternion.  

\begin{Lem} 
For all $n$, $S(q^n)\in\Syl_2(\Spin_7(q^n))$.  
\end{Lem}

\begin{proof}  By \cite[6.23]{Suzuki}, $A$ is contained in a cyclic 
subgroup of order $q^n-1$ or $q^n+1$ (depending on which of them is 
divisible by 4).  Also, the normalizer of this cyclic subgroup is a 
quaternion group of order $2(q^n\pm1)$, and the formula 
$|SL_2(q^n)|=q^n(q^{2n}-1)$ shows that this quaternion group has odd 
index.  Thus by construction, $Q(q^n)$ is a Sylow 2--subgroup of 
$SL_2(q^n)$.  Hence $\omega(Q(q^n)^3)$ is a Sylow 2--subgroup of $H_0(q^n)$, 
so $\omega(Q(q^\infty)^3)\cap\Spin_7(q^n)$ is a Sylow 2--subgroup of 
$H(q^n)$. It follows that $S(q^n)$ is a Sylow 2--subgroup of 
$H(q^n){\cdot}\gen{\tau}$, and hence also of $\Spin_7(q^n)$ by Proposition 
\ref{omega}(b).
\end{proof}

Following the notation of Definition \ref{types}, we say that an 
elementary abelian 2--subgroup $E\le\Spin_7(q^n)$ has type I if its 
eigenspaces all have square discriminant, and has type II otherwise.  
Let $\cale_r$ be the set of elementary abelian subgroups of rank $r$ in 
$\Spin_7(q^n)$ which contain $z$, and let $\cale_r^{I}$ and $\cale_r^{II}$ 
be the sets of those of type I or II, respectively.  In Proposition 
\ref{elem.abel.}, we show that there are two conjugacy classes of 
subgroups in $\cale_4^I$ and one conjugacy class of subgroups in 
$\cale_4^{II}$.  In Proposition \ref{E4II}, an invariant 
$x_{\calc}(E)\in{}E$ is defined, for all $E\in\cale_4$ (and where $\calc$ 
is one of the conjugacy classes in $\cale_4^I$) as a tool for determining 
the conjugacy class of a subgroup.  More precisely, $E$ has type I if 
and only if $x_{\calc}(E)\in\gen{z}$, and $E\in\calc$ if and only if 
$x_{\calc}(E)=1$.  The next lemma provides some more detailed information 
about the rank four subgroups and these invariants.

Recall that we define $\AAA=\trp[A,A,A]$ and $\BBB=\trp[B,B,B]$.

\newcommand{\EUfour}{\cale^U_4}

\begin{Lem} \label{E4-props}
Fix $n\ge1$, set $E_*=\gen{z,z_1,\AAA,\BBB}\le{}S(q^n)$, and let $\calc$ 
be the $\Spin_7(q^n)$--conjugacy class of $E_*$.  Let $\EUfour$ be the set 
of all elementary abelian subgroups $E\le{}S(q^n)$ of rank $4$ which contain 
$U=\gen{z,z_1}$.  Fix a generator $X\in{}C(q^n)$ (the 2--power torsion in 
$C_{SL_2(q^n)}(A)$), and choose $Y\in{}C(q^{2n})$ such that $Y^2=X$.  Then 
the following hold.
\begin{enumerate}  
\item $E_*$ has type I.

\item $\EUfour = \bigl\{ E_{ijk},E'_{ijk} \,|\, i,j,k\in\Z \bigr\}$ (a 
finite set), where 
        $$ E_{ijk}=\gen{z,z_1,\AAA,\trp[X^iB,X^jB,X^kB]} $$
and
        $$ E'_{ijk}=\gen{z,z_1,\AAA,\trp[X^iYB,X^jYB,X^kYB]}. $$

\item $x_{\calc}(E_{ijk})=\trp[(-I)^i,(-I)^j,(-I)^k]$ and
$x_{\calc}(E'_{ijk})=\trp[(-I)^i,(-I)^j,(-I)^k]{\cdot}\AAA$.

\item All of the subgroups $E'_{ijk}$ have type II.  The subgroup 
$E_{ijk}$ has type I if and only if $i\equiv{}j$ (mod $2$), and lies in 
$\calc$ (is conjugate to $E_*$) if and only if $i\equiv{}j\equiv{}k$ (mod 
$2$).  The subgroups $E_{000}$, $E_{001}$, and $E_{100}$ thus represent 
the three conjugacy classes of rank four elementary abelian subgroups of 
$\Spin_7(q^n)$ (and $E_*=E_{000}$).

\item For any $\varphi\in\Gamma_n\le\Aut(H(q^n))$ (see Definition 
\ref{Gamma-def}), if $E',E''\in\EUfour$ are such that $\varphi(E')=E''$, 
then $\varphi(x_{\calc}(E'))=x_{\calc}(E'')$.
\end{enumerate}
\end{Lem}

\begin{proof}  
\textbf{(a) } The set
	$$ \bigl\{ (I,0)\,, (A,0)\,,(B,0)\,,(AB,0)\,,
	(0,A)\,,(0,B)\,,(0,AB) \bigr\} $$
is a basis of eigenvectors for the action of $E_*$ on 
$V_n=M_2(\F_{q^n})\oplus{}M_2^0(\F_{q^n})$. (Since the matrices $A$, $B$, 
and $AB$ all have order 4 and determinant one, each has as eigenvalues the 
two distinct fourth roots of unity, and hence they all have trace zero.)  
Since all of these have determinant one, $E_*$ has type I by definition.

\smallskip

\noindent\textbf{(b) } Consider the subgroups
        $$ R_0=\omega(C(q^\infty)^3)\cap{}S(q^n) =
        \bigl\{\trp[X^i,X^j,X^k],\trp[X^iY,X^jY,X^kY]\,\big|\,
        i,j,k\in\Z\bigr\} $$
and
        $$ R_1=C_{S(q^n)}(\gen{U,\AAA})=R_0{\cdot}\gen{\BBB}. $$
Clearly, each subgroup $E\in\EUfour$ is contained in 
	$$ C_{S(q^n)}(U) = S_0(q^n) =
	R_0{\cdot}\gen{\trp[B^i,B^j,B^k]}. $$
All involutions in this subgroup are contained in 
$R_1=R_0{\cdot}\gen{\trp[B,B,B]}$, and thus $E\le{}R_1$.  Hence 
$E\cap{}R_0$ has rank $3$, which implies that $E\ge\gen{z,z_1,\AAA}$ (the 
2--torsion in $R_0$).  Since all elements of order two in the coset 
$R_0{\cdot}\BBB$ have the form
	$$ \trp[X^iB,X^jB,X^kB] \qquad\textup{or}\qquad 
	\trp[X^iYB,X^jYB,X^kYB] $$
for some $i,j,k$, this shows that $E$ must be one of the groups $E_{ijk}$ 
or $E'_{ijk}$.  (Note in particular that $E_*=E_{000}$.)

\smallskip

\noindent\textbf{(c) } By Proposition \ref{E4II}(a), the element 
$x_{\calc}(E)\in{}E$ is characterized uniquely by the property that 
$x_{\calc}(E)=g^{-1}\psi^{q^n}(g)$ for some $g\in\Spin_7(q^\infty)$ such 
that $gEg^{-1}\in\calc$.  We now apply this explicitly to the subgroups 
$E_{ijk}$ and $E'_{ijk}$.

For each $i$, $Y^{-i}(X^iB)Y^i=Y^{-2i}X^iB=B$.  Hence for each $i,j,k$,
        $$ \trp[Y^i,Y^j,Y^k]^{-1}{\cdot}E_{ijk} {\cdot}\trp[Y^i,Y^j,Y^k] = 
        E_* $$
and
        $$ \psi^{q^n}(\trp[Y^i,Y^j,Y^k]) = 
        \trp[Y^i,Y^j,Y^k]{\cdot}\trp[(-I)^i,(-I)^j,(-I)^k]. $$
Hence
        \beq x_{\calc}(E_{ijk})=\trp[(-I)^i,(-I)^j,(-I)^k]. \eeq
Similarly, if we choose $Z\in{}C_{SL_2(q^\infty)}(A)$ such that $Z^2=Y$, 
then for each $i$, 
	$$ (Y^iZ)^{-1}(X^iYB)(Y^iZ)=B. $$
Hence for each $i,j,k$,
        $$ \trp[Y^iZ,Y^jZ,Y^kZ]^{-1}{\cdot}E'_{ijk} 
        {\cdot}\trp[Y^iZ,Y^jZ,Y^kZ] = 
        E_*. $$
Since $\psi^{q^n}(Z)=\pm{}ZA$,
        $$ \psi^{q^n}(\trp[Y^iZ,Y^jZ,Y^kZ])=
        \trp[Y^iZ,Y^jZ,Y^kZ]{\cdot}\trp[(-I)^iA,(-I)^jA,(-I)^kA], $$
and hence 
        \beq x_{\calc}(E'_{ijk})=\trp[(-I)^iA,(-I)^jA,(-I)^kA]. \eeq

\smallskip

\noindent\textbf{(d) } This now follows immediately from point (c) and 
Proposition \ref{E4II}(b,c).

\smallskip

\noindent\textbf{(e) } By Definition \ref{Gamma-def}, $\Gamma_n$ is 
generated by $\Inn(H(q^n))$ and the permutations of the three factors in 
$H(q^\infty)\cong(SL_2(q^\infty)^3)/\{\pm(I,I,I)\}$.  If 
$\varphi\in\Gamma_n$ is a permutation automorphism, then it permutes the 
elements of $\EUfour$, and preserves the elements $x_{\calc}(-)$ by the 
formulas in (c).  If $\varphi\in\Inn(H(q^n))$ and $\varphi(E')=E''$ for 
$E',E''\in\EUfour$, then $\varphi(x_{\calc}(E'))=x_{\calc}(E'')$ by 
definition of $x_{\calc}(-)$; and so the same property holds for all 
elements of $\Gamma_n$.
\end{proof}

Following the notation introduced in Section 1, $\Hom_{\Spin_7(q^n)}(P,Q)$ 
(for $P,Q\le{}S(q^n)$) denotes the set of homomorphisms from $P$ to 
$Q$ induced by conjugation by some element of $\Spin_7(q^n)$.  Also, if 
$P,Q\le{}S(q^n)\cap{}H(q^n)$, $\Hom_{\Gamma_n}(P,Q)$ denotes the set of 
homomorphisms induced by restriction of an element of $\Gamma_n$.  Let 
$\calf_n=\calf_{\Sol}(q^n)$ be the fusion system over $S(q^n)$ generated by 
$\Spin_7(q^n)$ and $\Gamma_n$.  In other words, for each $P,Q\le{}S(q^n)$, 
$\Hom_{\calf_n}(P,Q)$ is the set of all composites 
        $$ P = P_0 \Right2{\varphi_1} P_1 \Right2{\varphi_2} P_2
        \Right2{} \cdots \Right2{}
        P_{k-1} \Right2{\varphi_k} P_k = Q, $$
where $P_i\le{}S(q^n)$ for all $i$, and each $\varphi_i$ lies in 
$\Hom_{\Spin_7(q^n)}(P_{i-1},P_i)$ or 
(if $P_{i-1},P_i\le{}H(q^n)$) $\Hom_{\Gamma_n}(P_{i-1},P_i)$.  This clearly 
defines a fusion system over $S(q^n)$.


\begin{Prop} \label{rk3}
Fix $n\ge1$.  Let $E\le{}S(q^n)$ be an elementary abelian subgroup of rank 
$3$ which contains $U$, and such that 
	$$ C_{S(q^n)}(E)\in\Syl_2(C_{\Spin_7(q^n)}(E)). $$
Then 
        \beq \{\varphi\in\Aut_{\calf_n}(C_{S(q^n)}(E))\,|\,\varphi(z)=z\} 
        = \Aut_{\Spin_7(q^n)}(C_{S(q^n)}(E)). \tag{1} \eeq
\end{Prop}

\begin{proof}  Set 
	$$ \Spin=\Spin_7(q^n), \qquad S=S(q^n), \qquad \Gamma=\Gamma_n,
	\qquad\textup{and}\qquad \calf=\calf_n $$ 
for short.  Consider the subgroups
        $$ R_0=R_0(q^n)\defeq\omega(C(q^\infty)^3)\cap{}S
        \quad\textup{and}\quad
        R_1=R_1(q^n)\defeq C_S(\gen{U,\AAA})=\gen{R_0,\BBB}. $$
Here,  $R_0$ is generated by elements of the form $\trp[X_1,X_2,X_3]$, where 
either $X_i\in{}C(q^n)$, or $X_1=X_2=X_3=X\in{}C(q^{2n})$ and 
$\psi^{q^n}(X)=-X$.  Also, $C(q^n)\in\Syl_2(C_{SL_2(q^n)}(A))$ is cyclic 
of order $2^k\ge4$, where $2^k$ is the largest power which divides 
$q^n\pm1$; and $C(q^{2n})$ is cyclic of order $2^{k+1}$.  So 
        $$ R_0 \cong (C_{2^k})^3 \qquad\textup{and}\qquad
        R_1 = R_0 \rtimes \gen{\BBB}, $$
where $\BBB=\trp[B,B,B]$ has order 2 and acts on $R_0$ via 
$(g\mapsto{}g^{-1})$.  Note that 
        $$ \gen{U,\AAA} = \gen{\trp[\pm{}I,\pm{}I,\pm{}I],\trp[A,A,A]} 
        \cong C_2^3 $$ 
is the 2--torsion subgroup of $R_0$.  

We claim that
	\beq \textup{$R_0$ is the only subgroup of $S$ isomorphic to 
	$(C_{2^k})^3$.} \tag{2} \eeq
To see this, let $R'\le{}S$ be any subgroup isomorphic to $(C_{2^k})^3$, 
and let $E'\cong{}C_2^3$ be its 2--torsion subgroup.  Recall that for any 
$2$--group $P$, the Frattini subgroup $\Fr(P)$ is the subgroup generated by 
commutators and squares in $P$.  Thus 
	$$ E' \le \Fr(R') \le \Fr(S) \le \gen{R_0,\trp[B,B,I]} $$
(note that $\trp[B,B,I]=(\tau{\cdot}\trp[B,I,I])^2$).  Any elementary 
abelian subgroup of rank 4 in $\Fr(S)$ would have to contain 
$\gen{U,\AAA}$ (the 2--torsion in $R_0\cong{}C_{2^k}^3$), and this is 
impossible since no element of the coset $R_0{\cdot}\trp[B,B,I]$ commutes 
with $\AAA$.  Thus, $\rk(\Fr(S))=3$.  Hence $U\le{}E'$, since otherwise 
$\gen{U,E'}$ would be an elementary abelian subgroup of $\Fr(S)$ of rank 
$\ge4$.  This in turn implies that $R'\le{}C_S(U)$, and hence that 
$E'\le\Fr(C_S(U))\le{}R_0$.  Thus $E'=\gen{U,\AAA}$ (the 2--torsion in 
$R_0$ again).  Hence $R'\le{}C_S(\gen{U,\AAA})=\gen{R_0,\BBB}$, and it 
follows that $R'=R_0$.  This finishes the proof of (2).

Choose generators $x_1,x_2,x_3\in{}R_0$ as follows.  Fix 
$X\in{}C_{SL_2(q^\infty)}(A)$ of order $2^k$, and 
$Y\in{}C_{SL_2(q^{2n})}(A)$ of order $2^{k+1}$ such that $Y^2=X$.  Set 
$x_1=\trp[I,I,X]$, $x_2=\trp[X,I,I]$, and $x_3=\trp[Y,Y,Y]$.  Thus, 
$x_1^{2^{k-1}}=z$, $x_2^{2^{k-1}}=z_1$, and 
$(x_3)^{2^{k-1}}=\AAA$.  

Now let $E\le{}S(q^n)$ be an elementary abelian subgroup of rank $3$ which
contains $U$, and such that $C_{S(q^n)}(E)\in\Syl_2(C_{\Spin}(E))$.  
In particular, $E\le{}R_1=C_{S(q^n)}(U)$.  There are two cases to
consider:  that where $E\le{}R_0$ and that where $E\nleq{}R_0$.  

\smallskip

\noindent\textbf{Case 1: }  
Assume $E\le{}R_0$.  Since $R_0$ is abelian of rank 3, we must
have $E=\gen{U,\AAA}$, the 2--torsion subgroup of $R_0$, and $C_S(E)=R_1$.  
Also, by (2), neither $R_0$ nor $R_1$ is isomorphic to any
other subgroup of $S$; and hence
	\beq \Aut_{\calf}(R_i) = \bigl\langle 
	\Aut_{\Spin}(R_i),\Aut_{\Gamma}(R_i) \bigr\rangle
	\qquad \textup{for $i=0,1$.} \tag{4} \eeq

By Proposition \ref{elem.abel.}, $\Aut_{\Spin}(E)$ is the group of all 
automorphisms of $E$ which send $z$ to itself.  In particular, since 
$H(q^n)=C_{\Spin}(U)$, $\Aut_{H(q^n)}(E)$ is the group of all 
automorphisms of $E$ which are the identity on $U$.  Also, 
$\Gamma=\Inn(H(q^n)){\cdot}\Sthree$, where $\Sthree$ sends 
$\AAA=\trp[A,A,A]$ to itself and permutes the nontrivial elements of 
$U=\{\trp[\pm{}I,\pm{}I,\pm{}I]\}$.  Hence $\Aut_\Gamma(E)$ is the group 
of all automorphisms which send $U$ to itself.  So if we identify 
$\Aut(E)\cong{}GL_3(\Z/2)$ via the basis $\{z,z_1,\AAA\}$, then
        $$ \Aut_{\Spin}(E)=T_1\defeq
        GL^1_2(\Z/2)=\bigl\{(a_{ij})\in{}GL_3(\Z/2)\,|\,a_{21}=a_{31}=0\bigr\} 
        $$
and
        $$ \Aut_{\Gamma}(E)=T_2\defeq
        GL^2_1(\Z/2)=\bigl\{(a_{ij})\in{}GL_3(\Z/2)\,|\,
        a_{31}=a_{32}=0 \bigr\}. $$
By (2) (and since $E$ is the 2--torsion in $R_0$),
	$$ N_{\Spin}(E)=N_{\Spin}(R_0) \qquad\textup{and}\qquad
	\{\gamma\in\Gamma\,|\,\gamma(E)=E\} = 
	\{\gamma\in\Gamma\,|\,\gamma(R_0)=R_0\}. $$
Since $C_{\Spin}(E)=C_{\Spin}(R_0){\cdot}\gen{\BBB}$, the only nonidentity 
element of $\Aut_{\Spin}(R_0)$ or of $\Aut_{\Gamma}(R_0)$ which is the 
identity on $E$ is conjugation by $\BBB$, which is $-I$.  Hence 
restriction from $R_0$ to $E$ induces isomorphisms
	$$ \Aut_{\Spin}(R_0)/\{\pm I\}\cong \Aut_{\Spin}(E)
	\qquad\textup{and}\qquad
	\Aut_{\Gamma}(R_0)/\{\pm I\}\cong \Aut_{\Gamma}(E). $$
Upon identifying $\Aut(R_0)\cong{}GL_3(\Z/2^k)$ via the basis 
$\{x_1,x_2,x_3\}$, these can be regarded as sections
        $$ \mu_i\: T_i \Right4{} GL_3(\Z/2^k)/\{\pm I\} = SL_3(\Z/2^k) 
        \times \{\lambda{}I\,|\,\lambda\in(\Z/2^k)^*\}/\{\pm I\} $$
of the natural projection from $GL_3(\Z/2^k)/\{\pm{}I\}$ to $GL_3(\Z/2)$, 
which agree on the group $T_0=T_1\cap{}T_2$ of upper triangular matrices.  

We claim that $\mu_1$ and $\mu_2$ both map trivially to the second factor. 
Since this factor is abelian, it suffices to show that $T_0$ is generated 
by $[T_1,T_1]\cap{}T_0$ and $[T_2,T_2]\cap{}T_0$, and that each $T_i$ is 
generated by $[T_i,T_i]$ and $T_0$ --- and this is easily checked. (Note 
that $T_1\cong{}T_2\cong\Sigma_4$.)

By carrying out the above procedure over the field $\F_{q^{2n}}$, we see 
that both of these sections $\mu_i$ can be lifted further to 
$SL_3(\Z/2^{k+1})$ (still agreeing on $T_0$).  So by Lemma \ref{liftAut}, 
there is a section 
	$$ \mu \: GL_3(\Z/2) \Right4{} SL_3(\Z/2^k) $$
which extends both $\mu_1$ and $\mu_2$.  By (4),
$\Aut_{\calf}(R_0)=\Im(\mu){\cdot}\gen{-I}$.  

We next identify $\Aut_{\calf}(R_1)$.  By Lemma \ref{E4-props}(a), 
$E_*\defeq\gen{z,z_1,\AAA,\BBB}\le\Spin_7(q^n)$ is a subgroup of rank 4 
and type I.  So by Proposition \ref{elem.abel.}, $\Aut_{\Spin}(E_*)$ 
contains all automorphisms of $E_*\cong{}C_2^4$ which send 
$z\in{}Z(\Spin)$ to itself.  Hence for any $x\in{}N_{\Spin}(R_1)$, since 
$c_x(z)=z$, there is $x_1\in{}N_{\Spin}(E_*)$ such that 
$c_{x_1}|_E=c_x|_E$ (ie, $xx_1^{-1}\in{}C_{\Spin}(E)$) and 
$c_{x_1}(\BBB)=\BBB$ (ie, $[x_1,\BBB]=1$).  Set $x_2=xx_1^{-1}$.  

Since $C_{\Spin}(U)=H(q^n)\le\Im(\omega)$, we see that 
$C_{\Spin}(E)=K_0{\cdot}\gen{\BBB}$, where 
	$$ K_0 = \omega(C_{SL_2(q^\infty)}(A)^3) \cap \Spin $$
is abelian, $R_0\in\Syl_2(K_0)$, and $\BBB$ acts on $K_0$ by inversion.  
Upon replacing $x_1$ by $\BBB{}x_1$ and $x_2$ by $x_2\BBB^{-1}$ if 
necessary, we can assume that $x_2\in{}K_0$.  Then
	$$ [x_2,\BBB] = x_2{\cdot}(\BBB x_2 \BBB^{-1})^{-1} = x_2^2, $$
while by the original choice of $x,x_1$ we have
	$$ [x_2,\BBB] = [xx_1^{-1},\BBB] = [x,\BBB]\in R_0. $$
Thus $x_2^2\in{}R_0\in\Syl_2(K_0)$, and hence $x_2\in{}R_0\le{}R_1$.  Since 
$x=x_2x_1$ was an arbitrary element of $N_{\Spin}(R_1)$, this shows that
$N_{\Spin}(R_1)\le{}R_1{\cdot}C_{\Spin}(\BBB)$, and hence that
	\beq \Aut_{\Spin}(R_1) = \Inn(R_1){\cdot}
        \{\varphi\in\Aut_{\Spin}(R_1)\,|\,\varphi(\BBB)=\BBB\}. \tag{5} 
        \eeq

Since $\Aut_\Gamma(R_1)$ is generated by its intersection with 
$\Aut_{\Spin}(R_1)$ and the group $\Sthree$ which permutes the three 
factors in $H(q^\infty)$ (and since the elements of $\Sthree$ all fix 
$\BBB$), we also have
        $$ \Aut_{\Gamma}(R_1) = \Inn(R_1) {\cdot} 
        \{\varphi\in\Aut_{\Gamma}(R_1)\,|\,\varphi(\BBB)=\BBB\}. $$
Together with (4) and (5), this shows that $\Aut_\calf(R_1)$ is generated 
by $\Inn(R_1)$ together with certain automorphisms of 
$R_1=R_0{\cdot}\gen{\BBB}$ which send $\BBB$ to itself.  In other words,
        \begin{align*}  
        \Aut_{\calf}(R_1) & = \Inn(R_1) {\cdot}
        \bigl\{\varphi\in\Aut(R_1) \,\big|\, \varphi(\BBB)=\BBB,\ 
        \varphi|_{R_0} \in \Aut_{\calf}(R_0) \bigr\} \\
        & = \Inn(R_1) {\cdot}
        \bigl\{\varphi\in\Aut(R_1)\,\big|\, \varphi(\BBB)=\BBB,\ 
        \varphi|_{R_0} \in \mu(GL_3(\Z/2)) \bigr\}.
        \end{align*}
Thus
	\begin{multline*}  
	\bigl\{\varphi\in\Aut_{\calf}(R_1) \,\big|\, \varphi(z)=z \bigr\} 
	\\
	= \Inn(R_1){\cdot} 
	\bigl\{\varphi\in\Aut(R_1) \,\big|\, \varphi(\BBB)=\BBB,\ 
	\varphi|_{R_0}\in\mu(T_1)=\Aut_{\Spin}(R_0) \bigr\} \\
	= \Aut_{\Spin}(R_1), 
	\end{multline*}
the last equality by (5); and (1) now follows.

\smallskip

\noindent\textbf{Case 2: }  Now assume that $E\nleq{}R_0$.  By assumption, 
$U\le{}E$ (hence $E\le{}C_S(E)\le{}C_S(U)$), and $C_S(E)$ is a Sylow 
subgroup of $C_{\Spin}(E)$.  Since $C_S(E)$ is not isomorphic to 
$R_1=C_S(\gen{z,z_1,\AAA})$ (by (2)), this shows that $E$ is not 
$\Spin$--conjugate to $\gen{z,z_1,\AAA}$.  By Proposition \ref{elem.abel.}, 
$\Spin$ contains exactly two conjugacy classes of rank 3 subgroups 
containing $z$, and thus $E$ must have type II.  Hence by Proposition 
\ref{elem.abel.}(d), $C_S(E)$ is elementary abelian of rank 4, and also 
has type II.  

Let $\calc$ be the $\Spin_7(q^n)$--conjugacy class of the subgroup 
$E_*=\gen{U,\AAA,\BBB}\cong{}C_2^4$, which by Lemma \ref{E4-props}(a) has 
type I.  Let $\cale'$ be the set of all subgroups of $S$ which are 
elementary abelian of rank 4, contain $U$, and are not in $\calc$.  By Lemma 
\ref{E4-props}(e), for any $\varphi\in\Iso_{\Gamma}(E',E'')$ and any 
$E'\in\cale'$, $E''\defeq\varphi(E')\in\cale'$, and $\varphi$ sends 
$x_{\calc}(E')$ to $x_{\calc}(E'')$.  The same holds for 
$\varphi\in\Iso_{\Spin}(E',E'')$ by definition of the elements 
$x_{\calc}(-)$ (Proposition \ref{E4II}).  Since $C_S(E)\in\cale'$, this 
shows that all elements of $\Aut_{\calf}(C_S(E))$ send the element 
$x_{\calc}(C_S(E))$ to itself.  By Proposition \ref{E4II}(c), 
$\Aut_{\Spin}(C_S(E))$ is the group of automorphisms which are the 
identity on the rank two subgroup $\gen{x_{\calc}(C_S(E)),z}$; and (1) now 
follows.
\end{proof}

One more technical result is needed.

\newcommand{\SLx}{SL^*}

\begin{Lem} \label{point(e)}
Fix $n\ge1$, and let $E,E'\le{}S(q^n)$ be two elementary abelian subgroups 
of rank three which contain $U$, and which are $\Gamma_n$--conjugate.  Then 
$E$ and $E'$ are $\Spin_7(q^n)$--conjugate.
\end{Lem}

\begin{proof}  By \cite[3.6.3(ii)]{Suzuki}, $-I$ is the only element of 
order 2 in $SL_2(q^\infty)$.  Consider the sets
	$$ \calj_1 = \bigl\{X\in SL_2(q^n) \,\big|\, X^2=-I \bigr\} $$
and
	$$ \calj_2 = \bigl\{X\in SL_2(q^{2n}) \,\big|\, \psi^{q^n}(X)=-X,\ 
	X^2=-I \bigr\}. $$
Here, as usual, $\psi^{q^n}$ is induced by the field automorphism 
$(x\mapsto{}x^{q^n})$.  All elements in $\calj_1$ are $SL_2(q)$--conjugate 
(this follows, for example, from \cite[3.6.23]{Suzuki}), and we claim the 
same is true for elements of $\calj_2$.

Let $\SLx_2(q^n)$ be the group of all elements $X\in{}SL_2(q^{2n})$ such 
that $\psi^{q^n}(X)=\pm{}X$.  This is a group which contains $SL_2(q^n)$ 
with index 2.  Let $k$ be such that the Sylow 2--subgroups of $SL_2(q^n)$ 
have order $2^k$; then $k\ge3$ since $|SL_2(q^n)|=q^n(q^{2n}-1)$.  Any 
$S\in\Syl_2(\SLx_2(q^n))$ is quaternion of order $2^{k+1}\ge16$ 
(see \cite[Theorem 2.8.3]{Gorenstein}) and its intersection with 
$SL_2(q^n)$ is quaternion of order $2^k$, so all elements in $S\cap\calj_2$ 
are $S$--conjugate.  It follows that all elements of $\calj_2$ are 
$\SLx_2(q^n)$--conjugate.  If $X,X'\in\calj_2$ and $X'=gXg^{-1}$ for 
$g\in\SLx_2(q^n)$, then either $g\in{}SL_2(q^n)$ or $gX\in{}SL_2(q^n)$, 
and in either case $X$ and $X'$ are conjugate by an element of $SL_2(q^n)$.

By Proposition \ref{omega}(a),
	$$ E,E' \le C_{\Spin_7(q^n)}(U)=H(q^n)
	\defeq \omega(SL_2(q^\infty)^3)\cap\Spin_7(q^n). $$
Thus $E=\gen{z,z_1,\trp[X_1,X_2,X_3]}$ and 
$E'=\gen{z,z_1,\trp[X'_1,X'_2,X'_3]}$, where the $X_i$ are all in 
$\calj_1$ or all in $\calj_2$, and similarly for the $X'_i$.  Also, since 
$E$ and $E'$ are $\Gamma_n$--conjugate (and each element of $\Gamma_n$ 
leaves $U=\gen{z,z_1}$ invariant), the $X_i$ and $X'_i$ must all be in the 
same set $\calj_1$ or $\calj_2$.  Hence they are all 
$SL_2(q^n)$--conjugate, and so $E$ and $E'$ are $\Spin_7(q^n)$--conjugate.
\end{proof}


We are now ready to show that the fusion systems $\calf_n$ are saturated, 
and satisfy the conditions listed in Theorem \ref{Sol(q)}.

\begin{Prop} \label{F_Sol(fqbar)}
For a fixed odd prime power $q$, let 
$S(q^n)\le{}S(q^\infty)\le\Spin_7(q^\infty)$ be as defined above.  Let 
$z\in{}Z(\Spin_7(q^\infty))$ be the central element of order 2.  Then for 
each $n$, $\calf_n=\calf_{\Sol}(q^n)$ is saturated as a fusion system over 
$S(q^n)$, and satisfies the following conditions:
\begin{enumerate}  
\item For all $P,Q\le{}S(q^n)$ which contain $z$, if $\alpha\in\Hom(P,Q)$ is 
such that $\alpha(z)=z$, then $\alpha\in\Hom_{\calf_n}(P,Q)$ if and only if 
$\alpha\in\Hom_{\Spin_7(q^n)}(P,Q)$.
\item $C_{\calf_n}(z)=\calf_{S(q^n)}(\Spin_7(q^n))$ as fusion systems over 
$S(q^n)$.
\item All involutions of $S(q^n)$ are $\calf_n$--conjugate.
\end{enumerate}
Furthermore, $\calf_m\subseteq\calf_n$ for $m|n$.  The union of 
the $\calf_n$ is thus a category $\calf_{\Sol}(q^\infty)$ 
whose objects are the finite subgroups of $S(q^\infty)$.
\end{Prop}

\begin{proof}  We apply Proposition \ref{construct-sFs}, where $p=2$, 
$G=\Spin_7(q^n)$, $S=S(q^n)$, $Z=\gen{z}=Z(G)$; and $U$ and 
$C_G(U)=H(q^n)$ are as defined above.  Also, 
$\Gamma=\Gamma_n\le\Aut(H(q^n))$.  Condition (a) in Proposition 
\ref{construct-sFs} (all noncentral involutions in $G$ are conjugate) 
holds since all subgroups in $\cale_2$ are conjugate (Proposition 
\ref{elem.abel.}), and condition (b) holds by definition of $\Gamma$.  
Condition (c) holds since
        $$ \{\gamma\in\Gamma\,|\,\gamma(z)=z\} = 
        \Inn(H(q^n)){\cdot}\gen{c_\tau} = \Aut_{N_G(U)}(H(q^n)) $$
by definition, since $H(q^n)=C_G(U)$, and by Proposition \ref{omega}(b).  
Condition (d) was shown in Proposition \ref{rk3}, and condition (e) in 
Lemma \ref{point(e)}.  So by Proposition \ref{construct-sFs}, $\calf_n$ is 
a saturated fusion system, and 
$C_{\calf_n}(Z)=\calf_{S(q^n)}(\Spin_7(q^n))$.  

The last statement is clear.
\end{proof}


\section{Linking systems and their automorphisms}
\label{sect:Sol2}

We next show the existence and uniqueness of centric linking systems 
associated to the $\calf_{\Sol}(q)$, and also construct certain 
automorphisms of these categories analogous to the automorphisms $\psi^q$ 
of the group $\Spin_7(q^n)$.  One more technical lemma about 
elementary abelian subgroups, this time about their $\calf$--conjugacy 
classes, is first needed.

\begin{Lem}  \label{rk4}
Set $\calf=\calf_{\Sol}(q)$.  For each $r\le3$, there is a unique 
$\calf$--conjugacy class of elementary abelian subgroups $E\le{}S(q)$ of 
rank $r$.  There are two $\calf$--conjugacy classes of rank four elementary 
abelian subgroups $E\le{}S(q)$:  one is the set $\calc$ of subgroups 
$\Spin_7(q)$--conjugate to $E_*=\gen{z,z_1,\AAA,\BBB}$, while the other 
contains the other conjugacy class of type I subgroups as well as all 
type II subgroups.  Furthermore,  $\Aut_{\calf}(E)=\Aut(E)$ for all 
elementary abelian subgroups $E\le{}S(q)$ \emph{except} when $E$ has rank 
four and is not $\calf$--conjugate to $E_*$, in which case
        $$ \Aut_{\calf}(E) = 
        \{\alpha\in\Aut(E)\,|\,\alpha(x_{\calc}(E))=x_{\calc}(E)\}. $$
\end{Lem}

\begin{proof}  By Lemma \ref{E4-props}(d), the three subgroups
\begin{small}  
	$$ E_*=\gen{z,z_1,\AAA,\trp[B,B,B]},\
	E_{001}=\gen{z,z_1,\AAA,\trp[B,B,XB]},\
	E_{100}=\gen{z,z_1,\AAA,\trp[XB,B,B]} $$
\end{small}
(where $X$ is a generator of $C(q)$) represent the three 
$\Spin_7(q)$--conjugacy classes of rank four subgroups.  Clearly, $E_{100}$ 
and $E_{001}$ are $\Gamma_1$--conjugate, hence $\calf$--conjugate; and by 
Lemma \ref{E4-props}(e), neither is $\Gamma_1$--conjugate to $E_*$.  This 
proves that there are exactly two $\calf$--conjugacy classes of such 
subgroups.

Since $E_*$ and $E_{001}$ both are of type I in $\Spin_7(q)$, their 
$\Spin_7(q)$--automorphism groups contain all automorphisms which fix $z$ 
(see Proposition \ref{elem.abel.}).  By Lemma 
\ref{E4-props}(e), $z$ is fixed by all $\Gamma$--automorphisms of 
$E_{001}$, and so $\Aut_{\calf}(E_{001})$ is the group of all 
automorphisms of $E_{001}$ which send $z=x_{\calc}(E_{001})$ to itself.  
On the other hand, $E_*$ contains automorphisms (induced by permuting the 
three coordinates of $H$) which permute the three elements $z,z_1,zz_1$; 
and these together with $\Aut_{\Spin}(E_*)$ generate $\Aut(E_*)$.

It remains to deal with the subgroups of smaller rank.  By Proposition 
\ref{elem.abel.} again, there is just one $\Spin_7(q)$--conjugacy class of 
elementary abelian subgroups of rank one or two.  There are two conjugacy 
classes of rank three subgroups, those of type I and those of type II.  
Since $E_{100}$ is of type II and $E_{001}$ of type I, all rank three 
subgroups of $E_{001}$ have type I, while some of the rank three subgroups 
of $E_{100}$ have type II.  Since $E_{001}$ is $\calf$--conjugate to 
$E_{100}$, this shows that some subgroup of rank three and type II is 
$\calf$--conjugate to a subgroup of type I, and hence all rank three 
subgroups are conjugate to each other.  Finally, 
$\Aut_{\calf}(E)=\Aut(E)$ whenver $\rk(E)\le3$ since any such group is 
$\calf$--conjugate to a subgroup of $E_*$ (and we have just seen that 
$\autf(E_*)=\Aut(E_*)$).
\end{proof}

To simplify the notation, we now define
	$$ \fspin(q^n) \defeq \calf_{S(q^n)}(\Spin_7(q^n)) $$
for all $1\le{}n\le\infty$:  the fusion system of the group $\Spin_7(q^n)$ 
at the Sylow subgroup $S(q^n)$.  By construction, this is a subcategory of 
$\fsol(q^n)$.  We write 
	$$ \osol(q^n)=\orb(\fsol(q^n)) \qquad\textup{and}\qquad
	\ospin(q^n)=\orb(\fspin(q^n)) $$
for the corresponding orbit categories:  both of these have as objects the 
subgroups of $S(q^n)$, and have as morphism sets
	$$ \Mor_{\osol(q^n)}(P,Q)=\Hom_{\fsol(q^n)}(P,Q)/\Inn(Q) 
	\subseteq \Rep(P,Q) $$
and
	$$ \Mor_{\ospin(q^n)}(P,Q)=\Hom_{\fspin(q^n)}(P,Q)/\Inn(Q)\,. $$
Let $\ocsol(q^n)\subseteq\osol(q^n)$ and $\ocspin(q^n)\subseteq\ospin(q^n)$ 
be the centric orbit categories; ie, the full subcategories whose objects 
are the $\fsol(q^n)$-- or $\fspin(q^n)$--centric subgroups of $S(q^n)$.  
(We will see shortly that these in fact have the same objects.)

The obstructions to the existence and uniqueness of linking systems 
associated to the fusion systems $\fsol(q^n)$, and to the existence and 
uniqueness of certain automorphisms of those linking systems, lie in 
certain groups which were identified in \cite{BLO2} and \cite{BLO1}.  It is 
these groups which are shown to vanish in the next lemma.

\begin{Lem} \label{lim*Z=0}
Fix a prime power $q$, and let 
	$$ \zsol(q)\:\ocsol(q)\Right3{}\Ab
	\qquad\textup{and}\qquad
	\zspin(q)\:\ocspin(q)\Right3{}\Ab $$
be the functors $\calz(P)=Z(P)$.  Then for all $i\ge0$,
	$$ \higherlim{\ocsol(q)}i(\zsol(q)) = 0
	= \higherlim{\ocspin(q)}i(\zspin(q)). $$
\end{Lem}

\begin{proof}  Set $\calf=\fsol(q)$ for short.  Let $P_1,\dots,P_k$ be 
$\calf$--conjugacy class representatives for all $\calf$--centric subgroups 
$P_i\le{}S(q)$, arranged such that $|P_i|\le|P_j|$ for $i\le{}j$.  For each 
$i$, let $\calz_i\subseteq\zsol(q)$ be the subfunctor defined by setting 
$\calz_i(P)=\zsol(q)(P)$ if $P$ is conjugate to $P_j$ for some $j\le{}i$ 
and $\calz_i(P)=0$ otherwise.  We thus have a filtration
	$$ 0 = \calz_0 \subseteq \calz_1 \subseteq \cdots \subseteq 
	\calz_k=\zsol(q) $$
of $\zsol(q)$ by subfunctors, with the property that for each $i$, the 
quotient functor $\calz_i/\calz_{i-1}$ vanishes except on the conjugacy 
class of $P_i$ (and such that $(\calz_i/\calz_{i-1})(P_i)=\zsol(q)(P_i)$).  
By \cite[Proposition 3.2]{BLO2}, 
	$$ \higherlim{}*(\calz_i/\calz_{i-1}) \cong 
	\Lambda^*(\outf(P_i);Z(P_i)) $$
for each $i$.  Here, $\Lambda^*(\Gamma;M)$ are certain graded groups, 
defined in \cite[\S5]{JMO} for all finite groups $\Gamma$ and all
finite $\zploc[\Gamma]$--modules $M$.  We will show that 
$\Lambda^*(\outf(P_i);Z(P_i))=0$ except when $P_i=S(q)$ or $S_0(q)$ 
(see Definition \ref{S-def}).

Fix an $\calf$--centric subgroup $P\le{}S(q)$.  For each $j\ge1$, let 
$\Omega_j(Z(P))=\{g\in{}Z(P)\,|\,g^{2^j}=1\}$, and set $E=\Omega_1(Z(P))$ 
--- the 2--torsion in the center of $P$.  For each $j\ge1$, let 
$\Omega_j(Z(P))=\{g\in{}Z(P)\,|\,g^{2^j}=1\}$, and set $E=\Omega_1(Z(P))$ 
--- the 2--torsion in the center of $P$.  We can assume $E$ is fully 
centralized in $\calf$ (otherwise replace $P$ and $E$ by appropriate 
subgroups in the same $\calf$--conjugacy classes).  

Assume first that $Q\defeq{}C_{S(q)}(E)\gneqq{}P$, and hence that 
$N_Q(P)\gneqq{}P$.  Then any $x\in{}N_Q(P){\sminus}P$ centralizes 
$E=\Omega_1(Z(P))$.  Hence for each $j$, $x$ acts trivially on 
$\Omega_j(Z(P))/\Omega_{j-1}(Z(P))$, since multiplication by $p^{j-1}$ 
sends this group $N_Q(P)/P$--linearly and monomorphically to $E$.  Since 
$c_x$ is a nontrivial element of $\outf(P)$ of $p$--power order, 
	$$ \Lambda^*(\outf(P);\Omega_j(Z(P))/\Omega_{j-1}(Z(P)))=0 $$
for all $j\ge1$ by \cite[Proposition 5.5]{JMO}, and thus
$\Lambda^*(\outf(P);Z(P))=0$.  

Now assume that $P=C_{S(q)}(E)=P$, the centralizer in $S(q)$ of a fully 
$\calf$--centralized elementary abelian subgroup.  Since there is a unique 
conjugacy class of elementary abelian subgroup of any rank $\le3$, 
$C_{S(q)}(E)$ always contains a subgroup $C_2^4$, and hence $P$ contains a 
subgroup $C_2^4$ which is self centralizing by Proposition 
\ref{elem.abel.}(a).  This shows that $Z(P)$ is elementary abelian, and 
hence that $Z(P)=E$.  

We can assume $P$ is fully normalized in $\calf$, so
	$$ \Aut_{S(q)}(P) \in \Syl_2(\Aut_\calf(P)) $$
by condition (I) in the definition of a saturated fusion system.  Since 
$P=C_{S(q)}(E)$ (and $E=Z(P)$), this shows that 
	$$ \Ker\bigl[\Out_{\calf}(P) \Right2{} \Aut_\calf(E)\bigr] $$
has odd order.  Also, since $E$ is fully centralized, any 
$\calf$--automorphism of $E$ extends to an $\calf$--automorphism of 
$P=C_{S(q)}(E)$, and thus this restriction map between automorphism groups 
is onto.  By \cite[Proposition 6.1(i,iii)]{JMO}, it now follows that
	\beq \Lambda^i(\Out_\calf(P);Z(P)) \cong 
	\Lambda^i(\Aut_\calf(E);E). \tag{1} \eeq

By Lemma \ref{rk4}, $\Aut_\calf(E)=\Aut(E)$, except when $E$ lies in one 
certain $\calf$--conjugacy class of subgroups
$E\cong{}C_2^4$; and in this case $P=E$ and $\Aut_\calf(E)$ is the group of 
automorphisms fixing the element $x_\calc(E)$.  
In this last (exceptional) case, $O_2(\Aut_\calf(E))\ne1$ (the subgroup of 
elements which are the identity on $E/\gen{x_\calc(E)}$), so
	\beq \Lambda^*(\Out_\calf(P);Z(P)) = \Lambda^*(\Aut_\calf(E);E) = 0 
	\tag{2} \eeq
by \cite[Proposition 6.1(ii)]{JMO}.  Otherwise, when 
$\Aut_\calf(E)=\Aut(E)$, by \cite[Proposition 6.3]{JMO} we have
	\beq \Lambda^i(\Aut_\calf(E);E) \cong
	\begin{cases}  
	\Z/2 & \textup{if $\rk(E)=2$, $i=1$} \\
	\Z/2 & \textup{if $\rk(E)=1$, $i=0$} \\
	0 & \textup{otherwise.}
	\end{cases} \tag{3} \eeq
By points (1), (2), and (3), the groups $\Lambda^*(\Out_\calf(P);Z(P))$ 
vanish except in the two cases $E=\gen{z}$ or $E=U$, and these correspond 
to $P=S(q)$ or $P=N_{S(q)}(U)=S_0(q)$.

We can assume that $P_k=S(q)$ and $P_{k-1}=S_0(q)$.  We have now shown that 
$\higherlim{}*(\calz_{k-2})=0$, and thus that $\zsol(q)$ has the same higher 
limits as $\calz_k/\calz_{k-2}$.  Hence $\higherlim{}j(\zsol(q))=0$ for all 
$j\ge2$, and there is an exact sequence
	\begin{multline*}  
	0 \Right2{} \higherlim{}0(\zsol(q)) \Right3{}
	\underset{\cong\Z/2}{\higherlim{}0(\calz_k/\calz_{k-1})} \Right3{}
	\underset{\cong\Z/2}{\higherlim{}1(\calz_{k-1}/\calz_{k-2})} \\
	\Right3{} \higherlim{}1(\zsol(q)) \Right2{} 0. 
	\end{multline*}
One easily checks that $\higherlim{}0(\zsol(q))=0$, and hence we also 
get $\higherlim{}1(\zsol(q))=0$.

The proof that $\higherlim{}i(\zspin(q))=0$ for all $i\ge1$ is similar, but 
simpler.  If $\calf=\fspin(q)$, then for any $\calf$--centric subgroup 
$P\lneqq{}S(q)$, there is an element $x\in{}N_S(P){\sminus}P$ such that 
$[x,P]=\gen{z}$, and $c_x$ is a nontrivial element of $O_2(\outf(P))$.  
Thus 
	$$ \Lambda^*(\outf(P);Z(P))=0 $$ 
for all such $P$ by \cite[Proposition 6.1(ii)]{JMO} again.
\end{proof}

We are now ready to construct classifying spaces $B\Sol(q)$ for these 
fusion systems $\fsol(q)$.  The following proposition finishes the proof 
of Theorem \ref{Sol(q)}, and also contains additional information about 
the spaces $B\Sol(q)$.  

To simplify notation, we write 
$\lcspin(q^n)=\call^c_{S(q^n)}(\Spin_7(q^n))$ ($n\ge1$) to denote the 
centric linking system for the group $\Spin_7(q^n)$.  The field 
automorphism $(x\mapsto{}x^q)$ induces an automorphism of $\Spin_7(q^n)$ 
which sends $S(q^n)$ to itself; and this in turn induces automorphisms 
$\psi^q_\calf=\psi^q_\calf(\Sol)$, $\psi^q_\calf(\Spin)$, and 
$\psi^q_\call(\Spin)$ of the fusion systems $\fsol(q^n)\supseteq\fspin(q^n)$ 
and of the linking system $\lcspin(q^n)$.  

\begin{Prop}  \label{Sol(q)+}
Fix an odd prime $q$, and $n\ge1$.  Let $S=S(q^n)\in\Syl_2(\Spin_7(q^n))$ 
be as defined above.  Let $z\in{}Z(\Spin_7(q^n))$ be the central element 
of order 2.  Then there is a centric linking system 
	$$ \call=\lcsol(q^n) \Right4{\pi} \fsol(q^n) $$
associated to the saturated fusion system $\calf\defeq\fsol(q^n)$ over 
$S$, which has the following additional properties.  
\begin{enumerate}
\item A subgroup $P\le{}S$ is $\calf$--centric if and only if it is 
$\fspin(q^n)$--centric.

\item $\lcsol(q^n)$ contains $\lcspin(q^n)$ as a subcategory, in such a way 
that $\pi|_{\lcspin(q^n)}$ is the usual projection to $\fcspin(q^n)$, 
and that the distinguished monomorphisms 
	$$ P\Right2{\delta_P}\Aut_\call(P) $$
for $\call=\lcsol(q^n)$ are the same as those for $\lcspin(q^n)$.

\item Each automorphism of $\lcspin(q^n)$ which covers the identity on 
$\fcspin(q^n)$ extends to an automorphism of $\lcsol(q^n)$ which covers 
the identity on $\fcsol(q^n)$.  Furthermore, such an extension is unique up to 
composition with the functor
	$$ C_z\: \lcsol(q^n) \Right4{} \lcsol(q^n) $$
which is the identity on objects and sends 
$\alpha\in\Mor_{\lcsol(q^n)}(P,Q)$ to $\q{z}\circ\alpha\circ\q{z}^{-1}$ 
(``conjugation by $z$'').

\item There is a unique automorphism 
$\psi^q_\call\in\Aut(\lcsol(q^n))$ which covers the 
automorphism of $\fsol(q^n)$ induced by the field automorphism 
$(x\mapsto{}x^q)$, which extends the automorphism of $\lcspin(q^n)$ 
induced by the field automorphism, and which is the identity on 
$\pi^{-1}(\fsol(q))$.
\end{enumerate}
\end{Prop}

\begin{proof}  By Proposition \ref{F_Sol(fqbar)}, $\calf=\fsol(q^n)$ is a 
saturated fusion system over $S=S(q^n)\in\Syl_2(\Spin_7(q^n))$, with the 
property that $C_{\calf}(z)=\fspin(q^n)$.  Point (a) follows as a special 
case of \cite[Proposition 2.5(a)]{BLO2}.

Since $\hilim{\ocsol(q^n)}i{\zsol(q^n)}=0$ for $i=2,3$ by Lemma 
\ref{lim*Z=0}, there is by \cite[Proposition 3.1]{BLO2} a centric linking 
system $\call=\lcsol(q^n)$ associated to $\calf$, which is unique up to 
isomorphism (an isomorphism which commutes with the projection to 
$\fsol(q^n)$ and with the distinguished monomorphisms).  Furthermore, 
$\pi^{-1}(\fspin(q^n))$ is a linking system associated to $\fspin(q^n)$, 
such a linking system is unique up to isomorphism since 
$\higherlim{}2(\zspin(q^n))=0$ (Lemma \ref{lim*Z=0} again), and this 
proves (b).

\smallskip

\noindent\textbf{(c) } By \cite[Theorem 6.2]{BLO1} (more precisely, by the 
same proof as that used in \cite{BLO1}), the vanishing of 
$\higherlim{}i(\zsol(q^n))$ for $i=1,2$ (Lemma \ref{lim*Z=0}) shows that 
each automorphism of $\calf=\fsol(q^n)$ lifts to an automorphism of 
$\call$, which is unique up to a natural isomorphism of functors; and any 
such natural isomorphism sends each object $P\le{}S$ to a isomorphism 
$\q{g}$ for some $g\in{}Z(P)$.  Similarly, the vanishing of 
$\higherlim{}i(\zspin(q^n))$ for $i=1,2$ shows that each automorphism of 
$\fspin(q^n)$ lifts to an automorphism of $\lcspin(q^n)$, also unique up 
to a natural isomorphism of functors.  Since $\lcsol(q^n)$ and 
$\lcspin(q^n)$ have the same objects by (a), this shows that each 
automorphism of $\lcspin(q^n)$ which covers the identity on $\fcspin(q^n)$ 
extends to a unique automorphism of $\lcsol(q^n)$ which covers the 
identity on $\fsol(q^n)$.  

It remains to show, for any $\Phi\in\Aut(\lcsol(q^n))$ which covers the 
identity on $\fcsol(q^n)$ and such that $\Phi|_{\lcspin(q^n)}=\Id$, that 
$\Phi$ is the identity or conjugation by $z$.  We have already noted that 
$\Phi$ must be naturally isomorphic to the identity; ie, that there are 
elements $\gamma(P)\in{}Z(P)$, for all $P$ in $\lcsol(q^n)$, such that
	$$ \Phi(\alpha)=\gamma(Q)\circ\alpha\circ\gamma(P)^{-1}
	\qquad\textup{for all $\alpha\in\Mor_{\lcsol(q^n)}(P,Q)$, all 
	$P,Q$.} $$
Since $\Phi$ is the identity on $\lcspin(q^n)$, the only possibilities are 
$\gamma(P)=1$ for all $P$ (hence $\Phi=\Id$), or $\gamma(P)=z$ for all $P$ 
(hence $\Phi$ is conjugation by $z$).

\smallskip

\noindent\textbf{(d) } Now consider the automorphism 
$\psi^q_\calf\in\Aut(\fsol(q^n))$ induced by the field automorphism 
$(x\mapsto{}x^q)$ of $\F_{q^n}$.  We have just seen that this lifts to an 
automorphism $\psi^q_\call$ of $\lcsol(q^n)$, which is unique up to 
natural isomorphism of functors.  The restriction of $\psi^q_\call$ to 
$\lcspin(q^n)$, and the automorphism $\psi^q_{\call}(\Spin)$ of 
$\lcspin(q^n)$ induced directly by the field automorphism, are two 
liftings of $\psi^q_\calf|_{\fspin(q^n)}$, and hence differ by a natural 
isomorphism of functors which extends to a natural isomorphism of functors 
on $\lcsol(q^n)$.  Upon composing with this natural isomorphism, we can 
thus assume that $\psi^q_\call$ does restrict to the automorphism of 
$\lcspin(q^n)$ induced by the field automorphism.  

Now consider the action of $\psi^q_\call$ on $\Aut_{\call}(S_0(q))$, which 
by assumption is the identity on $\Aut_{\lcspin(q)}(S_0(q))$, and in 
particular on $\delta(S_0(q))$ itself.  Thus, with respect to the extension
	$$ 1 \Right2{} S_0(q) \Right3{} \Aut_{\call}(S_0(q)) \Right3{} 
	\Sigma_3 \Right2{} 1, $$
$\psi^q_\call$ is the identity on the kernel and on the quotient, and 
hence is described by a cocycle 
	$$ \eta \in Z^1(\Sigma_3;Z(S_0(q))) \cong Z^1(\Sigma_3;(\Z/2)^2). $$
Since $H^1(\Sigma_3;(\Z/2)^2)=0$, $\eta$ must be a coboundary, and thus 
the action of $\psi^q_\call$ on $\Aut_{\call}(S_0(q))$ is conjugation by 
an element of $Z(S_0(q))$.  Since it is the identity on 
$\Aut_{\lcspin(q)}(S_0(q))$, it must be conjugation by $1$ or $z$.  If it 
is conjugation by $z$, then we can replace $\psi^q_\call$ (on the whole 
category $\call$) by its composite with $z$; ie, by its composite with 
the functor which is the identity on objects and sends 
$\alpha\in\Mor_\call(P,Q)$ to $\q{z}\circ\alpha\circ\q{z}$.  

In this way, we can assume that $\psi^q_\call$ is the identity on 
$\Aut_{\call}(S_0(q))$.  By construction, every morphism in $\fsol(q)$ is 
a composite of morphisms in $\fspin(q)$ and restrictions of automorphisms 
in $\fsol(q)$ of $S_0(q)$.  Since $\psi^q_\call$ is the identity on 
$\pi^{-1}(\fspin(q))$, this shows that it is the identity on 
$\pi^{-1}(\fsol(q))$.

It remains to check the uniqueness of $\psi^q_\call$.  If $\psi'$ is 
another functor with the same properties, then by (e), 
$(\psi')^{-1}\circ\psi^q_\call$ is either the identity or conjugation by 
$z$; and the latter is not possible since conjugation by $z$ is not the 
identity on $\pi^{-1}(\fsol(q))$.
\end{proof}

This finishes the construction of the classifying spaces 
$B\Sol(q)=|\call_{\Sol}(q)|\doscom$ for the fusion systems constructed in 
Section \ref{sect:Sol}.  We end the section with an explanation of why 
these are not the fusion systems of finite groups.  

\begin{Prop} \label{G<>Sol} For any odd prime power $q$, there is no 
finite group $G$ whose fusion system is isomorphic to that of 
$\calf_{\Sol}(q)$. 
\end{Prop}

\begin{proof}  Let $G$ be a finite group, fix $S\in\Syl_2(G)$, and assume 
that $S\cong{}S(q)\in\Syl_2(\Spin_7(q))$, and that the fusion system 
$\calf_{S}(G)$ satisfies conditions (a) and (b) in Theorem \ref{Sol(q)}.  
In particular, all involutions in $G$ are conjugate, and the centralizer 
of any involution $z\in{}G$ has the fusion system of $\Spin_7(q)$.  When 
$q\equiv\pm3$ (mod $8$), Solomon showed \cite[Theorem 3.2]{Solomon} that 
there is no finite group whose fusion system has these properties.  When 
$q\equiv\pm1$ (mod $8$), he showed (in the same theorem) that there is no 
such $G$ such that $\widehat{H}\defeq{}C_G(z)/O_{2'}(C_G(z))$ is 
isomorphic to a subgroup of $\Aut(\Spin_7(q))$ which contains $\Spin_7(q)$ 
with odd index. (Here, $O_{2'}(-)$ means largest odd order normal 
subgroup.)  

Let $G$ be a finite group whose fusion system is isomorphic to $\fsol(q)$, 
and again set $\widehat{H}\defeq{}C_G(z)/O_{2'}(C_G(z))$ for some 
involution $z\in{}G$.  Set $H=O^{2'}(\widehat{H}/\gen{z})$:  the smallest 
normal subgroup of $\widehat{H}/\gen{z}$ of odd index.  Then $H$ has the 
fusion system of $\Omega_7(q)\cong\Spin_7(q)/Z(\Spin_7(q))$.  We will show 
that $H\cong\Omega_7(q')$ for some odd prime power $q'$.  It then follows 
that $O^{2'}(\widehat{H})\cong\Spin_7(q')$, thus contradicting Solomon's 
theorem and proving our claim.

The following ``classification free'' argument for proving that 
$H\cong\Omega_7(q')$ for some $q'$ was explained to us by Solomon.  We 
refer to the appendix for general results about the groups 
$\Spin_n^\pm(q)$ and $\Omega_n^\pm(q)$.  Fix $S\in\Syl_2(H)$.  Thus $S$ is 
isomorphic to a Sylow 2--subgroup of $\Omega_7(q)$, and has the same 
fusion.  

We first claim that $H$ must be simple.  By definition 
($H=O^{2'}(\widehat{H}/\gen{z})$), $H$ has no proper normal subgroup of 
odd index, and $H$ has no proper normal subgroup of odd order since any 
such subgroup would lift to an odd order normal subgroup of 
$\widehat{H}=C_G(z)/O_{2'}(C_G(z))$.  Hence for any proper normal subgroup 
$N\nsg{}H$, $Q\defeq{}N\cap{}S$ is a proper normal subgroup of $S$, which 
is strongly closed in $S$ with respect to $H$ in the sense that no element 
of $Q$ can be $H$--conjugate to an element of $S{\sminus}Q$.  Using Lemma 
\ref{SO-invol}(a), one checks that the group $\Omega_7(q)$ contains three 
conjugacy classes of involutions, classified by the dimension of their 
$(-1)$--eigenspace.  It is not hard to see (by taking products) that any 
subgroup of $S$ which contains all involutions in one of these conjugacy 
classes contains all involutions in the other two classes as well.  
Furthermore, $S$ is generated by the set of all of its involutions, and 
this shows that there are no proper subgroups which are strongly closed in 
$S$ with respect to $H$.  Since we have already seen that the intersection 
with $S$ of any proper normal subgroup of $H$ would have to be such a 
subgroup, this shows that $H$ is simple.

Fix an isomorphism 
	$$ S \RIGHT6{\varphi}{\cong} S' \in \Syl_2(\Omega_7(q)) $$ which 
preserves fusion.  Choose $x'\in{}S'$ whose $(-1)$--eigenspace is 
4--dimen\-sion\-al, and such that $\gen{x'}$ is fully centralized in 
$\calf_{S'}(\Omega_7(q))$.  Then 
	$$ C_{O_7(q)}(x') \cong O_4^+(q) \times O_3(q) $$ 
by Lemma \ref{SO-invol}(c).  Since $\Omega_4^+(q)\le{}O_4^+(q)$ and 
$\Omega_3(q)\le{}O_3(q)$ both have index $4$, $C_{\Omega_7(q)}(x')$ is 
isomorphic to a subgroup of $O_4^+(q)\times{}O_3(q)$ of index 4, and 
contains a normal subgroup $K'\cong\Omega_4^+(q)\times\Omega_3(q)$ of 
index 4.  Since $\gen{x'}$ is fully centralized, $C_{S'}(x')$ is a Sylow 
2--subgroup of $C_{\Omega_7(q)}(x')$, and hence $S'_0\defeq{}S'\cap{}K'$ is 
a Sylow 2--subgroup of $K'$.  

Set $x=\varphi^{-1}(x')\in{}S$.  Since $S\cong{}S'$ have the same fusion 
in $H$ and $\Omega_7(q)$, $C_S(x)\cong{}C_{S'}(x')$ have the same fusion 
in $C_H(x)$ and $C_{\Omega_7(q)}(x')$.  Hence 
	$$ H_1(C_H(x);\Z_{(2)})\cong{}H_1(C_{\Omega_7(q)}(x');\Z_{(2)}) $$ 
(homology is determined by fusion), both have order $4$, and thus $C_H(x)$ 
also has a unique normal subgroup $K\nsg{}H$ of index $4$.  Set 
$S_0=K\cap{}S$.  Thus $\varphi(S_0)=S'_0$, and using Alperin's fusion 
theorem one can show that this isomorphism is fusion preserving with 
respect to the inclusions of Sylow subgroups $S_0\le{}K$ and 
$S'_0\le{}K'$.  

Using the isomorphisms of Proposition \ref{SO3-4}: 
	$$ \Omega_4^+(q) \cong SL_2(q)\times_{\gen{x}}SL_2(q) 
	\qquad\textup{and}\qquad \Omega_3(q) \cong  PSL_2(q), $$ 
we can write $K'=K'_1\times_{\gen{x'}}K'_2$, where $K'_1\cong{}SL_2(q)$ 
and $K'_2\cong{}SL_2(q)\times{}PSL_2(q)$.  Set 
$S'_i=S'\cap{}K'_i\in\Syl_2(K'_i)$; thus $S'_0=S'_1\times_{\gen{x'}}S'_2$. 
Set $S_i=\varphi^{-1}(S'_i)$, so that $S_0=S_1\times_{\gen{x}}S_2$ is 
normal of index $4$ in $C_S(x)$.  The fusion system of $K$ thus splits as a 
central product of fusion systems, one of which is isomorphic to the fusion 
system of $SL_2(q)$.  

We now apply a theorem of Goldschmidt, which says very roughly that under 
these conditions, the group $K$ also splits as a central product.  To make 
this more precise, let $K_i$ be the normal closure of $S_i$ in 
$K\nsg{}C_H(x)$.  By \cite[Corollary A2]{Goldschmidt}, since $S_1$ and 
$S_2$ are strongly closed in $S_0$ with respect to $K$, 
	$$ [K_1,K_2] \le \gen{x}{\cdot}O_{2'}(K). $$ 
Using this, it is not hard to check that $S_i\in\Syl_2(K_i)$.  Thus $K_1$ 
has same fusion as $SL_2(q)$ and is subnormal in $C_H(x)$ 
($K_1\nsg{}K\nsg{}C_H(x)$), and an argument similar to that used above to 
prove the simplicity of $H$ shows that $K_1/(\gen{x}{\cdot}O_{2'}(K_1))$ 
is simple.  Hence $K_1$ is a 2--component of $C_H(x)$ in the sense described 
by Aschbacher in \cite{Asch2}.  By \cite[Corollary III]{Asch2}, 
this implies that $H$ must be isomorphic to a Chevalley group of odd 
characteristic, or to $M_{11}$.  It is now straightforward to check that 
among these groups, the only possibility is that $H\cong\Omega_7(q')$ for 
some odd prime power $q'$. 
\end{proof}


\section{Relation with the Dwyer-Wilkerson space}

We now want to examine the relation between the spaces $B\Sol(q)$ which we 
have just constructed, and the space $BDI(4)$ constructed by Dwyer and 
Wilkerson in \cite{DW:DI4}.  Recall that this is a 2--complete space 
characterized by the property that its cohomology is the Dickson algebra 
in four variables over $\F_2$; ie, the ring of invariants 
$\F_2[x_1,x_2,x_3,x_4]^{GL_4(2)}$.  We show, for any odd prime power $q$, 
that $BDI(4)$ is homotopy equivalent to the 2--completion of the union of 
the spaces $B\Sol(q^n)$, and that $B\Sol(q)$ is homotopy equivalent to the 
homotopy fixed point set of an Adams map from $BDI(4)$ to itself.  

We would like to define an infinite ``linking system'' $\lcsol(q^\infty)$ 
as the union of the finite categories $\lcsol(q^n)$, and then set 
$B\Sol(q^\infty)=|\lcsol(q^\infty)|\doscom$.  The difficulty with this 
approach is that a subgroup which is centric in the fusion system 
$\fsol(q^m)$ need not be centric in a larger fusion system $\fsol(q^n)$ 
(for $m|n$).  To get around this problem, we define 
$\lccsol(q^n)\subseteq\lcsol(q^n)$ to be the full subcategory whose 
objects are those subgroups of $S(q^n)$ which are 
$\fsol(q^\infty)$--centric; or equivalently $\fsol(q^k)$--centric for
all $k\in{}n\Z$.  Similarly, we define $\lccspin(q^n)$ to be the full 
subcategory of $\lcspin(q^n)$ whose objects are those subgroups of 
$S(q^n)$ which are $\fspin(q^\infty)$--centric.  We can then define 
$\lcsol(q^\infty)$ and $\lcspin(q^\infty)$ to be the unions of these 
categories.  

For these definitions to be useful, we must first show that 
$|\lccsol(q^n)|\doscom$ has the same homotopy type as 
$|\lcsol(q^n)|\doscom$.  This is done in the following lemma.

\begin{Lem} \label{Lcc=Lc}
For any odd prime power $q$ and any $n\ge1$, the inclusions
	$$ |\lccsol(q^n)|\doscom \subseteq |\lcsol(q^n)|\doscom 
	\qquad\textup{and}\qquad
	|\lccspin(q^n)|\doscom \subseteq |\lcspin(q^n)|\doscom $$
are homotopy equivalences.
\end{Lem}

\begin{proof}  It clearly suffices to show this when $n=1$.  

Recall, for a fusion system $\calf$ over a $p$--group $S$, that a subgroup 
$P\le{}S$ is \emph{$\calf$--radical} if $\outf(P)$ is $p$--reduced; ie, if 
$O_p(\outf(P))=1$.  We will show that
	\begin{small}  
	\beq \textup{all $\fsol(q)$--centric $\fsol(q)$--radical subgroups of 
	$S(q)$ are $\fsol(q^\infty)$--centric} \tag{1} \eeq
	\end{small}%
and similarly
	\begin{small}  
	\beq \textup{all $\fspin(q)$--centric $\fspin(q)$--radical subgroups of 
	$S(q)$ are $\fspin(q^\infty)$--centric.} \tag{2} \eeq
	\end{small}%
In other words, (1) says that for each $P\le{}S(q)$ which is an object of 
$\lcsol(q)$ but not of $\lccsol(q)$, $O_2\bigl(\Out_{\fsol(q)}(P)\bigr)\ne1$.  
By \cite[Proposition 6.1(ii)]{JMO}, this implies that 
	$$ \Lambda^*(\Out_{\fsol(q)}(P);H^*(BP;\F_2))=0. $$
Hence by \cite[Propositions 3.2 and 2.2]{BLO2} (and the spectral sequence 
for a homotopy colimit), the inclusion $\lccsol(q)\subseteq\lcsol(q)$ 
induces an isomorphism 
	$$ H^*\bigl(|\lcsol(q)|;\F_2\bigr) \Right5{\cong}
	H^*\bigl(|\lccsol(q)|;\F_2\bigr) , $$
and thus $|\lccsol(q)|\doscom\simeq|\lcsol(q)|\doscom$.  The proof that 
$|\lccspin(q)|\doscom\simeq|\lcspin(q)|\doscom$ is similar, using (2).

Point (2) is shown in Proposition \ref{Spin-cent-rad}, so it remains only 
to prove (1).  Set $\calf=\fsol(q)$, and set $\calf_k=\fsol(q^k)$ for all 
$1\le{}k\le\infty$.  Let $E\le{}Z(P)$ be the 2--torsion in the center of 
$P$, so that $P\le{}C_{S(q)}(E)$.  Set 
	$$ E' = \begin{cases}  
	\gen{z} & \textup{if $\rk(E)=1$} \\
	\gen{z,z_1} & \textup{if $\rk(E)=2$} \\
	\gen{z,z_1,\AAA} & \textup{if $\rk(E)=3$} \\
	E & \textup{if $\rk(E)=4$}
	\end{cases} $$
in the notation of Definition \ref{S-def}.  In all cases, $E$ is 
$\calf$--conjugate to $E'$ by Lemma \ref{rk4}.  We claim that $E'$ is fully 
centralized in $\calf_k$ for all $k<\infty$.  This is clear when 
$\rk(E')=1$ ($E'=Z(S(q^k))$), follows from Proposition \ref{omega}(a) when 
$\rk(E')=2$, and from Proposition \ref{elem.abel.}(a) (all rank 4 
subgroups are self centralizing) when $\rk(E')=4$.  If $\rk(E')=3$, then 
by Proposition \ref{elem.abel.}(d), the centralizer in $\Spin_7(q^k)$ 
(hence in $S(q^k)$) of any rank $3$ subgroup has an abelian subgroup of 
index $2$; and using this (together with the construction of $S(q^k)$ in 
Definition \ref{S-def}), one sees that $E'$ is fully centralized in 
$\calf_k$.

If $E'\ne{}E$, choose $\varphi\in\homf(E,S(q))$ such that $\varphi(E)=E'$; 
then $\varphi$ extends to $\widebar{\varphi}\in\homf(C_{S(q)}(E),S(q))$ by 
condition (II) in the definition of a saturated fusion system, and we can 
replace $P$ by $\widebar{\varphi}(P)$ and $E$ by $\varphi(E)$.  We can 
thus assume that $E$ is fully centralized in $\calf_k$ for each 
$k<\infty$.  So by \cite[Proposition 2.5(a)]{BLO2}, $P$ is 
$\calf_k$--centric if and only if it is $C_{\calf_k}(E)$--centric; and this 
also holds when $k=\infty$.  Furthermore, since 
$\Out_{C_\calf(E)}(P)\nsg\outf(P)$, $O_2(\Out_{C_\calf(E)}(P))$ is a 
normal $2$--subgroup of $\outf(P)$, and thus
	$$ O_2\bigl(\Out_{C_\calf(E)}(P)\bigr) \le O_2(\outf(P)). $$
Hence $P$ is $C_\calf(E)$--radical if it is $\calf$--radical.  So it 
remains to show that 
	\beq \textup{all $C_\calf(E)$--centric $C_\calf(E)$--radical subgroups 
	of $S(q)$ are also $C_{\calf_\infty}(E)$--centric.} \tag{3} \eeq

If $\rk(E)=1$, then $C_\calf(E)=\fspin(q)$ and $C_{\calf_\infty}(E)= 
\fspin(q^\infty)$, and (3) follows from (2).  If $\rk(E)=4$, then 
$P=E=C_{S(q^\infty)}(E)$ by Proposition \ref{elem.abel.}(a), so $P$ is 
$\calf_\infty$--centric, and the result is clear.

If $\rk(E)=3$, then by Proposition \ref{elem.abel.}(d), 
$C_\calf(E)\subseteq{}C_{\calf_\infty}(E)$ are the fusion systems of a 
pair of semidirect products $A{\rtimes}C_2\le{}A_\infty{\rtimes}C_2$, 
where $A\le{}A_\infty$ are abelian and $C_2$ acts on $A_\infty$ by 
inversion.  Also, 
$E$ is the full 2--torsion subgroup of $A_\infty$, since otherwise 
$\rk(A_\infty)>3$ would imply $A_\infty{\rtimes}C_2\le\Spin_7(q^\infty)$ 
contains a subgroup $C_2^5$ (contradicting Proposition \ref{elem.abel.}).  
If $P\le{}A$, then either 
$\Out_{C_\calf(E)}(P)$ has order 2, which contradicts the assumption that 
$P$ is radical; or $P$ is elementary abelian and $\Out_{C_\calf(E)}(P)=1$, 
in which case $P\le{}Z(A{\rtimes}C_2)$ is not centric.  Thus $P\nleq{}A$; 
$P\cap{}A\ge{}E$ contains all 2--torsion in $A_\infty$, and hence $P$ is 
centric in $A_\infty{\rtimes}C_2$.  

If $\rk(E)=2$, then by Proposition \ref{omega}(a), $C_{\calf_\infty}(E)$ 
and $C_{\calf}(E)$ are the fusion systems of the groups 
	\beq H(q^\infty)\cong{}SL_2(q^\infty)^3/\{\pm(I,I,I)\} \tag{4} \eeq
and
	$$ H(q)=H(q^\infty)\cap\Spin_7(q)
	\ge H_0(q) \defeq SL_2(q)^3/\{\pm(I,I,I)\}. $$
If $P\le{}S(q)$ is centric and radical in the fusion system of $H(q)$, 
then by Lemma \ref{str-p-rad}(c), its intersection with 
$H_0(q)\cong{}SL_2(q)^3/\{\pm(I,I,I)\}$ is centric and radical in the 
fusion system of that group.  So by Lemma \ref{str-p-rad}(a,f), 
	\beq P\cap{}H_0(q)\cong(P_1\times{}P_2\times{}P_3)/\{\pm(I,I,I)\} 
	\tag{5} \eeq
for some $P_i$ which are centric and radical in the fusion system of 
$SL_2(q)$.  Since the Sylow 2--subgroups of $SL_2(q)$ are quaternion 
\cite[Theorem 2.8.3]{Gorenstein}, the $P_i$ must be nonabelian and 
quaternion, so each $P_i/\{\pm{}I\}$ is centric
in $PSL_2(q^\infty)$.  Hence $P\cap{}H_0(q)$ is centric in 
$SL_2(q)^3/\{\pm(I,I,I)\}$ by (5), and so $P$ is 
centric in $H(q^\infty)$ by (4).  
\end{proof}

We would like to be able to regard $B\Spin_7(q)$ as a subcomplex of 
$B\Sol(q)$, but there is no simple natural way to do so.  Instead, we set 
	$$ \bspin_7(q)=|\lccspin(q)|\doscom \subseteq |\lccsol(q)|\doscom 
	\subseteq B\Sol(q); $$
then $\bspin_7(q)\simeq{}B\Spin_7(q)\doscom$ by \cite[Proposition 1.1]{BLO1} 
and Lemma \ref{Lcc=Lc}.  Also, we write
	$$ \bsol(q)=|\lccsol(q)|\doscom \subseteq B\Sol(q) \defeq 
	|\lcsol(q)|\doscom $$
to denote the subcomplex shown in Lemma \ref{Lcc=Lc} to be equivalent to 
$B\Sol(q)$; and set 
	$$ \bspin_7(q^\infty)=|\lcspin(q^\infty)|\doscom. $$
 From now on, when we talk about the inclusion of $B\Spin_7(q)$ into 
$B\Sol(q)$, as long as it need only be well defined up to homotopy, we 
mean the composite
	$$ B\Spin_7(q) \simeq \bspin_7(q) \subseteq \bsol(q) $$
(for some choice of homotopy equivalence).  Similarly, if we talk about 
the inclusion of $B\Sol(q^m)$ into $B\Sol(q^n)$ (for $m|n$) where it need 
only be defined up to homotopy, we mean these spaces identified with their 
equivalent subcomplexes $\bsol(q^m)\subseteq\bsol(q^n)$. 

\begin{Lem} \label{H*-pullback}
Let $q$ be any odd prime.  Then for all $n\ge1$,
        \beq \begin{diagram}[w=55pt]
        H^*(B\Sol(q^n);\F_2) & \rTo & H^*(BH(q^n);\F_2)^{C_3} \\
        \dTo && \dTo \\
        H^*(B\Spin_7(q^n);\F_2) & \rTo & H^*(BH(q^n);\F_2)
        \end{diagram} \tag{1} \eeq
(with all maps induced by inclusions of groups or spaces) is a pullback 
square.
\end{Lem}

\begin{proof}  By \cite[Theorem B]{BLO2}, $H^*(B\Sol(q^n);\F_2)$ is the ring 
of elements in the cohomology of $S(q^n)$ which are stable relative to 
the fusion.  By the construction in Section 2, the fusion in $\Sol(q^n)$ is 
generated by that in $\Spin_7(q^n)$, together with the permutation action of 
$C_3$ on the subgroup $H(q^n)\le\Spin_7(q^n)$, and hence (1) is a pullback 
square.
\end{proof}

\begin{Prop} \label{BSol(fqbar)}
For each odd prime $q$, there is a category $\lcsol(q^\infty)$, 
together with a functor
        $$ \pi\: \lcsol(q^\infty) \Right5{} \fsol(q^\infty), $$
such that the following hold:
\begin{enumerate}  
\item For each $n\ge1$, 
$\pi^{-1}(\fsol(q^n))\cong\lccsol(q^n)$.

\item There is a homotopy equivalence 
	$$ B\Sol(q^\infty) \defeq |\lcsol(q^\infty)|\doscom
	\RIGHT6{\eta}{\simeq} BDI(4) $$ 
such that the following square commutes up to homotopy
	\beq \begin{diagram}[w=55pt]
	\bspin_7(q^\infty)\doscom & \rTo^{\delta(q^\infty)} & 
	B\Sol(q^\infty) \\
	\dTo<{\eta_0}>{\simeq} && \dTo<{\eta}>{\simeq} \\
	B\Spin(7)\doscom & \rTo^{\widehat{\delta}} & BDI(4) \rlap{\,.}
	\end{diagram} \tag{1} \eeq
Here, $\eta_0$ is the homotopy equivalence of \cite{FM}, induced by some 
fixed choice of embedding of the Witt vectors for $\fqbar$ into $\C$, while 
$\delta(q^\infty)$ is the union of the inclusions 
$|\lccspin(q^n)|\doscom\subseteq|\lccsol(q^n)|\doscom$, and
$\widehat{\delta}$ is the inclusion arising from the construction of 
$BDI(4)$ in \cite{DW:DI4}.
\end{enumerate}
Furthermore, there is an automorphism 
$\psi^q_\call\in\Aut(\lcsol(q^\infty))$ of categories which 
satisfies the conditions:
\begin{enumerate} \setcounter{enumi}{2}
\item the restriction of $\psi^q_\call$ to each subcategory 
$\lccsol(q^n)$ is equal to the restriction of 
$\psi^q_\call\in\Aut(\lcsol(q^n))$ as defined in Proposition 
\ref{Sol(q)+}(d);
\item $\psi^q_\call$ covers the automorphism $\psi^q_\calf$ of 
$\fsol(q^\infty)$ induced by the field automorphism $(x\mapsto{}x^q)$; and
\item for each $n$, $(\psi^q_\call)^n$ fixes $\lccsol(q^n)$.
\end{enumerate}
\end{Prop}

\begin{proof}  By Proposition \ref{F_Sol(fqbar)}, the inclusions 
$\Spin_7(q^m)\le\Spin_7(q^n)$ for all $m|n$ induce inclusions of fusion 
systems $\fsol(q^m)\subseteq\fsol(q^n)$.  Since the restriction of a 
linking system over $\fccsol(q^n)$ is a linking system over $\fccsol(q^m)$, 
the uniqueness of linking systems (Proposition \ref{Sol(q)+}) implies that 
we get inclusions $\lccsol(q^m)\subseteq\lccsol(q^n)$.  We define 
$\lcsol(q^\infty)$ to be the union of the finite categories 
$\lccsol(q^n)$.  (More precisely, fix a sequence of positive integers 
$n_1|n_2|n_3|\cdots$ such that every positive integer divides some $n_i$, 
and set
	$$ \lcsol(q^\infty) = \bigcup_{i=1}^\infty \lccsol(q^{n_i}). $$
Then by uniqueness again, we can identify $\lccsol(q^n)$ for each $n$ with 
the appropriate subcategory.)  

Let $\pi\:\lcsol(q^\infty)\Right1{}\fsol(q^\infty)$ be the union of the 
projections from $\lccsol(q^{n_i})$ to 
$\fsol(q^{n_i})\subseteq\fsol(q^\infty)$.  Condition (a) is clearly 
satisfied.  Also, using Proposition \ref{Sol(q)+}(d) and Lemma 
\ref{Lcc=Lc}, we see that there is an automorphism $\psi^q_\call$ of 
$\lcsol(q^\infty)$ which satisfies conditions (c,d,e) above.  (Note that 
by the fusion theorem as shown in \cite[Theorem A.10]{BLO2}, morphisms in 
$\lcsol(q^n)$ are generated by those between radical subgroups, and hence 
by those in $\lccsol(q^n)$.)

It remains only to show that $|\lcsol(q^\infty)|\doscom\simeq{}BDI(4)$, 
and to show that square (1) commutes.  The space $BDI(4)$ is 
2--complete by its construction in \cite{DW:DI4}.  By Lemma \ref{Lcc=Lc}, 
	$$ H^*(B\Sol(q^\infty);\F_2) \cong
	\invlim{n}{} H^*\bigl(|\lcsol(q^n)|;\F_2\bigr) =
	\invlim{n}{} H^*\bigl(B\Sol(q^n);\F_2\bigr). $$
Hence by Lemma \ref{H*-pullback} (and since the inclusions 
$B\Spin_7(q^n)\Right1{}B\Sol(q^n)$ commute with the maps 
induced by inclusions of fields $\F_{q^m}\subseteq\F_{q^n}$), there is a 
pullback square
        \beq \begin{diagram}[w=55pt]
        H^*(B\Sol(q^\infty);\F_2) & \rTo & H^*(BH(q^\infty);\F_2)^{C_3} \\
        \dTo && \dTo \\
        H^*(B\Spin_7(q^\infty);\F_2) & \rTo & H^*(BH(q^\infty);\F_2) 
        \rlap{\,.}
        \end{diagram} \tag{2} \eeq
Also, by \cite[Theorem 1.4]{FM}, there are maps
	$$ B\Spin_7(q^\infty) \Right2{} B\Spin(7)
	\quad\textup{and}\quad
	BH(q^\infty) \Right2{} B\bigl(SU(2)^3/\{\pm(I,I,I)\}\bigr) $$
which induce isomorphisms of $\F_2$--cohomology, and hence homotopy 
equivalences after 2--completion.  So by Propositions \ref{polynomials} and 
\ref{pullbk} (or more directly by the computations in \cite[\S3]{DW:DI4}), 
the pullback of the above square is the ring of Dickson invariants in the 
polynomial algebra $H^*(BC_2^4;\F_2)$, and thus isomorphic to 
$H^*(BDI(4);\F_2)$.  

Point (b), including the commutativity of (1), now follows from the following 
lemma.
\end{proof}

\newcommand{\bfA}{\mathbf{A}}

\begin{Lem} \label{Spin7->DI4}
Let $X$ be a 2--complete space such that $H^*(X;\F_2)$ is the Dickson 
algebra in 4 variables.  Assume further that there is a map 
$B\Spin(7)\Right1{f}X$ such that $H^*(f|_{BC_2^4};\F_2)$ is the inclusion 
of the Dickson invariants in the polynomial algebra $H^*(BC_2^4;\F_2)$.  
Then $X\simeq{}BDI(4)$.  More precisely, there is a homotopy equivalence 
between these spaces such that the composite
	$$ B\Spin(7) \Right5{f} X \simeq BDI(4) $$
is the inclusion arising from the construction in \cite{DW:DI4}.
\end{Lem}

\begin{proof}  In fact, Notbohm \cite[Theorem 1.2]{Notbohm} has proven that 
the lemma holds even without the assumption about $B\Spin(7)$ (but with the
more precise assumption that $H^*(X;\F_2)$ is isomorphic as an algebra over 
the Steenrod algebra to the Dickson algebra).  The result as stated above 
is much more elementary (and also implicit in \cite{DW:DI4}), so we sketch 
the proof here.  

Since $H^*(X;\F_2)$ is a polynomial algebra, $H^*(\Omega{}X;\F_2)$ is 
isomorphic as a graded vector space to an exterior algebra on the same 
number of variables, and in particular is finite.  Hence $X$ is a 
2--compact group.  By \cite[Theorem 8.1]{DW-cent} (the centralizer 
decomposition for a $p$--compact group), there is an $\F_2$--homology 
equivalence
	$$ \hclim{\bfA}(\alpha) \Right5{\simeq} X. $$
Here, $\bfA$ is the category of pairs $(V,\varphi)$, where $V$ is a 
nontrivial elementary abelian 2--group, and $\varphi:BV\Right1{}X$ makes 
$H^*(BV;\F_2)$ into a finitely generated module over $H^*(X;\F_2)$ (see 
\cite[Proposition 9.11]{DW:p-compact}).  Morphisms in $\bfA$ are defined 
by letting $\Mor_{\bfA}((V,\varphi),(V',\varphi'))$ be the set of 
monomorphisms $V\Right1{}V'$ of groups which make the obvious triangle 
commute up to homotopy.  Also, 
	$$ \alpha\:\bfA\op\Right1{}\Top 
	\quad\textup{is the functor}\quad 
	\alpha(V,\varphi)=\map(BV,X)_\varphi. $$
By \cite[Lemma 1.6(1)]{DW:DI4} and \cite[Th\'eor\`eme 0.4]{Lannes}, $\bfA$ 
is equivalent to the category of elementary abelian 2--groups $E$ with 
$1\le\rk(E)\le4$, whose morphisms consist of all group monomorphisms.  
Also, if $BC_2\Right1{\varphi}X$ is the restriction of $f$ to any subgroup 
$C_2\le\Spin(7)$, then in the notation of Lannes,
	$$ T_{C_2}(H^*(X;\F_2);\varphi^*) \cong H^*(B\Spin(7);\F_2) $$
by \cite[Lemmas 16.(3), 3.10 and 3.11]{DW:DI4}, and hence
	$$ H^*(\map(BC_2,X)_\varphi;\F_2) \cong H^*(B\Spin(7);\F_2) $$
by Lannes \cite[Th\'eor\`eme 3.2.1]{Lannes}.  This shows that 
	$$ \Bigl(\map(BC_2,X)_\varphi\Bigr){}\doscom
	\simeq{}B\Spin(7)\doscom, $$
and thus that the centralizers of other elementary abelian 2--groups are 
the same as their centralizers in $B\Spin(7)\doscom$.  In other words, 
$\alpha$ is equivalent in the homotopy category to the diagram used in 
\cite{DW:DI4} to define $BDI(4)$.  By \cite[Proposition 7.7]{DW:DI4} (and 
the remarks in its proof), this homotopy functor has a unique homotopy 
lifting to spaces.  So by definition of $BDI(4)$,
	\beq X \simeq \bigl(\hclim{\bfA}(\alpha)\bigr){}\doscom 
	\simeq BDI(4). \tag*{\qed} \eeq
\def\qed{}\end{proof}

Set $B\psi^q\defeq|\psi^q_\call|$, a self homotopy equivalence of 
$B\Sol(q^\infty)\simeq{}BDI(4)$.  By construction, the restriction of 
$B\psi^q$ to the maximal torus of $B\Sol(q^\infty)$ is the map induced by 
$x\mapsto{}x^q$, and hence this is an ``Adams map'' as defined by 
Notbohm \cite{Notbohm}.  In fact, by \cite[Theorem 3.5]{Notbohm}, there is 
an Adams map from BDI(4) to itself, unique up to homotopy, of degree any 
2--adic unit.

Following Benson \cite{Benson}, we define $BDI_4(q)$ for any odd prime 
power $q$ to be the homotopy fixed point set of the $\Z$--action on 
$B\Sol(q^\infty)\simeq{}BDI(4)$ induced by the Adams map $B\psi^q$.  By 
``homotopy fixed point set'' in this situation, we mean that the following 
square is a homotopy pullback:
        \begin{diagram}[w=50pt]
        BDI_4(q) & \rTo & B\Sol(q^\infty) \\
        \dTo && \dTo<{\Delta} \\
        B\Sol(q^\infty) & \rTo^{(\Id,B\psi^q)} & 
        B\Sol(q^\infty)\times{}B\Sol(q^\infty).  
        \end{diagram}
The actual pullback of this square is the subspace $B\Sol(q)$ of elements 
fixed by $B\psi^q$, and we thus have a natural map 
$B\Sol(q)\Right1{\delta_0}BDI_4(q)$.

\begin{Thm} \label{BSol(q)=BDI4(q)}
For any odd prime power $q$, the natural map
        $$ B\Sol(q) \RIGHT6{\delta_0}{\simeq} BDI_4(q) $$
is a homotopy equivalence.
\end{Thm}

\begin{proof}  Since $BDI(4)$ is simply connected, the square used to 
define $BDI_4(q)$ remains a homotopy pullback square after 2--completion by 
\cite[II.5.3]{BK}.  Thus $BDI_4(q)$ is 2--complete.  Also, 
$B\Sol(q)\!\defeq\!|\lcsol(q)|\doscom$ is 2--complete since $|\lcsol(q)|$ is 
2--good \cite[Proposition 1.12]{BLO2}, and hence it suffices to prove that 
the map between these spaces is an $\F_2$--cohomology equivalence.  By 
Lemma \ref{H*-pullback}, this means showing that the following commutative 
square is a pullback square:
        \beq \begin{diagram}[w=60pt]
        H^*(BDI_4(q);\F_2) & \rTo & H^*(BH(q);\F_2)^{C_3} \\
        \dTo && \dTo \\
        H^*(B\Spin_7(q);\F_2) & \rTo & H^*(BH(q);\F_2)\,.
        \end{diagram} \tag{1} \eeq
Here, the maps are induced by the composite
        $$ B\Spin_7(q)\simeq\bspin_7(q)\doscom \subseteq B\Sol(q) 
        \Right5{} BDI_4(q) $$
and its restriction to $BH(q)$. Also, by Proposition \ref{BSol(fqbar)}(b), 
the following diagram commutes up to homotopy:
	\beq \begin{diagram}[w=45pt]
	B\Spin_7(q) & \rTo^{\incl} & B\Spin_7(q^\infty) & \rTo^{\eta_0} & 
	B\Spin(7) \\
	\dTo<{\delta(q)} && \dTo<{\delta(q^\infty)} && 
	\dTo<{\widehat{\delta}} \\
	B\Sol(q) & \rTo^{\incl} & B\Sol(q^\infty) & \rTo^{\eta} & BDI(4)
	\end{diagram} \tag{2} \eeq

By \cite[Theorem 12.2]{Friedlander}, together with \cite[\S1]{FM}, for any 
connected reductive Lie group 
$\gg$ and any algebraic epimorphism $\psi$ on $\gg(\fqbar)$ with finite 
fixed subgroup, there is a homotopy pullback square
        \beq \begin{diagram}  
        B(\gg(\fqbar)^\psi)\doscom & \rTo^{\incl} & B\gg(\fqbar)\doscom \\
        \dTo<{\incl} && \dTo<{\Delta} \\
        B\gg(\fqbar)\doscom & \rTo^{(\Id,B\psi)} & 
        B\gg(\fqbar)\doscom\times{}B\gg(\fqbar)\doscom \rlap{\,.}
        \end{diagram} \tag{3} \eeq
We need to apply this when $\gg=\Spin_7$ or 
$\gg=H=(SL_2)^3/\{\pm(I,I,I)\}$.  
In particular, if $\psi=\psi^q$ is the automorphism induced by the field 
automorphism $(x\mapsto{}x^q)$, then 
$\Spin_7(\fqbar)^\psi=\Spin_7(q)$ by Lemma \ref{Galois-Spin}, and 
$H(\fqbar)^\psi=H(q)\defeq{}H(\fqbar)\cap\Spin_7(q)$.  
We thus get a description of $B\Spin_7(q)$ and $BH(q)$ as homotopy 
pullbacks.  

By \cite[Theorem 1.4]{FM}, $B\gg(\fqbar)\doscom\simeq{}B\gg(\C)\doscom$.  
Also, we can replace the complex Lie groups $\Spin_7(\C)$ and $H(\C)$ by 
maximal compact subgroups $\Spin(7)$ and 
$\widebar{H}\defeq{}SU(2)^3/\{\pm(I,I,I)\}$, since these have the same 
homotopy type.  

If we set $\RR=H^*(B\gg(\fqbar);\F_2)\cong{}H^*(B\gg(\C);\F_2)$, then
there are Eilenberg-Moore spectral sequences
        \beq E_2 = \Tor^*_{\RR\otimes\RR\op}(\RR,\RR) \Longrightarrow 
        H^*(B(\gg(\fqbar)^\psi);\F_2); \eeq
where the $(\RR\otimes\RR\op)$--module structure on $\RR$ is defined by 
setting $(a\otimes{}b){\cdot}x=a{\cdot}x{\cdot}B\psi(b)$.  When 
$\gg=\Spin_7$ or $H$, then $\RR$ is a polynomial algebra by Proposition 
\ref{polynomials} and the above remarks, and $B\psi$ acts on $\RR$ via the 
identity.  The above spectral sequence thus satisfies the hypotheses of 
\cite[Theorem II.3.1]{Smith}, and hence collapses.  (Alternatively, note 
that in this case, $E_2$ is generated multiplicatively by $E_2^{0,*}$ and 
$E_2^{-1,*}$ by (5) below.)  Similarly, when $\RR=H^*(BDI(4);\F_2)$, there 
is an analogous spectral sequence which converges to $H^*(BDI_4(q);\F_2)$, 
and which collapses for the same reason. By the above remarks, these 
spectral sequences are natural with respect to the inclusions 
$BH(-)\subseteq{}B\Spin_7(-)$, and (using the naturality of $\psi^q$ shown 
in Proposition \ref{Sol(q)+}(d)) of $B\Spin_7(-)$ into $B\Sol(-)$ or 
$BDI(4)$.

To simplify the notation, we now write 
        $$ \A\defeq{}H^*(BDI(4);\F_2),\quad 
        \B\defeq{}H^*(B\Spin(7);\F_2),\quad\textup{and}\quad 
        \CC\defeq{}H^*(\widebar{H};\F_2) $$ 
to denote these cohomology rings.  The Frobenius automorphism $\psi^q$ 
acts via the identity on each of them.  We claim that the square 
        \beq \begin{diagram}[w=55pt] 
        \Tor^*_{\A\otimes\A\op}(\A,\A) & \rTo & 
        \Tor^*_{\CC\otimes\CC\op}(\CC,\CC)^{C_3} \\ 
        \dTo && \dTo \\ 
        \Tor^*_{\B\otimes\B\op}(\B,\B) & \rTo & 
        \Tor^*_{\CC\otimes\CC\op}(\CC,\CC) 
        \end{diagram} \tag{4} \eeq 
is a pullback square.  Once this has been shown, it then follows that in 
each degree, square (1) has a finite filtration under which each quotient 
is a pullback square.  Hence (1) itself is a pullback.

For any commutative $\F_2$--algebra $\RR$, let $\Omega_{\RR/\F_2}$ denote 
the $\RR$--module generated by 
elements $dr$ for $r\in\RR$ with the relations $dr=0$ if $r\in\F_2$, 
        $$ d(r+s)=dr+ds \qquad\textup{and}\qquad d(rs)=r{\cdot}ds+s{\cdot}dr. 
        $$
Let $\Omega^*_{\RR/\F_2}$ denote the \emph{ring of K\"ahler 
differentials}:  the exterior algebra (over $\RR$) of 
$\Omega_{\RR/\F_2}=\Omega^1_{\RR/\F_2}$.  When $\RR$ is a polynomial 
algebra, there are natural identifications
        \beq \Tor^*_{\RR\otimes\RR\op}(\RR,\RR) \cong HH_*(\RR;\RR) \cong 
        \Omega^*_{\RR/\F_2}. \tag{5} \eeq
The first isomorphism holds for arbitrary algebras, and is shown, e.g., in 
\cite[Lemma 9.1.3]{Weibel}.  The second holds for smooth algebras over a 
field \cite[Theorem 9.4.7]{Weibel} (and polynomial algebras are smooth as 
shown in \cite[\S9.3.1]{Weibel}).  In particular, the isomorphisms (5) 
hold for $\RR=\A,\B,\CC$, which are shown to be polynomial algebras in 
Proposition \ref{polynomials} below.  Thus, square (4) is isomorphic to 
the square
        \beq \begin{diagram}  
        \Omega^*_{\A/\F_2} & \rTo & \bigl(\Omega^*_{\CC/\F_2}\bigr){}^{C_3} \\
        \dTo && \dTo \\
        \Omega^*_{\B/\F_2} & \rTo & \Omega^*_{\CC/\F_2} \rlap{\,,}
        \end{diagram} \tag{6} \eeq
which is shown to be a pullback square in Propositions \ref{polynomials} 
and \ref{pullbk} below.
\end{proof}

It remains to prove that square (6) in the above proof is a pullback 
square.  In what follows, we let $D_i(x_1,\dots,x_n)$ denote the $i$-th 
Dickson invariant in variables $x_1,\dots,x_n$.  This is the 
$(2^n{-}2^{n-i})$-th symmetric polynomial in the elements (equivalently 
in the nonzero elements) of the vector space 
$\gen{x_1,\dots,x_n}_{\F_2}$. 
We refer to \cite{Wilker} for more detail.  Note that what he denotes
$c_{n,i}$ is what we call $D_{n-i}(x_1,\dots,x_n)$.

\begin{Lem}  \label{Dickson-rel}
For any $n$, 
        \begin{align*}  
        D_1(x_1,\dots,x_{n+1}) &= 
        \prod_{x\in\gen{x_1,\dots,x_n}_{\F_2}}(x_{n+1}+x) + 
        D_1(x_1,\dots,x_n)^2 
        \\ &= x_{n+1}^{2^n}+ 
        \sum_{i=1}^{n}x_{n+1}^{2^{n-i}}D_i(x_1,\dots,x_n) + 
        D_1(x_1,\dots,x_n)^2. 
        \end{align*}
\end{Lem}

\begin{proof}  The first equality is shown in \cite[Proposition 
1.3(b)]{Wilker}; here we prove them both simultaneously.  Set 
$V_n=\gen{x_1,\dots,x_n}_{\F_2}$.  Since $\sigma_i(V_n)=0$ whenever $2^n-i$ is 
not a power of $2$ (cf \cite[Proposition 1.1]{Wilker}),
        \begin{align*}  
        D_1(x_1,\dots,x_{n+1}) &= \sum_{i=0}^{2^n} 
        \sigma_i(V_n){\cdot}\sigma_{2^n-i}(x_{n+1}+V_n) \\
        &= \prod_{x\in{}V_n}(x_{n+1}+x) + 
        \sum_{i=1}^{n}D_i(x_1,\dots,x_n){\cdot}
        \sigma_{2^{n-i}}(x_{n+1}+V_n). 
        \end{align*}
Also, since $\sigma_i(V_n)=0$ for $0<i<2^{n-1}$ as noted above,
        $$ \sigma_k(x_{n+1}+V_n) = \sum_{i=0}^k x_{n+1}^{k-i}{\cdot}
        \begin{smallpmatrix}2^n-i\\k-i\end{smallpmatrix} \sigma_{i}(V_n)
        = \begin{cases} 0 & \textup{if $0<k<2^{n-1}$} \\
        D_1(x_1,\dots,x_n) & \textup{if $k=2^{n-1}$.} \end{cases} $$
This proves the first equality, and the second follows since
        \beq \prod_{x\in{}V_n}(x_{n+1}+x) = x_{n+1}^{2^n} + 
        \sum_{i=1}^{2^n}x_{n+1}^{2^n-i}\sigma_i(V_n) 
        = x_{n+1}^{2^n} + 
        \sum_{i=1}^{n}x_{n+1}^{2^{n-i}}D_i(x_1,\dots,x_n). \tag*{\qed} \eeq
\def\qed{}\end{proof}

In the following proposition (and throughout the rest of the section), we 
work with the polynomial ring $\F_2[x,y,z,w]$, with the natural action of 
$GL_4(\F_2)$.  Let
        $$ GL^2_2(\F_2), \ GL^3_1(\F_2) \le GL_4(\F_2) $$
be the subgroups of automorphisms of $V\defeq\gen{x,y,z,w}_{\F_2}$ which leave 
invariant the subspaces $\gen{x,y}$ and $\gen{x,y,z}$, respectively.  Also, let
$GL^2_{2'}(\F_2)\le{}GL^2_2(\F_2)$ be the subgroup of automorphisms 
which are the identity modulo $\gen{x,y}$.  Thus, when described 
in terms of block matrices (with respect to the given basis $\{x,y,z,w\}$),
        $$ GL^3_1(\F_2)=
        \left\{\begin{smallpmatrix}A&X\\0&1\end{smallpmatrix}\right\}
        ,\quad
        GL^2_2(\F_2)= 
        \left\{\begin{smallpmatrix}B&Y\\0&C\end{smallpmatrix}\right\},
        \quad \textup{and}\quad
        GL^2_{2'}(\F_2)= 
        \left\{\begin{smallpmatrix}B&Y\\0&I\end{smallpmatrix}\right\}, $$
for $A\in{}GL_3(\F_2)$, $X$ a column vector, $B,C\in{}GL_2(\F_2)$, and
$Y\in{}M_2(\F_2)$.  

We need to make more precise the relation between $V$ (or the polynomial 
ring $\F_2[x,y,z,w]$) and the cohomology of $\Spin(7)$.  To do this,
let $W\le\Spin(7)$ be the inverse image of the elementary abelian subgroup
        \begin{multline*}  
        {\bigl\langle}\diag(-1,-1,-1,-1,1,1,1), \diag(-1,-1,1,1,-1,-1,1), 
        \\ \diag(-1,1,-1,1,-1,1,-1) {\bigr\rangle}  \le \textup{SO}(7). 
        \end{multline*}
Thus, $W\cong{}C_2^4$.  Fix a basis $\{\eta,\eta',\xi,\zeta\}$ for $W$, 
where $\zeta\in{}Z(\Spin(7))$ is the nontrivial element.  Identify $V=W^*$ 
in such a way that 
$\{x,y,z,w\}\subseteq{}V$ is the dual basis to 
$\{\eta,\eta',\xi,\zeta\}$.  This gives an identification
        $$ H^*(BW;\F_2) = \F_2[x,y,z,w], $$
arranged such that the action of $N_{\Spin(7)}(W)/W$ on $V=\gen{x,y,z,w}$ 
consists of all automorphisms which leave $\gen{x,y,z}$ invariant, and 
thus can be identified with the action of $GL^3_1(\F_2)$.  Finally, set
        $$ \widebar{H} = C_{\Spin(7)}(\xi) \cong 
        \Spin(4)\times_{C_2}\Spin(3) \cong SU(2)^3/\{\pm(I,I,I)\} $$
(the central product). Then in the same way, the action of 
$N_{\widebar{H}}(W)/W$ on $H^*(BW;\F_2)$ can be identified with that of 
$GL^2_{2'}(\F_2)$.  

\begin{Prop} \label{polynomials}
The inclusions
        $$ BW \Right4{} B\widebar{H} \Right4{} B\Spin(7) \Right4{} 
        BDI(4) $$
as defined above, together with the identification 
$H^*(BW;\F_2)=\F_2[x,y,z,w]$, induce isomorphisms
        \beq \begin{split}  
        \A &\defeq{}H^*(BDI(4);\F_2) = \F_2[x,y,z,w]^{GL_4(\F_2)} 
        =\F_2[a_8,a_{12},a_{14},a_{15}] \\
        \B &\defeq{}H^*(B\Spin(7);\F_2) = \F_2[x,y,z,w]^{GL^3_1(\F_2)} = 
        \F_2[b_4,b_6,b_7,b_8] \\
        \CC &\defeq{}H^*(B\widebar{H};\F_2) = 
        \F_2[x,y,z,w]^{GL^2_{2'}(\F_2)} = \F_2[c_2,c_3,c'_4,c''_4] \,;
        \end{split} \tag{$*$} \eeq
where 
	\begin{multline*}  
        a_8=D_1(x,y,z,w),\  a_{12}=D_2(x,y,z,w),\ \\
        a_{14}=D_3(x,y,z,w),\  a_{15}=D_4(x,y,z,w); 
	\end{multline*}
        $$ b_4=D_1(x,y,z),\quad b_6=D_2(x,y,z),\quad
        b_7=D_3(x,y,z),\quad b_8=\prod_{\alpha\in\gen{x,y,z}}(w+\alpha); $$
and
        $$ c_2=D_1(x,y),\quad c_3=D_2(x,y), \quad 
        c'_4=\prod_{\alpha\in\gen{x,y}}(z+\alpha), \quad 
        c''_4=\prod_{\alpha\in\gen{x,y}}(w+\alpha). $$
Furthermore,
\begin{enumerate}  
\item the natural action of $\Sigma_3$ on 
$\widebar{H}\cong{}SU(2)^3/\{\pm(I,I,I)\}$ induces the action on $\CC$ 
which fixes $c_2,c_3$ and permutes $\{c'_4,c''_4,c'_4+c''_4\}$; and

\item the above variables satisfy the relations
        \begin{align*}  
        a_8&=b_8+b_4^2 &  a_{12}&=b_8b_4+b_6^2 & 
        a_{14}&=b_8b_6+b_7^2 &  a_{15}&=b_8b_7 \\
        b_4&=c'_4+c_2^2 &  b_6&=c_2c'_4+c_3^2 & 
        b_7&=c_3c'_4 &  b_8&=c''_4(c'_4+c''_4) \,. 
        \end{align*}
\end{enumerate}
\end{Prop}

\begin{proof}  The formulas for $\A=H^*(BDI(4);\F_2)$ are shown in 
\cite{DW:DI4}.  From \cite[Lemmas 3.10 and 3.11]{DW:DI4}, we see there are 
(some) identifications
	\begin{small}  
        $$ H^*(B\Spin(7);\F_2)\cong\F_2[x,y,z,w]^{GL^3_1(\F_2)}
        \quad\textup{and}\quad
        H^*(B\widebar{H};\F_2)\cong\F_2[x,y,z,w]^{GL^2_{2'}(\F_2)}. $$
	\end{small}%
{}From the explicit choices of subgroups $W\le\widebar{H}\le\Spin(7)$ as
described above (and by the descriptions in 
Proposition \ref{elem.abel.} of the automorphism groups), the images of 
$H^*(B\Spin(7);\F_2)$ and 
$H^*(B\widebar{H};\F_2)$ in $\F_2[x,y,z,w]$ are seen to be contained in 
the rings of invariants, and hence these isomorphisms actually are 
equalities as claimed. 

We next prove the equalities in ($*$) between the given rings of 
invariants and polynomial algebras.  The following argument was shown to 
us by Larry Smith.  If $k$ is a field and $V$ is an $n$--dimensional vector 
space over $k$, then a \emph{system of parameters} in the polynomial 
algebra $k[V]$ is a set of $n$ homogeneous elements $f_1,\dots,f_n$ such 
that $k[V]/(f_1,\dots,f_n)$ is finite dimensional over $k$.  By 
\cite[Proposition 5.5.5]{Smith2}, if $V$ is an $n$--dimensional 
$k[G]$--representation, and $f_1,\dots,f_n\in{}k[V]^G$ is a system of 
parameters the product of whose degrees is equal to $|G|$, then $k[V]^G$ 
is a polynomial algebra with $f_1,\dots,f_n$ as generators.  By 
\cite[Proposition 8.1.7]{Smith2}, $\F_2[x,y,z,w]$ is a free finitely 
generated module over the ring generated by its Dickson invariants (this 
holds for polynomial algebras over any $\F_p$), and thus 
$\F_2[x,y,z,w]/(a_8,a_{12},a_{14},a_{15})$ is finite.  (This can also be 
shown directly using the relation in Lemma \ref{Dickson-rel}.)  So 
assuming the relations in point (b), the quotients 
$\F_2[x,y,z,w]/(b_4,b_6,b_7,b_8)$ and $\F_2[x,y,z,w]/(c_2,c_3,c'_4,c''_4)$ 
are also finite.  In each case, the product of the degrees of the 
generators is clearly equal to the order of the group in question, and 
this finishes the proof of the last equality in the second and third lines 
of ($*$).

It remains to prove points (a) and (b).  Using Lemma \ref{Dickson-rel}, 
the $c_i$ are expressed as polynomials in $x,y,z,w$ as follows:
	\begin{small}  
        \beq \begin{split}  
        c_2 &= D_1(x,y)=x^2+xy+y^2 \\
        c_3 &= D_2(x,y)=xy(x+y) \\
        c'_4 &= D_1(x,y,z)+D_1(x,y)^2 = z^4 + z^2D_1(x,y) + zD_2(x,y) 
        = z^4 + z^2c_2 + zc_3 \\
        c''_4 &= D_1(x,y,w)+D_1(x,y)^2 = w^4 + w^2D_1(x,y) + wD_2(x,y) 
        = w^4 + w^2c_2 + wc_3 \,.
        \end{split} \tag{1} \eeq
	\end{small}%
In particular,
        \beq c'_4+c''_4 = (z+w)^4 + (z+w)^2D_1(x,y) + (z+w)D_2(x,y) 
        = \prod_{\alpha\in\gen{x,y}}(z+w+\alpha). \tag{2} \eeq
Furthermore, by (1), we get
        \beq \begin{split}  
        Sq^1(c_2)&=c_3 \\
        Sq^1(c_3)&=Sq^1(c'_4)=Sq^1(c''_4)=0 \\
        Sq^2(c_3)&=x^2y^2(x+y)+xy(x+y)^3=c_2c_3 \\
        Sq^2(c'_4) &= z^4c_2 + z^2c_2^2 + zc_2c_3 = c_2c'_4 \\
        Sq^3(c'_4) &= Sq^1(c_2c'_4) = c_3c'_4 \\
        Sq^2(c''_4) &= c_2c''_4 \\
        Sq^3(c''_4) &= c_3c''_4. 
        \end{split} \tag{3} \eeq

The permutation action of $\Sigma_3$ on 
$\widebar{H}\cong{}SU(2)^3/\{\pm(I,I,I)\}$ permutes the three elements 
$\zeta,\xi,\zeta+\xi$ of $Z(\widebar{H})\subseteq{}W$, and thus (via the 
identification $V=W^*$ described above) induces the identity on 
$x,y\in{}V$ and permutes the elements $\{z,w,z+w\}$ modulo $\gen{x,y}$.  
Hence the induced action of $\Sigma_3$ on $\CC=\F_2[V]^{GL^2_{2'}(\F_2)}$ 
is the restriction of the action on $\F_2[V]=\F_2[x,y,z,w]$ which fixes 
$x,y$ and permutes $\{z,w,z+w\}$.  So by (1) and (2), we see that this 
action fixes $c_2,c_3$ and permutes the set $\{c'_4,c''_4,c'_4+c''_4\}$.  
This proves (a). 

It remains to prove the formulas in (b).  From (1) and (3) we get
        \begin{align*}  
        b_4 &= D_1(x,y,z) = c'_4+c_2^2, \\
        b_6 &= D_2(x,y,z) = Sq^2(b_4) = c_2c'_4+c_3^2, \\
        b_7 &= D_3(x,y,z) = Sq^1(b_6) = c_3c'_4. 
        \end{align*}
Also, by (1) and (2),
        $$ b_8 = \prod_{\alpha\in\gen{x,y,z}}(w+\alpha) =
        \Bigl(\prod_{\alpha\in\gen{x,y}}(w+\alpha) \Bigr) {\cdot}
        \Bigl(\prod_{\alpha\in\gen{x,y}}(w+z+\alpha) \Bigr) 
        = c''_4(c'_4+c''_4). $$
This proves the formulas for the $b_i$ in terms of $c_i$.  Finally, we have
        \begin{align*}  
        a_8 &= D_1(x,y,z,w) = b_8+b_4^2, \\
        a_{12} &= D_2(x,y,z,w) = Sq^4(b_8+b_4^2) = 
        Sq^4(c''_4(c'_4+c''_4)+(c'_4+c_2^2)^2) \\
        &= c'_4c''_4(c'_4+c''_4) + c_2^2c''_4(c'_4+c''_4) + c_2^2c'_4{}^2 + 
        c_3^4 = b_8b_4+b_6^2 \\
        a_{14} &= D_3(x,y,z,w) = Sq^2(a_{12}) = c_2c'_4c''_4(c'_4+c''_4) 
        + c_3^2c''_4(c'_4+c''_4) + c_3^2c'_4{}^2 \\
        &= b_8b_6+b_7^2 \\
        a_{15} &= D_4(x,y,z,w) = Sq^1(a_{14}) = c_3c'_4c''_4(c'_4+c''_4) 
        = b_8b_7 \,;
        \end{align*}
and this finishes the proof of the proposition.
\end{proof}

\begin{Lem} \label{inA}
Let $\kappa\in\Aut(\CC)$ be the algebra involution which exchanges $c'_4$ 
and $c''_4$ and leaves $c_2$ and $c_3$ fixed.  An element of $\CC$ will be 
called ``$\kappa$--invariant'' if it is fixed by this involution.  Then the 
following hold:
\begin{enumerate}  
\item If $\beta\in{}\B$ is $\kappa$--invariant, then $\beta\in{}\A$.
\item If $\beta\in{}\B$ is such that $\beta{\cdot}c'_4{}^i$ is 
$\kappa$--invariant, then $\beta=\beta'{\cdot}b_8^i$ for some $\beta'\in{}\A$.
\end{enumerate}
\end{Lem}

\begin{proof}  Point (a) follows from Proposition \ref{polynomials} upon 
regarding $\A$, $\B$, and $\CC$ as the fixed subrings of the groups 
$GL_4(\F_2)$, $GL^3_1(\F_2)$ and $GL^2_{2'}(\F_2)$ acting on 
$\F_2[x,y,z,w]$, but also follows from the following direct argument.  Let 
$m$ be the degree of $\beta$ as a polynomial in $b_8$; we argue by 
induction on $m$.  Write $\beta=\beta_0+b_8^m{\cdot}\beta_1$, where 
$\beta_1\in\F_2[b_4,b_6,b_7]$, and where $\beta_0$ has degree $<m$ (as a 
polynomial in $b_8$).  If $m=0$, then 
$\beta=\beta_1\in\F_2[b_4,b_6,b_7]\subseteq\F_2[c_2,c_3,c'_4]$, and hence 
$\beta\in\F_2[c_2,c_3]$ since it is $\kappa$--invariant.  But from the 
formulas in Proposition \ref{polynomials}(b), we see that 
$\F_2[b_4,b_6,b_7]\cap\F_2[c_2,c_3]$ contains only constant polynomials 
(hence it is contained in $\A$).  

Now assume $m\ge1$.  Then, expressed as a polynomial in 
$c_2,c_3,c'_4,c''_4$, the largest power of $c''_4$ which occurs in $\beta$ 
is $c''_4{}^{2m}$.  Since $\beta$ is $\kappa$--invariant, the highest power 
of $c'_4$ which occurs is $c'_4{}^{2m}$; and hence by Proposition 
\ref{polynomials}(b), the total degree of each term in $\beta_1$ (its 
degree as a polynomial in $b_4,b_6,b_7$) is at most $m$.  So for each 
term $b_4^rb_6^sb_7^t$ in $\beta_1$, 
        $$ b_4^rb_6^sb_7^tb_8^m - 
        a_8^{m-r-s-t}a_{12}^ra_{14}^sa_{15}^t $$
is a sum of terms which have degree $<m$ in $b_8$, and thus lies in $\A$ 
by the induction hypothesis. 

To prove (b), note first that since $\beta{\cdot}c'_4{}^i$ is 
$\kappa$--invariant and divisible by $c'_4{}^i$, it must also be divisible 
by $c''_4{}^i$, and hence $c''_4{}^i|\beta$.  Furthermore, by Proposition 
\ref{polynomials}, all elements of $\B$ as well as $c'_4$ are invariant 
under the involution which fixes $c'_4$ and sends 
$c''_4\mapsto{}c'_4+c''_4$.  Thus $(c'_4+c''_4)^i|\beta$. Since 
$b_8=c''_4(c'_4+c''_4)$, we can now write $\beta=\beta'{\cdot}b_8^i$ for 
some $\beta\in{}\B$.  Finally, since
        $$ \beta{\cdot}c'_4{}^i = 
        \beta'{\cdot}c'_4{}^i{\cdot}c''_4{}^i{\cdot}(c'_4+c''_4)^i $$
is $\kappa$--invariant, $\beta'$ is also $\kappa$--invariant, and hence 
$\beta'\in\A$ by (a).  
\end{proof}

Note that $C_3\le\Sigma_3=GL_2(\F_2)$ act on 
$\CC=\F_2[x,y,z,w]^{GL^2_{2'}(\F_2)}$:  via the action of the group
$GL^2_2(\F_2)/GL^2_{2'}(\F_2)$, or equivalently by permuting $c'_4$, 
$c''_4$, and $c'_4+c''_4$ (and fixing $c_2,c_3$).  Thus 
$\A=\B\cap{}\CC^{C_3}$, since $GL_4(\F_2)$ is generated by the subgroups 
$GL^3_1(\F_2)$ and $GL^2_2(\F_2)$.  This is also shown directly in the 
following lemma.

\begin{Prop} \label{pullbk}
The following square is a pullback square, where all maps are induced by 
inclusions between the subrings of $\F_2[x,y,z,w]$:
        \begin{diagram}  
        \Omega^*_{\A/\F_2} & \rTo & \bigl(\Omega^*_{\CC/\F_2}\bigr)^{C_3} \\
        \dTo && \dTo \\
        \Omega^*_{\B/\F_2} & \rTo & \Omega^*_{\CC/\F_2} \rlap{\,.}
        \end{diagram}
\end{Prop}

\begin{proof}  Let $\kappa$ be the involution of Lemma \ref{inA}:  the 
algebra involution of $\CC$ which
exchanges $c'_4$ and $c''_4$ and leaves $c_2$ and $c_3$ fixed.  By 
construction, all elements in the image of $\Omega^*_{\B/\F_2}$ are 
invariant under the involution which fixes $c'_4$ (and $c_2,c_3$), and 
sends $c''_4$ to $c'_4+c''_4$.  Hence elements in the image of 
$\Omega^*_{\B/\F_2}$ are fixed by $C_3$ if and only if they are fixed by 
$\Sigma_3$, if and only if they are $\kappa$--invariant.  So it will 
suffice to show that all of the above maps are injective, and that all 
$\kappa$--invariant elements in the image of $\Omega^*_{\B/\F_2}$ lie in 
the image of $\Omega^*_{\A/\F_2}$.  The injectivity is clear, and the 
square is a pullback for $\Omega^0_{-/\F_2}$ by Lemma \ref{inA}.

Fix a $\kappa$--invariant element 
	\beq \begin{split}  
        \omega &= P_1\,db_4 + P_2\,db_6 + P_3\,db_7 + P_4\,db_8 \\
        &= P_2c'_4\,dc_2 + P_3c'_4\,dc_3 + P_4c'_4\,dc''_4 +
        (P_1+P_2c_2+P_3c_3+P_4c''_4)\,dc'_4
        \in \Omega^1_{\B/\F_2},
	\end{split} \tag{1} \eeq
where $P_i\in{}\B$ for each $i$.  
By applying $\kappa$ to (1) and comparing the coefficients
of $dc_2$ and $dc_3$, we see that $P_2c'_4$ and $P_3c'_4$ are 
$\kappa$--invariant.  
Also, upon comparing the
coefficients of $dc'_4$, we get the equation
	\beq P_1 + P_2c_2 + P_3c_3 + P_4c''_4 = \kappa(P_4)c''_4. \tag{2} \eeq
So by Lemma \ref{inA}, $P_2=P'_2b_8$ and $P_3=P'_3b_8$ 
for 
some $P'_2,P'_3\in{}\A$.  Upon subtracting 
        $$ P'_2\,da_{14}+P'_3\,da_{15} = P_2\,db_6 + P_3\,db_7 + 
        (P'_2b_6+P'_3b_7)\,db_8 
        $$
from $\omega$ and introducing an appropriate modification to $P_4$, we can 
now assume that $P_2=P_3=0$.  With this assumption and (2), we have
        $$ P_1+P_4c''_4 = \kappa(P_4c'_4) = \kappa(P_4){\cdot}c''_4, $$
so that 
        \beq P_1c'_4 = (P_4+\kappa(P_4))c'_4c''_4 \tag{3} \eeq
is $\kappa$--invariant.  This now shows that $P_1=P'_1b_8$ for some 
$P'_1\in{}\A$, and upon subtracting $P'_1\,da_{12}$ from $\omega$ we can 
assume that $P_1=0$. This leaves $\omega=P_4\,db_8=P_4\,da_8$.  By (3) 
again, $P_4$ is $\kappa$--invariant, so $P_4\in{}\A$ by Lemma \ref{inA} again, 
and thus $\omega\in\Omega^1_{\A/\F_2}$.

The remaining cases are proved using the same techniques, and so we sketch 
them more briefly.  To prove the result in degree two, fix a 
$\kappa$--invariant element
        \begin{align*}  
        \omega &= P_1\,db_4\,db_6 + P_2\,db_4\,db_7 + P_3\,db_4\,db_8 + 
        P_4\,db_6\,db_7 + 
        P_5\,db_6\,db_8 + 
        P_6\,db_7\,db_8 \\
        &= P_4c'_4{}^2\,dc_2\,dc_3 +
        (P_1c'_4+P_4c_3c'_4+P_5c'_4c''_4)\,dc_2\,dc'_4 +
        P_5c'_4{}^2\,dc_2\,dc''_4 \\
        &\qquad + (P_2c'_4+P_4c_2c'_4+P_6c'_4c''_4)\,dc_3\,dc'_4 +
        P_6c'_4{}^2\,dc_3\,dc''_4 \\
        &\qquad + (P_3c'_4+P_5c_2c'_4+P_6c_3c'_4)\,dc'_4\,dc''_4 \in 
        \Omega^2_{\B/\F_2}.
        \end{align*}
Using Lemma \ref{inA}, we see that $P_4=P'_4b_8^2$, and hence can assume 
that $P_4=0$.  One then eliminates $P_1$ and $P_2$, then $P_5$ and $P_6$, 
and finally $P_3$.  

If 
        \begin{align*}  
        \omega &= P_1\,db_4\,db_6\,db_7 + P_2\,db_4\,db_6\,db_8 + 
        P_3\,db_4\,db_7\,db_8 + P_4\,db_6\,db_7\,db_8 \\
        &= (P_1c'_4{}^2+P_4c'_4{}^2c''_4)\,dc_2\,dc_3\,dc'_4 + 
        (P_2c'_4{}^2+P_4c_3c'_4{}^2)\,dc_2\,dc'_4\,dc''_4 \\
        &\qquad + (P_3c'_4{}^2+P_4c_2c'_4{}^2)\,dc_3\,dc'_4\,dc''_4 + 
        P_4c'_4{}^3\,dc_2\,dc_3\,dc''_4 \in 
        \Omega^3_{\B/\F_2} 
        \end{align*}
is $\kappa$--invariant, then we eliminate successively $P_1$, then $P_4$, then 
$P_2$ and $P_3$.  

Finally, if 
        $$ \omega = P\,db_4\,db_6\,db_7\,db_8 = 
        Pc'_4{}^3\,dc_2\,dc_3\,dc'_4\,dc''_4 \in \Omega^4_{\B/\F_2} $$
is $\kappa$--invariant, then $P=P'b_8^3$ for some $P'\in{}\A$ by Lemma 
\ref{inA} 
again, and so 
	\beq \omega=P'\,da_8\,da_{12}\,da_{14}\,da_{15}
	\in\Omega^4_{\A/\F_2}. \tag*{\qed} \eeq
\def\qed{}\end{proof}

\vglue 0.3in
\noindent$\phantom{XX}$

\appendix

\section{Appendix : Spinor groups over finite fields}
\label{Appx:Spin}

Let $F$ be any field of characteristic $\ne2$.  Let $V$ be a vector space 
over $F$, and let $\bb\:V\Right1{}F$ be a nonsingular quadratic form.  As 
usual, $O(V,\bb)$ denotes the group of isometries of $(V,\bb)$, and 
$SO(V,\bb)$ the subgroup of isometries of determinant 1.  We will be 
particularly interested in elementary abelian 2--subgroups of such 
orthogonal groups.

\begin{Lem}  \label{2^n<SO}
Fix an elementary abelian 2--subgroup $E\le{}O(V,\bb)$.  For each irreducible 
character $\chi\in{}\Hom(E,\{\pm1\})$, let $V_\chi\subseteq{}V$ denote 
the corresponding eigenspace:  the subspace of elements $v\in{}V$ such that 
$g(v)=\chi(g){\cdot}v$ for all $g\in{}E$.  Then the restriction of $\bb$ to 
each subspace $V_\chi$ is nonsingular, and $V$ is the orthogonal direct 
sum of the $V_\chi$. 
\end{Lem}

\begin{proof}  Elementary.
\end{proof}

We give a very brief sketch of the definition of spinor groups via Clifford 
algebras; for more details we refer to \cite[\S{}II.7]{Dieudonne} or
\cite[\S22]{Aschbacher}.  Let $T(V)$ denote the tensor algebra of $V$, and 
set
	$$ C(V,\bb) = T(V)/\gen{(v\otimes{}v)-\bb(v)}\,: $$
the Clifford algebra of $(V,\bb)$.  To simplify the notation, we regard 
$F$ as a subring of $C(V,\bb)$, and $V$ as a subgroup of its additive 
group; thus the class of $v_1\otimes\cdots\otimes{}v_k$ will be written 
$v_1{\cdots}v_k$.  Note that $vw+wv=0$ if $v,w\in{}V$ and $v\perp{}w$.  
Hence if $\dim_F(V)=n$, and $\{v_1,\dots,v_n\}$ is an orthogonal basis, 
then the set of $1$ and all $v_{i_1}\cdots{}v_{i_k}$ for $i_1<\cdots<i_k$ 
($1\le{}k\le{}n$) is an $F$--basis for $C(V,\bb)$.  

Write $C(V,\bb)=C_0\oplus{}C_1$, where $C_0$ and $C_1$ consist of classes 
of elements of even or odd degree, respectively.  Let $G\le{}C(V,\bb)^*$ 
denote the group of invertible elements $u$ such that $uVu^{-1}=V$, and 
let $\pi\:G\Right1{}O(V,\bb)$ be the homomorphism
	$$ \pi(u) = \begin{cases}  
	(v\mapsto-uvu^{-1}) & \textup{if $u\in{}C_1$} \\
	(v\mapsto uvu^{-1}) & \textup{if $u\in{}C_0$\,.} 
	\end{cases} $$
In particular, for any nonisotropic element $v\in{}V$ (ie, 
$\bb(v)\ne0$), $v\in{}G$ and $\pi(v)$ is the reflection in the hyperplane 
$v^\perp$.  By \cite[\S{}II.7]{Dieudonne}, $\pi$ is surjective and 
$\Ker(\pi)=F^*$.  

Let $J$ be the antiautomorphism of $C(V,\bb)$ induced by the 
antiautomorphism 
$v_1\otimes\cdots\otimes{}v_k\mapsto{}v_k\otimes\cdots\otimes{}v_1$ of 
$T(V)$.  Since $O(V,\bb)$ is generated by hyperplane reflections, $G$ is 
generated by $F^*$ and nonisotropic elements $v\in{}V$.  In particular, for 
any $u=\lambda{\cdot}v_1\cdots{}v_k\in{}G$,
	$$ J(u){\cdot}u = \lambda^2\cdot v_k\cdots{}v_1\cdot{}v_1\cdots v_k
	= \lambda^2\cdot\bb(v_1)\cdots\bb(v_k) \in F^*=\Ker(\pi); $$
implying that $\pi(J(u))=\pi(u)^{-1}$ for all $u\in G$.  There is thus a 
homomorphism
	$$ \widetilde{\theta}\: G \Right6{} F^* 
	\qquad\textup{defined by}\qquad 
	\widetilde{\theta}(u)=u{\cdot}J(u). $$
In particular, $\widetilde{\theta}(\lambda)=\lambda^2$ for 
$\lambda\in{}F^*\le{}G$, while for any set of nonisotropic elements 
$v_1,\dots,v_k$ of $V$, 
	$$ \widetilde{\theta}(v_1\cdots v_k) = (v_1\cdots v_k)(v_k\cdots 
	v_1) = \bb(v_1)\cdots\bb(v_k). $$
Hence $\widetilde{\theta}$ factors through a homomorphism
	$$ \theta_{V,\bb} \: O(V,\bb) \Right6{} F^*/F^{*2}
	= F^*/\{u^2\,|\,u\in{}F^*\}, $$ 
called the \emph{spinor norm}.

Set $G^+=\pi^{-1}(SO(V,\bb))=G\cap{}C_0$, and define
	$$ \Spin(V,\bb)=\Ker(\widetilde{\theta}|_{G^+})
	\qquad\textup{and}\qquad
	\Omega(V,\bb)=\Ker(\theta_{V,\bb}|_{SO(V,\bb)}). $$
In particular, $\Omega(V,\bb)$ has index 2 in $SO(V,\bb)$ if $F$ is a 
finite field, and $\Omega(V,\bb)=SO(V,\bb)$ if $F$ is algebraically closed 
(all units are squares).  We thus get a commutative diagram
\stepcounter{Thm}\def\bigdiag{\hbox{\theThm}}

\smallskip

	\hbox to \textwidth{ (\theThm)\hfill
	$\xymx[@C=13pt@R=13pt]{
	& 1 \ar[d] && 1 \ar[d] && 1 \ar[d] \\
	1 \ar[r] & \{\pm1\} \ar[rr] \ar[dd] && F^* 
	\ar[rr]^{\lambda\mapsto\lambda^2} \ar[dd] && 
	F^{*2} \ar[dd] \ar[r] & 1 \\ \\
	1 \ar[r] & \Spin(V,\bb) \ar[rr] \ar[dd] && G^+ 
	\ar[rr]^{\widetilde{\theta}} 
	\ar[dd]_{\pi} && F^* \ar[dd] \ar[r] & 1 \\ \\
	1 \ar[r] & \Omega(V,\bb) \ar[rr] \ar[d] && SO(V,\bb) 
	\ar[rr]^{\theta_{V,\bb}} \ar[d] && 
	F^*/F^{*2} \ar[r] \ar[d] & 1 \\
	& 1 && 1 && 1 
	}$ \hfill } 
	\edef\bigdiag{\theThm}%

\smallskip

\noindent 
where all rows and columns are short exact, and where all columns are 
central extensions of groups. If $\dim(V)\ge3$ (or if $\dim(V)=2$ and the 
form $\bb$ is hyperbolic), then $\Omega(V,\bb)$ is the commutator subgroup 
of $SO(V,\bb)$ \cite[\S{}II.8]{Dieudonne}.  

The following lemma follows immediately from this description of 
$\Spin(V,\bb)$, together with the analogous description of the 
corresponding spinor group over the algebraic closure of $F$.  

\begin{Lem} \label{Galois-Spin}
Let $\widebar{F}$ be the algebraic closure of $F$, and set 
$\widebar{V}=\widebar{F}\otimes_{F}V$ and 
$\widebar{\bb}=\Id_{\widebar{F}}\otimes\bb$.  Then $\Spin(V,\bb)$ is the 
subgroup of $\Spin(\widebar{V},\widebar{\bb})$ consisting of those 
elements fixed by all Galois automorphisms 
$\psi\in\textup{Gal}(\widebar{F}/F)$. 
\hfill \qed
\end{Lem}

For any nonsingular quadratic form $\bb$ on a vector space $V$, the 
\emph{discriminant} of $\bb$ (or of $V$) is the determinant of the 
corresponding symmetric bilinear form $B$, related to $\bb$ by the formulas 
	$$ \bb(v)=B(v,v) \qquad\textup{and}\qquad
	B(v,w)=\tfrac12\bigl(\bb(v+w)-\bb(v)-\bb(w)\bigr). $$
Note that the discriminant is well defined only modulo squares in $F^*$.  
When $W\subseteq{}V$ is a subspace, then we define the discriminant of $W$ 
to mean the discriminant of $\bb|_W$. In what follows, we say that the 
discriminant of a quadratic form is a square or a nonsquare to mean that 
it is the identity or not in the quotient group $F^*/F^{*2}$.  

\begin{Lem}  \label{SO-invol}
Fix an involution $x\in{}SO(V,\bb)$, and let $V=V_+\oplus{}V_-$ be its 
eigenspace decomposition.  Then the following hold.
\begin{enumerate}  
\item $x\in{}\Omega(V,\bb)$ if and only if the discriminant of $V_-$ is a 
square.
\item If $x\in\Omega(V,\bb)$, then it lifts to an element of order 2 in 
$\Spin(V,\bb)$ if and only if $\dim(V_-)\in4\Z$.
\item If $x\in\Omega(V,\bb)$, and if $\alpha\in\Omega(V,\bb)$ is such that 
$[x,\alpha]=1$, then $\alpha=\alpha_+\oplus\alpha_-$, where 
$\alpha_{\pm}\in{}O(V_{\pm},\bb)$.  Also, the liftings of $x$ and $\alpha$ 
commute in $\Spin(V,\bb)$ if and only if $\det(\alpha_-)=1$.
\end{enumerate}
\end{Lem}

\begin{proof}  Let $\{v_1,\dots,v_k\}$ be an orthogonal basis for $V_-$ 
($k$ is even).  Then $x=\pi(v_1{\cdots}v_k)$ in the above notation, since 
$\pi(v_i)$ is the reflection in the hyperplane $v_i{}^\perp$.  Hence by the 
commutativity of Diagram (\bigdiag),
	$$ \theta_{V,\bb}(x) \equiv \bb(v_1){\cdots}\bb(v_k) = 
	\det(\bb|_{V_-}) \pmod{F^{*2}}. $$
Thus $x\in\Omega(V,\bb)=\Ker(\theta_{V,\bb})$ if and only if $V_-$ has 
square discriminant.  

In particular, if $x\in\Omega(V,\bb)$, then the product of the $\bb(v_i)$ 
is a square, and hence (upon replacing $v_1$ by a scalar multiple) we can 
assume that 
$\bb(v_1){\cdots}\bb(v_k)=1$.  Then $\widetilde{x}\defeq{}v_1\cdots{}v_k 
\in\Spin(V,\bb)=\Ker(\widetilde{\theta})$.  Since $vw=-wv$ in the Clifford 
algebra whenever $v\perp{}w$, and since $(v_i)^2=\bb(v_i)$ for each $i$,
	$$ \widetilde{x}{}^2 = (-1)^{k(k-1)/2}{\cdot}(v_1)^2{\cdots} 
	(v_k)^2 = (-1)^{k(k-1)/2} = 
	\begin{cases}  1 & \textup{if $k\equiv0$ (mod $4$)} \\
	-1 & \textup{if $k\equiv2$ (mod $4$)\,.} 
	\end{cases} $$
This proves (b).

It remains to prove (c).  The first statement 
($\alpha=\alpha_+\oplus\alpha_-$) is clear.  Fix liftings 
$\widetilde{\alpha}_\pm\in{}C(V_\pm,\bb)^*$.  Rather than defining the 
direct sum of an element of $C(V_+,\bb)$ with an element of $C(V_-,\bb)$, 
we regard the groups $C(V_\pm,\bb)^*$ as (commuting) subgroups of 
$C(V,\bb)^*$, and set
	$$ \widetilde{\alpha}= \widetilde{\alpha}_+\circ\widetilde{\alpha}_-
	= \widetilde{\alpha}_-\circ\widetilde{\alpha}_+ \in \Spin(V,\bb). $$
Let $\widetilde{x}=v_1{\cdots}v_k$ be as above.  
Clearly, $\widetilde{x}$ commutes with all elements of 
$C(V_+,\bb)$.  Since
	$$ (v_1{\cdots}v_k){\cdot}v_i 
	= (-1)^{k-1}{\cdot}v_i{\cdot}(v_1{\cdots}v_k)
	= -v_i{\cdot}(v_1{\cdots}v_k) $$
for $i=1,\dots,k$, we have
$\widetilde{x}{\cdot}\beta=(-1)^i{\cdot}\beta{\cdot}\widetilde{x}$ for all 
$\beta\in{}C_i(V_-,\bb)$ ($i=0,1$).  In particular, since 
$[\widetilde{\alpha}_+,\widetilde{\alpha}_-]=1$, 
$[\widetilde{x},\widetilde{\alpha}]=[\widetilde{x},\widetilde{\alpha}_-] 
=\det(\alpha_-)$, and this finishes the proof.
\end{proof}

We will need explicit isomorphisms which describe $\Spin_3(F)$ and  
$\Spin_4(F)$ in terms of $SL_2(F)$.  These are constructed in the 
following proposition, where $M_2^0(F)$ denotes the vector space of 
matrices of trace zero.  Note that the determinant is a nonsingular 
quadratic form on $M_2(F)$ and on $M_2^0(F)$, in both cases with square 
discriminant.  

\newcommand{\Tr}{\textup{Tr}}

\begin{Prop}  \label{SO3-4}
Define 
        $$ \rho_3\:SL_2(F) \Right5{} \Omega(M_2^0(F),\det) $$
and 
        $$ \rho_4\:SL_2(F)\times{}SL_2(F)
        \Right5{} \Omega(M_2(F),\det) $$
by setting
        $$ \rho_3(A)(X)=AXA^{-1} \qquad\textup{and}\qquad
        \rho_4(A,B)(X)=AXB^{-1}. $$
Then $\rho_3$ and $\rho_4$ are both epimorphisms, and lift to unique 
isomorphisms
        $$ SL_2(F) \RIGHT5{\widetilde{\rho}_3}{\cong}
        \Spin(M_2^0(F),\det) $$
and
        $$ SL_2(F)\times{}SL_2(F) \RIGHT5{\widetilde{\rho}_4}{\cong}
        \Spin(M_2(F),\det). $$
\end{Prop}

\begin{proof}  See \cite[pages 142, 199]{Taylor} for other ways of defining 
these isomorphisms.  By Lemma \ref{Galois-Spin}, it suffices to prove this 
(except for the uniqueness of the lifting) when $F$ is algebraically closed.  
In particular, 
	$$ \Omega(M_2^0(F),\det)=SO(M_2^0(F),\det)
	\quad\textup{and}\quad
	\Omega(M_2(F),\det)=SO(M_2(F),\det) $$ 
in this case.

For general $V$ and $\bb$, the group $SO(V,\bb)$ is generated by 
reflections fixing nonisotropic subspaces (ie, of nonvanishing 
discriminant) of codimension 2 (cf \cite[\S{}II.6(1)]{Dieudonne}).  Hence 
to see that $\rho_3$ and $\rho_4$ are surjective, it suffices to show that 
such elements lie in their images.  A codimension 2 reflection in 
$SO(M_2^0(F),\det)$ is of the form $R_X$ (the reflection fixing the line 
generated by $X$) for some $X\in{}M_2^0(F)$ which is nonisotropic (ie, 
$\det(X)\ne0$).  Since $F$ is algebraically closed, we can assume 
$X\in{}SL_2(F)$.  Then $X^2=-I$ (since $\Tr(X)=0$ and $\det(X)=1$), and 
$R_X=\rho_3(X)$ since it has order 2 and fixes $X$.  Thus $\rho_3$ is onto.

Similarly, any 2--dimensional nonisotropic subspace $W\subseteq{}V$ has an 
orthonormal basis $\{Y,Z\}$, and $ZY^{-1}$ and $Y^{-1}Z$ have trace zero 
(since they are orthogonal to the identity matrix) and determinant one.  
Hence their square is $-I$, and one repeats the above argument to show 
that $R_W=\rho_4(ZY^{-1},Y^{-1}Z)$.  So $\rho_4$ is onto.

The liftings $\widetilde{\rho}_m$ exist and are unique since $SL_2(F)$ is 
the universal central extension of $PSL_2(F)$ (or universal among central 
extensions by 2--groups if $F=\F_3$).
\end{proof}

We now restrict to the case $F=\F_q$ where $q$ is an odd prime power.  We 
refer to \cite[\S21]{Aschbacher} for a description of quadratic forms in 
this situation, and the notation for the associated orthogonal groups.  If 
$n$ is odd and $\bb$ is any nonsingular quadratic form on $\F_q^n$, then 
every nonsingular quadratic form is isomorphic to $u\bb$ for some 
$u\in\F_q^*$, and hence one can write 
$SO_n(q)=SO(\F_q^n,\bb)=SO(\F_q^n,u\bb)$ without ambiguity.  If $n$ is 
even, then there are exactly two isomorphism classes of quadratic forms on 
$\F_q^n$; and one writes $SO^+_n(q)=SO(\F_q^n,\bb)$ when $\bb$ is the 
hyperbolic form (equivalently, has discriminant $(-1)^{n/2}$ modulo 
squares), and $SO^-_n(q)=SO(\F_q^n,\bb)$ when $\bb$ is not hyperbolic 
(equivalently, has discriminant $(-1)^{n/2}{\cdot}u$ for $u\in{}\F_q^*$ 
not a square).  This notation extends in the obvious way to 
$\Omega^{\pm}_n(q)$ and $\Spin^{\pm}_n(q)$.

The following lemma does, in fact, hold for for orthogonal representations 
over arbitrary fields of characteristic $\ne2$.  But to simplify the proof 
(and since we were unable to find a reference), we state it only in the case 
of finite fields.

\begin{Lem} \label{2-irred}
Assume $F=\F_q$, where $q$ is a power of an odd prime.  Let $V$ be an 
$F$--vector space, and let $\bb$ be a nonsingular quadratic form on $V$.  
Let $P\le{}O(V,\bb)$ be a 2--subgroup which is orthogonally irreducible; 
ie, such that $V$ has no splitting as an orthogonal direct sum of 
nonzero $P$--invariant subspaces.  Then $\dim_F(V)$ is a power of 2; and if 
$\dim(V)>1$ then $\bb$ has square discriminant.
\end{Lem}

\begin{proof}  This means showing that each orthogonal group 
$O(\F_q{}^n,\bb)$, such that either $n$ is not a power of $2$, or 
$n=2^k\ge2$ and the quadratic form $\bb$ has nonsquare discriminant, 
contains some subgroup $O_{m}^{\pm}(q)\times{}O_{n-m}^{\pm}(q)$ (for 
$0<m<n$) of odd index.  We refer to the standard formulas for the orders 
of these groups (see \cite[p.165]{Taylor}):  if $\epsilon=\pm1$ then
	$$ |O_{2n}^\epsilon(q)| = 
	2q^{n(n-1)}(q^n-\epsilon)\prod_{i=1}^{n-1}(q^{2i}-1) 
	\quad\textup{and}\quad
	|O_{2n+1}(q)| = 2q^{n^2}\prod_{i=1}^{n}(q^{2i}-1). $$
We will also use repeatedly the fact that for all $0<i<2^k$ ($k\ge1$), the 
largest powers of $2$ dividing $(q^{2^k+i}-1)$ and $(q^i-1)$ are the same.  
In other words, $(q^{2^k+i}-1)/(q^i-1)$ is invertible in $\Z_{(2)}$.  

For any $n\ge1$, 
	$$ \frac{|O_{2n+1}(q)|}{|O_{2n}^\epsilon(q)|{\cdot}|O_1(q)|} = 
	q^n{\cdot}\frac{q^n+\epsilon}2 $$
is odd for an appropriate choice of $\epsilon$.  Thus, there are no 
irreducibles of odd dimension.  

Assume $n$ is not a power of $2$, and write $n=2^k+m$ where $0<m<2^k$ and 
$k\ge1$.  Then
	$$ \frac{|O_{2n}^\epsilon(q)|}{|O_{2^{k+1}}^+(q)|{\cdot}
	|O_{2m}^\epsilon(q)|} = q^{m2^{k+1}}{\cdot}
	\left(\prod_{i=1}^{m-1} \frac{q^{2(2^k+i)}-1}{q^{2i}-1}\right)
	{\cdot}\left(\frac{q^{2^k+m}-\epsilon}{q^m-\epsilon}\right)
	{\cdot}\left(\frac{q^{2^{k}}+1}{2}\right), $$
and each factor is invertible in $\Z_{(2)}$.  When $n=2m=2^k$ and $k\ge1$, 
then $O_{2n}^-(q)$ is the orthogonal group for the quadratic form with 
nonsquare discriminant, and
	$$ \frac{|O_{2n}^-(q)|}{|O_{2m}^+(q)|{\cdot}
	|O_{2m}^-(q)|} = q^{2m^2}{\cdot}
	\left(\prod_{i=1}^{m-1} \frac{q^{2(m+i)}-1}{q^{2i}-1}\right)
	{\cdot} \frac{q^{2m}+1}{2}, $$
and again each factor is invertible in $\Z_{(2)}$.  Finally, 
	$$ \frac{|O_2^\epsilon(q)|}{|O_1(q)|{\cdot}|O_1(q)|} = 
	\frac{q-\epsilon}2 $$
is odd whenever $q\equiv1$ (mod $4$) and $\epsilon=-1$, or $q\equiv3$ (mod 
$4$) and $\epsilon=+1$; and these are precisely the cases where the 
quadratic form on $\F_q{}^2$ has nonsquare discriminant.
\end{proof}

We must classify the conjugacy classes of those elementary abelian 
2--subgroups of $\Spin_7(q)$ which contain its center.  The following 
definition will be useful when doing this.

\begin{Defi}  \label{types}
Fix an odd prime power $q$.  Identify $SO_7(q)=SO(\F_q^7,\bb)$ and 
$\Spin_7(q)=\Spin(\F_q^7,\bb)$, where $\bb$ is a nonsingular quadratic 
form with square discriminant.  An elementary abelian 2--subgroup of 
$SO_7(q)$ or of $\Spin_7(q)$ will be called of \emph{type I} if its 
eigenspaces all have square discriminant (with respect to $\bb$), and of 
\emph{type II} otherwise.  Let $\cale_n$ be the set of elementary abelian 
2--subgroups in $\Spin_7(q)$ which contain $Z(\Spin_7(q))\cong{}C_2$ and 
have rank $n$.  Let $\cale_n^I$ and $\cale_n^{II}$ be the subsets of 
$\cale_n$ consisting of those subgroups of types I and II, respectively.
\end{Defi}

In the following two propositions, we collect together the information 
which will be needed about elementary abelian 2--subgroups of $\Spin_7(q)$. 
We fix $\Spin_7(q)=\Spin(V,\bb)$, where $V\cong\F_q^7$, and $\bb$ is a 
nonsingular quadratic form \emph{with square discriminant}.  Let 
$z\in{}Z(\Spin_7(q))$ be the generator.  For any subgroup $H\le\Spin_7(q)$ 
or any element $g\in\Spin_7(q)$, let $\widebar{H}$ and $\widebar{g}$ 
denote their images in $\Omega_7(q)\le{}SO_7(q)$.  For each elementary 
abelian 2--subgroup $E\le\Spin_7(q)$, and each character 
$\chi\in\Hom(\widebar{E},\{\pm1\})$, $V_\chi\subseteq{}V$ denotes the 
eigenspace of $\chi$ (and $V_1$ denotes the eigenspace of the trivial 
character).  Also (when $z\in{}E$), $\Aut(E,z)$ denotes the group of all 
automorphisms of $E$ which send $z$ to itself.

\newcommand{\discr}{\textup{discr}}

\begin{Prop} \label{elem.abel.} 
For any odd prime power $q$, the following table describes the numbers of 
$\Spin_7(q)$--conjugacy classes in each of the sets $\cale_n^I$ and 
$\cale_n^{II}$, the dimensions and discriminants of the eigenspaces of 
subgroups in these sets, and indicates in which cases 
$\Aut_{\Spin_7(q)}(E)$ contains all automorphisms which fix $z$.

\smallskip
\centerline{\upshape
\renewcommand{\arraystretch}{1.25}
\begin{tabular}{|l||c|c|c|c|c|}
\hline
Set of subgroups & $\cale_2^I$ & $\cale_3^I$ & $\cale_3^{II}$ & $\cale_4^I$ & 
$\cale_4^{II}$ \\ \hline\hline
Nr. conj. classes & 1 & 1 & 1 & 2 & 1 \\ \hline
$\dim(V_1)$ & 3 & \multicolumn{2}{c|}{1} & \multicolumn{2}{c|}{0} \\ \hline
$\dim(V_\chi)$, $\chi\ne1$ & 4 & \multicolumn{2}{c|}{2} & 
\multicolumn{2}{c|}{1} \\ \hline
$\discr(V_1,\bb)$ & square & square & nonsq. & --- & --- \\ \hline
$\discr(V_\chi,\bb)$, $\chi\ne1$ & square & square & nonsq. & square & 
both \\ \hline
$\Aut_{\Spin_7(q)}(E)=\Aut(E,z)$ & yes & yes & yes & yes & no \\ \hline
\end{tabular}
}

\smallskip\noindent There are no subgroups in $\cale_2$ of type II, and no 
subgroups of rank $\ge5$.  Furthermore, we have:
\begin{enumerate}  
\item For all $E\in\cale_4$, $C_{\Spin_7(q)}(E)=E$.  

\item If $E,E'\in\cale_4^I$, then $\widebar{E'}=g\widebar{E}g^{-1}$ for 
some $g\in{}SO_7(q)$, and $E$ and $E'$ are $\Spin_7(q)$--conjugate if and 
only if $g\in\Omega_7(q)$.  

\item If $E\in\cale_4^{II}$, then there is a unique element 
$1\ne\widebar{x}=\widebar{x}(E)\in\widebar{E}$ with the property that for 
$1\ne\chi\in\Hom(\widebar{E},\{\pm1\})$, $V_\chi$ has square discriminant 
if $\chi(\widebar{x})=1$ and nonsquare discriminant if 
$\chi(\widebar{x})=-1$.  Also, the image of $\Aut_{\Spin_7(q)}(E)$ in 
$\Aut(\widebar{E})$ is the group of all automorphisms which send 
$\widebar{x}$ to itself; and if $X\le{}E$ denotes the inverse image of 
$\gen{\widebar{x}}\le\widebar{E}$, then $\Aut_{\Spin_7(q)}(E)$ contains 
all automorphisms of $\widebar{E}$ which are the identity on $X$ and the 
identity modulo $\gen{z}$.

\item If $E\in\cale_3$, then $C_{\Spin_7(q)}(E)=A{\rtimes}C_2$, where $A$ 
is abelian and $C_2$ acts on $A$ by inversion.  If $E\in\cale_3^{II}$, 
then the Sylow 2--subgroups of $C_{\Spin_7(q)}(E)$ are elementary abelian 
of rank $4$ (and type II). 
\end{enumerate}
\end{Prop}

\begin{proof}  Write $\Spin=\Spin_7(q)$ for short.  Fix an elementary 
abelian subgroup $E\le\Spin$ such that $z\in{}E$.  

\smallskip

\noindent\textbf{Step 1}\qua 
We first show that $\rk(E)\le4$, and that the dimensions of the eigenspaces 
$V_\chi$ for $\chi\in\Hom(\widebar{E},\{\pm1\})$ are as described in the 
table.

By Lemma \ref{SO-invol}, every involution in $\widebar{E}$ has a 
4--dimensional $(-1)$--eigenspace.  In particular, if $\rk(E)=2$, 
($\widebar{E}\cong{}C_2$), then $\dim(V_\chi)=4$ for 
$1\ne\chi\in\Hom(\widebar{E},\{\pm1\})$, while $\dim(V_1)=3$.  

Now assume $\rk(E)=n$ for some $n>2$.  Assume we have shown, for all 
$E'\in\cale_{n-1}$, that the eigenspace of the trivial character of 
$\widebar{E}'$ is $r$--dimensional.  For each 
$1\ne\chi\in\Hom(\widebar{E},\{\pm1\})$, let $E_\chi\in\cale_{n-1}$ be the 
subgroup such that $\widebar{E}_\chi=\Ker(\chi)$; then $V_1\oplus{}V_\chi$ 
is the eigenspace of the trivial character of 
$\widebar{E}_\chi=\Ker(\chi)$, and thus $\dim(V_1)+\dim(V_\chi)=r$.  Hence 
all nontrivial characters of $E$ have eigenspaces of the same dimension.  
Since there are $2^{n-1}-1$ nontrivial characters of $\widebar{E}$, we have
$\dim(V_1)+(2^{n-1}-1)\dim(V_\chi)=7$, and these two equations completely 
determine $\dim(V_1)$ and $\dim(V_\chi)$.  Using this procedure, the 
dimensions of the eigenspaces are shown inductively to be 
equal to those given by the table.  Also, this shows that 
$\rk(E)\le4$, since otherwise $\rk(\widebar{E})\ge4$, so the $V_\chi$ for 
$\chi\ne1$ must be trivial (they cannot all have dimension $\ge1$), so $E$ 
acts on $V$ via the identity, which contradicts the assumption that 
$E\le\Spin_7(q)$.

\smallskip

\noindent\textbf{Step 2}\qua 
We next show that $\cale_2^{II}=\emptyset$, describe the discriminants of 
the eigenspaces of characters of $\widebar{E}$ for $E\in\cale_n$ (for 
all $n$), and show that subgroups of rank $4$ are self centralizing.  In 
particular, this proves (a) together with the first statement of (c).  

If $E\in\cale_2$, then $E=\gen{z,g}$ for some noncentral involution 
$g\in\Spin_7(q)$, and the eigenspaces of $\widebar{E}=\gen{\widebar{g}}$ 
have square discriminant by Lemma \ref{SO-invol}(a) (and since the ambient 
space $V$ has square discriminant by assumption).  Thus 
$\cale_2^{II}=\emptyset$.  

If $E\in\cale_3$, then the sum of any two eigenspaces of $\widebar{E}$ is 
an eigenspace of $\widebar{g}$ for some $g\in{}E{\sminus}\gen{z}$.  Hence 
the sum of any two eigenspaces of $\widebar{E}$ has square discriminant, so 
either all of the eigenspaces have square discriminant ($E\in\cale_3^I$), 
or all of the eigenspaces have nonsquare discriminant ($E\in\cale_3^{II}$).

Assume $\rk(E)=4$.  We have seen that all eigenspaces of $\widebar{E}$ are 
$1$--dimensional.  By Lemma \ref{SO-invol}(c), for each 
$a\in{}C_{\Spin_7(q)}(E)$, $\widebar{a}(V_\chi)=V_\chi$ for each 
$\chi\ne1$, and since $\dim(V_\chi)=1$ it must act on each $V_\chi$ via 
$\pm\Id$.  Thus $\widebar{a}\in\Omega_7(q)$ has order $2$; let $V_\pm$ be 
its eigenspaces.  Then $\dim(V_-)$ is even since $\det(\widebar{a})=1$, 
and $V_-$ has square discriminant by Lemma \ref{SO-invol}(a).  If 
$\dim(V_-)=4$, then $|a|=2$ (Lemma \ref{SO-invol}(b)), and hence $a\in{}E$ 
since otherwise $\gen{E,a}$ would have rank 5.  If $\dim(V_-)=2$, then 
$V_-$ is the sum of the eigenspaces of two distinct characters 
$\chi_1,\chi_2$ of $\widebar{E}$, there is some $g\in{}E$ such that 
$\chi_1(\widebar{g})\ne\chi_2(\widebar{g})$, hence
$\det(\widebar{g}|_{V_-})=\chi_1(\widebar{g})\chi_2(\widebar{g})=-1$, 
so $[g,a]=z$ by Lemma \ref{SO-invol}(c), and this contradicts the 
assumption that $[a,E]=1$.  If $\dim(V_-)=6$, then $V_-$ is the sum of the 
eigenspaces of all but one of the nontrivial characters of $\widebar{E}$, 
and this gives a similar contradiction to the assumption $[a,E]=1$.  Thus, 
$C_{\Spin_7(q)}(E)=E$.

Now assume that $E\in\cale_4^{II}$, and let $\widebar{x}\in{}O_7(q)$ be 
the element which acts via $-\Id$ on eigenspaces with nonsquare 
discriminant, and via the identity on those with square discriminant.  
Since $\bb$ has square discriminant on $V$, the number of eigenspaces of 
$\widebar{E}$ on which the discriminant is
nonsquare is even, so $\widebar{x}\in\Omega_7(q)$ by Lemma 
\ref{SO-invol}(a), and lifts to an element $x\in\Spin_7(q)$.  Also, for 
each $g\in{}E$, the $(-1)$--eigenspace of $\widebar{g}$ has square 
discriminant (Lemma \ref{SO-invol}(a) again), hence contains an even 
number of eigenspaces of $\widebar{E}$ of nonsquare discriminant, and by 
Lemma \ref{SO-invol}(c) this shows that $[g,x]=1$.  Thus 
$x\in{}C_{\Spin_7(q)}(E)=E$, and this proves the first statement in (c).

\smallskip

\noindent\textbf{Step 3}\qua We next check the numbers of conjugacy classes 
of subgroups in each of the sets $\cale_n^{I}$ and $\cale_n^{II}$, and 
describe $\Aut_{\Spin}(E)$ in each case.  This finishes the proof of (b) 
and (c), and of all points in the above table.

{}From the above description, we see immediately that if $E$ and $E'$ have 
the same rank and type, then any isomorphism 
$\alpha\in\Iso(\widebar{E},\widebar{E}')$, such that 
$\alpha(\widebar{x}(E))=\widebar{x}(E')$ if $E,E'\in\cale_4^{II}$, has the 
property that for all $\chi\in\Hom(\widebar{E}',\{\pm1\})$, $V_\chi$ and 
$V_{\chi\circ\alpha}$ have the same dimension and the same discriminant 
(modulo squares).  Hence for any such $\alpha$, there is an element 
$g\in{}O_7(q)$ such that $g(V_{\chi\circ\alpha})=V_\chi$ for all $\chi$; 
and $\alpha=c_g\in\Iso(\widebar{E},\widebar{E}')$ for such $g$.  Upon 
replacing $g$ by $-g$ if necessary, we can assume that $g\in{}SO_7(q)$.  
This shows that
	\beq \textup{$E,E'$ have the same rank and type \,$\Longrightarrow$\,
	$\widebar{E}$ and $\widebar{E}'$ are $SO_7(q)$--conjugate} 
	\tag{1} \eeq
and also that
	\beq \Aut_{SO_7(q)}(\widebar{E}) = \begin{cases}
	\Aut(\widebar{E}) & \textup{if $E\notin\cale_4^{II}$} \\ 
	\Aut(\widebar{E},\widebar{x}(E)) & \textup{if $E\in\cale_4^{II}$.}
	\end{cases} \tag{2} \eeq

We next claim that
	\beq E\notin\cale_4^I \,\Longrightarrow\, 
	\exists \gamma\in{}SO_7(q){\sminus}\Omega_7(q) \textup{ such that 
	}[\gamma,\widebar{E}]=1\,. \tag{3} \eeq
To prove this, choose 1--dimensional nonisotropic summands 
$W\subseteq{}V_\chi$ and $W'\subseteq{}V_\psi$, where $\chi,\psi$ are two 
distinct characters of $\widebar{E}$, and where $W$ has square 
discriminant and $W'$ has nonsquare discriminant.  Let
$\gamma\in{}SO_7(q)$ be the involution with $(-1)$--eigenspace 
$W\oplus{}W'$.  Then $[\gamma,\widebar{E}]=1$, since $\gamma$ sends each
eigenspace of $\widebar{E}$ to itself, and $\gamma\notin\Omega_7(q)$ since 
its $(-1)$--eigenspace has nonsquare discriminant (Lemma \ref{SO-invol}(a)).

If $E$ has rank 4 and type I, then $\Aut(\widebar{E})\cong{}GL_3(\F_2)$ is 
simple, and in particular has no subgroup of index 2.  Hence by (2), each 
element of $\Aut(\widebar{E})$ is induced by conjugation by an element of 
$\Omega_7(q)$.  Also, if $g\in{}SO_7(q)$ centralizes $\widebar{E}$, then 
$g(V_\chi)=V_\chi$ for all $\chi\in\Hom(\widebar{E},\{\pm1\})$, 
$g$ acts via $-\Id$ on an even number of eigenspaces (since it has 
determinant $+1$), and hence $g\in\Omega_7(q)$ by Lemma 
\ref{SO-invol}(a).  Thus 
	\beq E\in\cale_4^I \,\Longrightarrow\, N_{SO_7(q)}(\widebar{E}) 
	\le \Omega_7(q) \tag{4} \eeq

If $E\notin\cale_4^I$, then by (3), for any $g\in{}SO_7(q)$, there is 
$\gamma\in{}SO_7(q){\sminus}\Omega_7(q)$ such that 
$c_g|_E=c_{g\gamma}|_E$, and either $g$ or $g\gamma$ lies in 
$\Omega_7(q)$.  Thus $\Iso_{SO_7(q)}(\widebar{E},\widebar{E}') = 
\Iso_{\Omega_7(q)}(\widebar{E},\widebar{E}')$ for any $E'$.  
Together with (1), this shows that $E$ is $\Spin$--conjugate to all other 
subgroups of the same rank and type, and together with (2) it shows that
	\beq \Im\bigl[\Aut_{\Spin}(E) \Right1{} \Aut(\widebar{E})\bigr] =
	\begin{cases}  
	\Aut(\widebar{E}) & \textup{if $E\notin\cale_4^{II}$} \\ 
	\Aut(\widebar{E},\widebar{x}(E)) & \textup{if $E\in\cale_4^{II}$.}
	\end{cases} \tag{5} \eeq
If $E\in\cale_4^I$, then by (4) and (2), $\Aut_{\Omega_7(q)}(\widebar{E})= 
\Aut_{SO_7(q)}(\widebar{E})=\Aut(\widebar{E})$, and so (5) also holds in 
this case.  Furthermore, for any $g\in{}SO_7(q){\sminus}\Omega_7(q)$, 
$\widebar{E}$ and $g\widebar{E}g^{-1}$ are representatives for two 
distinct $\Omega_7(q)$--conjugacy classes --- since by (4), no element of 
the coset $g{\cdot}\Omega_7(q)$ normalizes $\widebar{E}$.  

We have now determined in all cases the number of conjugacy classes of 
subgroups of a given rank and type, and the image of $\Aut_{\Spin}(E)$ in 
$\Aut(\widebar{E})$.  We next claim that if $\rk(E)<4$ or $E\in\cale_4^{I}$, 
then
	\beq E\notin\cale_4^{II} \,\Longrightarrow\, 
	\Aut_{\Spin}(E) \ge \bigl\{\alpha\in\Aut(E)\,\big|\,
	\alpha(z)=z,\ 
	\textup{$\alpha\equiv\Id$ (mod $\gen{z}$) } \bigr\} \,.
	\tag{6} \eeq
Together with (5), this will finish the proof that $\Aut_{\Spin}(E)$ is 
the group of all automorphisms of $E$ which send $z$ to itself.  We also 
claim that
	\beq E\in\cale_4^{II} \,\Longrightarrow\, 
	\Aut_{\Spin}(E) \ge \bigl\{\alpha\in\Aut(E)\,\big|\,
	\alpha|_X=\Id_X,\ 
	\textup{$\alpha\equiv\Id$ (mod $\gen{z}$)} \bigr\} \,,
	\tag{7} \eeq
where $X\le{}E$ denotes the inverse image of 
$\gen{\widebar{x}(E)}\le\widebar{E}$, and this will finish the proof of 
(c).

We prove (6) and (7) together.  Fix $\alpha\in\Aut(E)$ ($\alpha\ne\Id$) 
which sends $z$ to itself, induces the identity on $\widebar{E}$, and such 
that $\alpha|_X=\Id_X$ if $E\in\cale_4^{II}$.  Then there is 
$1\ne\chi\in\Hom(\widebar{E},\{\pm1\})$ such that $\alpha(g)=g$ when 
$\chi(\widebar{g})=1$ and $\alpha(g)=zg$ when $\chi(\widebar{g})=-1$.  
Choose any character $\psi$ such that $V_\psi\ne0$ and $V_{\psi\chi}\ne0$, 
and let $W\subseteq{}V_\psi$ and $W'\subseteq{}V_{\psi\chi}$ be 
1--dimensional nonisotropic subspaces with the same discriminant (this is 
possible when $E\in\cale_4^{II}$ since $\widebar{x}(E)\in\Ker(\chi)$).  
Let $\widebar{g}\in{}O_7(q)$ be the involution whose $(-1)$--eigenspace is 
$W\oplus{}W'$.  Then $\widebar{g}\in\Omega_7(q)$ by Lemma 
\ref{SO-invol}(a), so $\widebar{g}$ lifts to $g\in\Spin_7(q)$, and using 
Lemma \ref{SO-invol}(c) one sees that $c_g=\alpha$.  

\smallskip

\noindent\textbf{Step 4}\qua It remains to prove (d). Assume $E\in\cale_3$.  
Let $1=\chi_1,\chi_2,\chi_3,\chi_4$ be the four characters of 
$\widebar{E}$, and set $V_i=V_{\chi_i}$. Then $\dim(V_1)=1$, $\dim(V_i)=2$ 
for $i=2,3,4$, and the $V_i$ either all have square discriminant or all 
have nonsquare discriminant.  For each $g\in{}C_{\Spin}(E)$, we can write 
$\widebar{g}=\bigoplus_{i=1}^4g_i$, where 
$g_i\in{}O(V_i,\bb_i)$.  For each pair of distinct indices 
$i,j\in\{2,3,4\}$, there is some $g\in{}E$ whose $(-1)$--eigenspace is 
$V_i\oplus{}V_j$, and hence $\det(g_i\oplus{}g_j)=1$ by Lemma 
\ref{SO-invol}(c).  This shows that the $g_i$ all have the same 
determinant.  Let $A\le{}C_{\Spin}(E)$ be the subgroup of index 2 
consisting of those $g$ such that $\det(g_i)=1$ for all $i$.

Now, $SO_1(\fqbar)=1$, while $SO_2(\fqbar)\cong\fqbar^*$ is the group of 
diagonal matrices of the form $\diag(u,u^{-1})$ with respect to a 
hyperbolic basis of $\fqbar{}^2$.  Thus $A$ is contained in a central 
extension of $C_2$ by $(\fqbar^*)^3$, and any such extension is abelian 
since $H_2((\fqbar^*)^3)=0$.  Hence $A$ is abelian.  The groups 
$O_2^\pm(q)$ are all dihedral (see \cite[Theorem 11.4]{Taylor}).  Hence 
for any $g\in{}C_{\Spin}(E){\sminus}A$, $\widebar{g}$ has order 
$2$ and $(-1)$--eigenspace of dimension 4 (its intersection with each $V_i$ 
is 1--dimensional), and hence $|g|=2$ by Lemma \ref{SO-invol}(b).  
Thus all elements of $C_{\Spin}(E){\sminus}A$ have order 2, so the 
centralizer must be a semidirect product of $A$ with a group of order 2 
which acts on it by inversion.

Now assume that $E\in\cale_3^{II}$; ie, that the $V_i$ all have 
nonsquare discriminant.  Then for $i=2,3,4$, $SO(V_i,\bb_i)$ has order 
$q\pm1$, whichever is not a multiple of 4 (see \cite[Theorem 11.4]{Taylor} 
again).  Thus if $g\in{}A\le{}C_{\Spin}(E)$ has 2--power order, then 
$g_i=\pm{}I$ for each $i$, the number of $i$ for which $g_i=\Id$ 
is even (since the $(-1)$--eigenspace of $\widebar{g}$ has square 
discriminant), and hence $g\in{}E$.  In other words,
$E\in\Syl_2(A)$.  A Sylow 2--subgroup of $C_{\Spin}(E)$ is thus generated by 
$E$ together with an element of order 2 which acts on $E$ by inversion; 
this is an elementary abelian subgroup of rank $4$, and is necessarily of 
type II.
\end{proof}


We also need some more precise information about the subgroups of 
$\Spin_7(q)$ of rank 4 and type II. Let $\psi^q\in\Aut(\Spin_7(\fqbar))$ 
denote the automorphism induced by the field automorphism 
$(x\mapsto{}x^q$).  By Lemma \ref{Galois-Spin}, $\Spin_7(q)$ is precisely 
the subgroup of elements fixed by $\psi^q$.  

\begin{Prop} \label{E4II}
Fix an odd prime power $q$, and let $z\in{}Z(\Spin_7(q))$ be the central 
involution.  Let $\calc$ and $\calc'$ denote the two conjugacy classes of 
subgroups $E\le\Spin_7(q)$ of rank 4 and type I.  Then the following hold.
\begin{enumerate}  
\item For each $E\in\cale_4$, there is an element $a\in\Spin_7(\fqbar)$ 
such that $aEa^{-1}\in\calc$.  For any such $a$, if we set
	$$ x_\calc(E) \defeq a^{-1}\psi^q(a), $$
then $x_\calc(E)\in{}E$ and is independent of the choice of $a$.  

\item $E\in\calc$ if and only if $x_{\calc}(E)=1$, and $E\in\calc'$ if and 
only if $x_{\calc}(E)=z$. 
 
\item Assume $E\in\cale_4^{II}$, and set $\tau(E)=\gen{z,x_{\calc}(E)}$.  
Then $\rk(\tau(E))=2$, and
        $$ \Aut_{\Spin_7(q)}(E) = \bigl\{\alpha\in\Aut(E) \,\big|\, 
        \alpha|_{\tau(E)}=\Id \bigr\}. $$
The four eigenspaces of $\widebar{E}$ contained in the $(-1)$--eigenspace 
of $\widebar{x_\calc(E)}$ all have nonsquare discriminant, and the other three 
eigenspaces all have square discriminant.
\end{enumerate}
\end{Prop}

\begin{proof}  
\textbf{(a) }  For all $E\in\cale_4$, $E$ has type I as a subgroup 
of $\Spin_7(q^2)$ since all elements of $\F_q$ are squares in $\F_{q^2}$.  
Hence by Proposition \ref{elem.abel.}(b), for all $E'\in\calc$, there is 
$\widebar{a}\in{}SO_7(q^2)\le\Omega_7(q^4)$ such that 
$\widebar{a}\widebar{E}\widebar{a}^{-1}=\widebar{E}'$.  Upon 
lifting $\widebar{a}$ to $a\in\Spin_7(q^4)$, this proves that there is 
$a\in\Spin_7(\fqbar)$ such that $aEa^{-1}\in\calc$.  

Fix any such $a$, and set
	$$ x = x_{\calc}(E)=a^{-1}\psi^q(a). $$
For all $g\in{}E$, $\psi^q(g)=g$ and $\psi^q(aga^{-1})=aga^{-1}$ since 
$E,aEa^{-1}\le\Spin_7(q)$, and hence
	$$ aga^{-1} = \psi^q(a){\cdot}g{\cdot}\psi^q(a^{-1}) = 
	a(xgx^{-1})a^{-1}. $$
Thus, $x\in{}C_{\Spin_7(\fqbar)}(E)$, and so $x\in{}E$ since it is self 
centralizing in each $\Spin_7(q^k)$ (Proposition \ref{elem.abel.}(a)).

We next check that $x_{\calc}(E)$ is independent of the choice of $a$.  
Assume $a,b\in\Spin_7(\fqbar)$ are such that $aEa^{-1}\in\calc$ and 
$bEb^{-1}\in\calc$.  Then by Proposition \ref{elem.abel.}(b), there is 
$g\in\Spin_7(q)$ such that $gbE(gb)^{-1}=aEa^{-1}$.  Set 
$E'=aEa^{-1}\in\calc$, then $gba^{-1}\in{}N_{\Spin_7(\fqbar)}(E')$.  
Furthermore, since $\Aut_{\Spin_7(q)}(E')$ contains all automorphisms 
which send $z$ to itself, and since $E'$ is self centralizing in each of 
the groups $\Spin_7(q^k)$ (both by Proposition \ref{elem.abel.} again), we 
see that $N_{\Spin_7(\fqbar)}(E')$ is contained in $\Spin_7(q)$.  Thus, 
$ba^{-1}\in\Spin_7(q)$, so $\psi^q(ba^{-1})=ba^{-1}$; and this proves that 
$x_{\calc}(E)=a^{-1}\psi^q(a)=b^{-1}\psi^q(b)$ is independent of the 
choice of $a$.  

\smallskip

\noindent\textbf{(b) }
If $E\in\calc$, then we can choose $a=1$, and so $x_{\calc}(E)=1$.  

If $E\in\calc'$, then by Proposition \ref{elem.abel.}(b), there is 
$a\in\Spin_7(q^2)$ such that $\widebar{a}\in{}SO_7(q){\sminus}\Omega_7(q)$ 
and $aEa^{-1}\in\calc$.  Then $\psi^q(a)\ne{}a$ since $a\notin\Spin_7(q)$ 
(Proposition \ref{Galois-Spin}), but $\psi^q(\widebar{a})=\widebar{a}$ 
since $\widebar{a}\in{}SO_7(q)$.  Thus, $x_{\calc}(E)=a^{-1}\psi^q(a)=z$ in 
this case.  

We have now shown that $x_\calc(E)\in\gen{z}$ if $E$ has type I, and it 
remains to prove the converse.  Fix $a\in\Spin_7(\fqbar)$ such that 
$aEa^{-1}\in\calc$.  If $x_\calc(E)\in\gen{z}$, then 
$\psi^q(a)\in\{a,za\}$, so $\psi^q(\widebar{a})=\widebar{a}$, and hence 
$\widebar{a}\in{}SO_7(q)$.  Conjugation by an element of $SO_7(q)$ sends 
eigenspaces with square discriminant to eigenspaces with square 
discriminant, so all eigenspaces of $E$ must have square discriminant 
since all eigenspaces of $aEa^{-1}$ do.  Hence $E$ has type I.

\smallskip

\noindent\textbf{(c) } Now write $\Spin=\Spin_7(q)$ for short.  Assume 
$E\in\cale_4^{II}$, and set $x=x_{\calc}(E)$ and $\tau(E)=\gen{z,x}$.  Then
$x\notin\gen{z}$ by (b), and thus $\tau(E)$ has rank 2.  

By (a) (the uniqueness of $x$ having the given properties), each element 
of $\Aut_{\Spin}(E)$ restricts to the identity on $\tau(E)$.  We have 
already seen (Proposition \ref{elem.abel.}(c)) that there is an element 
$\widebar{x}(E)\in\widebar{E}$ such that the image in $\Aut(\widebar{E})$ 
of $\Aut_{\Spin}(E)$ is the group of automorphisms which fix 
$\widebar{x}(E)$, and this shows that $\widebar{x}(E)=\widebar{x}$:  the 
image in $\widebar{E}$ of $x$.  Since we already showed (Proposition 
\ref{elem.abel.}(c) again) that $\Aut_{\Spin}(E)$ contains all 
automorphisms which are the identity on $\tau(E)$ and the identity modulo 
$\gen{z}$, this finishes the proof that $\Aut_{\Spin}(E)$ is the group of 
all automorphisms which are the identity on $\tau(E)$.  The last statement 
(about the discriminants of the eigenspaces) follows directly from the 
first statement of Proposition \ref{elem.abel.}(c).
\end{proof}


Throughout the rest of the section, we collect some more technical results 
which will be needed in Sections 2 and 4.  
 
\begin{Lem}  \label{liftAut}
Fix $k\ge2$.  Let $A=e_{13}(2^{k-1})\in{}GL_3(\Z/2^k)$ be the elementary 
matrix which has off diagonal entry $2^{k-1}$ in position $(1,3)$.  
Let $T_1$ and $T_2$ be the two maximal parabolic subgroups of $GL_3(2)$:
        $$ T_1=GL^1_2(\Z/2)=\bigl\{(a_{ij})\in{}GL_3(2)\,|\,
        a_{21}=a_{31}=0 \bigr\} $$
and
        $$ T_2=GL^2_1(\Z/2)=\bigl\{(a_{ij})\in{}GL_3(2)\,|\,
        a_{31}=a_{32}=0 \bigr\}. $$
Set $T_0=T_1\cap{}T_2$: the group of upper triangular matrices in $GL_3(2)$.  
Assume that 
        $$ \mu_i\: T_i \Right4{} SL_3(\Z/2^k) $$
are lifts of the inclusions (for $i=1,2$) such that 
$\mu_1|_{T_0}=\mu_2|_{T_0}$.  Then 
there is a homomorphism 
        $$ \mu\:GL_3(2)\Right2{}SL_3(\Z/2^k) $$ 
such that $\mu|_{T_1}=\mu_1$, and either $\mu|_{T_2}=\mu_2$, or 
$\mu|_{T_2}=c_A\circ\mu_2$.  
\end{Lem}

\begin{proof}  We first claim that any two liftings
$\sigma,\sigma'\:T_2\Right1{}SL_3(\Z/2^k)$ are conjugate by an element of 
$SL_3(\Z/2^k)$.  This clearly holds when $k=1$, and so we can assume 
inductively that $\sigma\equiv\sigma'$ (mod $2^{k-1}$).  Let $M_3^0(\F_2)$ 
be the group of $3\times3$ matrices of trace zero, and define 
$\rho\:T_2\Right1{}M_3^0(\F_2)$ via the formula
        $$ \sigma'(B)=(I+2^{k-1}\rho(B)){\cdot}\sigma(B) $$
for $B\in{}T_2$.  Then $\rho$ is a 1--cocycle.  Also, 
$H^1(T_2;M_3^0(\F_2))=0$ by \cite[Lemma 4.3]{DW:DI4} (the module is 
$\F_2[T_2]$--projective), so $\rho$ is the coboundary of some 
$X\in{}M_3^0(\F_2)$, and $\sigma$ and $\sigma'$ differ by conjugation by 
$I+2^{k-1}X$.  

By \cite[Theorem 4.1]{DW:DI4}, there exists a section $\mu$ defined on 
$GL_3(2)$ such that $\mu|_{T_1}=\mu_1$.  Let $B\in{}SL_3(\Z/2^k)$ be such 
that $\mu|_{T_2}=c_B\circ\mu_2$.  Since $\mu|_{T_0}=\mu_2|_{T_0}$, $B$ 
must commute with all elements in $\mu(T_0)$, and one easily checks that 
the only such elements are $A=e_{13}(2^{k-1})$ and the identity.
\end{proof}

Recall that a $p$--subgroup $P$ of a finite group $G$ is \emph{$p$--radical} 
if $N_G(P)/P$ is $p$--reduced; ie, if $O_p(N_G(P)/P)=1$.  (Here, $O_p(-)$ 
denotes the largest normal $p$--subgroup.)  We say here that $P$ is 
\emph{$\calf_p(G)$--radical} if $\Out_G(P)$ ($=\Out_{\calf_p(G)}(P)$) is 
$p$--reduced.  In Section 4, some information will be needed involving the 
$\calf_2(\Spin_7(q))$--radical subgroups of $\Spin_7(q)$ which are also 
$2$--centric.  We first note the following general result.

\begin{Lem} \label{str-p-rad}
Fix a finite group $G$ and a prime $p$.  Then the following hold for any 
$p$--subgroup $P\le{}G$ which is $p$--centric and $\calf_p(G)$--radical.
\begin{enumerate}  
\item If $G=G_1\times{}G_2$, then $P=P_1\times{}P_2$, where $P_i$ is 
$p$--centric in $G_i$ and $\calf_p(G_i)$--radical.

\item If $P\le{}H\nsg{}G$, then $P$ is $p$--centric in $H$ and 
$\calf_p(H)$--radical.

\item If $H\nsg{}G$ has $p$--power index, then $P\cap{}H$ is $p$--centric in 
$H$ and $\calf_p(H)$--radical.

\item If $G\nsg\widebar{G}$ has $p$--power index, then $P=G\cap\widebar{P}$ 
for some $\widebar{P}\le\widebar{G}$ which is $p$--centric in $\widebar{G}$ and 
$\calf_p(\widebar{G})$--radical.

\item If $Q\nsg{}G$ is a central $p$--subgroup, then $Q\le{}P$, and $P/Q$ 
is $p$--centric in $G/Q$ and $\calf_p(G/Q)$--radical.

\item If $\widetilde{G}\Onto2{\alpha}G$ is an epimorphism such that
$\Ker(\alpha)\le{}Z(\widetilde{G})$, then $\alpha^{-1}(P)$ is $p$--centric 
in $\widetilde{G}$ and $\calf_p(\widetilde{G})$--radical. 
\end{enumerate}
\end{Lem}

\begin{proof}  Point (a) follows from \cite[Proposition 1.6(ii)]{JMO}:  
$P=P_1\times{}P_2$ for $P_i\le{}G_i$ since $P$ is $p$--radical, and $P_i$ 
must be $p$--centric in $G_i$ and $\calf_p(G_i)$--radical since 
	$$ C_G(P) = C_{G_1}(P_1)\times{}C_{G_2}(P_2)
	\qquad\textup{and}\qquad
	\Out_G(P)\cong\Out_{P_1}(G_1)\times\Out_{P_2}(G_2). $$
Point (b) holds since $C_H(P)\le{}C_G(P)$ and 
$O_p(\Out_H(P))\le{}O_p(\Out_G(P))$.  

It remains to prove the other four points.

\smallskip

\noindent\textbf{(e) } Fix a central $p$--subgroup $Q\le{}Z(G)$.  Then 
$P\ge{}Q$, since otherwise $1\ne{}N_{QP}(P)/P\le{}O_p(N_G(P)/P)$.  Also, 
$P/Q$ is $p$--centric in $G/Q$, since otherwise there would be 
$x\in{}G{\sminus}P$ of $p$--power order such that
	$$ 1\ne [c_x] \in 
	\Ker\bigl[\Out_G(P)\Right2{}\Out_{G/Q}(P/Q)\times\Out_G(Q)\bigr]
	\le O_p(\Out_G(P)). $$
It remains only to prove that $P/Q$ is $\calf_p(G/Q)$--radical, and to do 
this it suffices to show that
	$$ \Out_{G/Q}(P/Q) \cong \Out_G(P). $$
Equivalently, since $P/Q$ and $P$ are $p$--centric, we must show that
	$$ \frac{N_{G/Q}(P/Q)}{C'_{G/Q}(P/Q)\times{}P/Q} \cong
	\frac{N_G(P)}{C'_G(P)\times{}P}; $$
and this is clear once we have shown that
	$$ C'_{G/Q}(P/Q) \cong C'_G(P). $$
Any $\widebar{x}\in{}C'_{G/Q}(P/Q)$ lifts to an element $x\in{}G$ of order 
prime to $p$, whose conjugation action on $P$ induces the identity on $Q$ 
and on $P/Q$.  By \cite[Corollary 5.3.3]{Gorenstein}, all such 
automorphisms of $P$ have $p$--power order, and thus $x$ centralizes $P$.  
Since $Q$ is a $p$--group and $C'_{G/Q}(P/Q)$ has order prime to $p$, this 
shows that the projection modulo $Q$ sends $C'_{G/Q}(P/Q)$ isomorphically 
to $C'_G(P)$.

\smallskip

\noindent\textbf{(f) } 
Let $\widetilde{G}\Onto2{\alpha}G$ be an epimorphism 
whose kernel is central.  
Clearly, $\alpha^{-1}P$ is $p$--centric in $\widetilde{G}$.  
It remains 
only to prove that $\alpha^{-1}P$ is 
$\calf_p(\widetilde{G})$--radical, and to do this it suffices to show that
	$$ \Out_{\widetilde{G}}(\alpha^{-1}P) \cong \Out_G(P). $$
Equivalently, since $P$ and $\alpha^{-1}(P)$ are $p$--centric, we 
must show that
	$$ \frac{N_{\widetilde{G}}(\alpha^{-1}P)}
	{C'_{\widetilde{G}}(\alpha^{-1}P)\times\alpha^{-1}P}
	\cong \frac{N_G(P)}{C'_G(P)\times{}P}; $$
and this is clear once we have shown that
	$$ C'_{\widetilde{G}}(\alpha^{-1}P) \cong C'_G(P). $$
This follows by exactly the same argument as in the proof of (e).

\smallskip

\noindent\textbf{(c) } Set $P'=P\cap{}H$ for short.  Let
	\[ N_H(P')\Onto5{\pi} \Out_H(P') \cong
	N_H(P')/(C_H(P'){\cdot}P') \] 
be the natural projection, and set 
	$$ K=\pi^{-1}(O_p(\Out_H(P'))) \le 
	N_H(P'). $$
Then $K\ge O_p(N_H(P'))$ is an extension of $C_H(P'){\cdot}P'$ by 
$O_p(\Out_H(P'))$.  It suffices to show that $p\nmid[K{:}P']$, 
since this implies that $O_p(\Out_H(P'))=1$ (ie, $P'$ is 
$\calf_p(H)$--radical), and that any Sylow $p$--subgroup of $C_H(P')$ is 
contained in $P'$ (hence $P'$ is $p$--centric in $H$).  

Assume otherwise:  that $p\big|[K{:}P']$.  Note first that $P'\nsg N_G(P)$, 
and that $N_G(P)\le N_G(K)$; ie, $N_G(P)$ normalizes 
$P'$ and $K$. The first statement is obvious, and the second is 
verified by observing directly that $N_G(P)$ normalizes $N_H(P')$ and 
$C_H(P')$. Thus the action of $N_G(P)$ on $K$ induces an action of 
$N_G(P)$, and in particular of $P$, on $K/P'$. Let $K_0/P'$ 
denote the fixed subgroup of this action of $P$. Since $p\big|[K{:}P']$ by 
assumption, and since $P$ is a $p$--group, $p\big||K_0/P'|$. 
A straightforward check also shows that $K_0\nsg N_G(P)$, and therefore 
that $PK_0\nsg N_G(P)$. Also, since $P'\le{}K_0\le{}H$,
	$$ PK_0/P \cong K_0/(P\cap{}K_0)=K_0/P' $$
is a normal subgroup of $N_G(P)/P$ of order a multiple of $p$.  Since $P$ 
is $p$--centric in $G$ by assumption, 
	$$ \Out_G(P) = N_G(P)/(C_G(P){\cdot}P) = N_G(P)/(C_G'(P)\times{}P),
	$$
and hence the image of $PK_0/P$ in $\Out_G(P)$ is a normal subgroup which 
also has order a multiple of $p$.  

By definition of $K$ as an extension of $C_H(P'){\cdot}P'$ by 
a $p$--group, if $x\in K$ has order prime to $p$, then $x\in C_H(P')$. 
Hence if $x\in K_0$ has order prime to $p$, then for every $z\in P$, 
$[x,z]\in P'$, so $x$ acts trivially on $P/P'$. Since $x$ 
also centralizes $P'$, it follows that $x$ centralizes $P$.  This 
shows that the image of $PK_0/P$ in $\Out_G(P)$ is a $p$--group, thus a 
nontrivial normal $p$--subgroup of $\Out_G(P)$, and this contradicts the 
original assumption that $P$ is $\calf_p(G)$--radical.

\smallskip

\noindent\textbf{(d) } 
Let $G\nsg\widebar{G}$ be a normal subgroup of $p$--power index and let 
$P\le G$ be a $p$--centric and $\calf_p(G)$--radical subgroup. Let
	\[ N_{\widebar{G}}(P) \Onto5{\pi} \Out_{\widebar{G}}(P) \cong
	N_{\widebar{G}}(P)\big/\bigl(C_{\widebar{G}}(P){\cdot}P\bigr) \]
be the natural surjection, and set 
	$$ K = \pi^{-1}\bigl(O_p(\Out_{\widebar{G}}(P))\bigr)
	\le N_{\widebar{G}}(P). $$
Then $K$ is an extension of $C_{\widebar{G}}(P){\cdot}P$ by 
$O_p(\Out_{\widebar{G}}(P))$.  Fix any $\widebar{P}\in \Syl_p(K)$.  We 
will show that $\widebar{P}\cap{}G=P$, and that $\widebar{P}$ is 
$p$--centric in $\widebar{G}$ and $\calf_p(\widebar{G})$--radical.

For each $x\in{}K\cap{}G\le N_G(P)$, 
	\[ \pi(x)\in O_p(\Out_{\widebar{G}}(P))\cap \Out_G(P) 
	\le O_p(\Out_G(P))= 1. \]
Hence
	\[ x\in\Ker\bigl[ N_G(P) \Right3{} 
	\Out_{\widebar{G}}(P) \bigr] =
	(C_{\widebar{G}}(P){\cdot}P)\cap G = 
	C_G(P)\cdot P \cong C_G'(P)\times P, \]
where $C_G'(P)\le C_G(P)$ is of order prime to $p$. Since the opposite 
inclusion is obvious, this shows that $K\cap G = C'_G(P)\times{}P$, and 
hence (since $\widebar{P}\in\Syl_p(K)$) that $\widebar{P}\cap G=P$. 

Next, note that $(K\cap G)\nsg K$ and $K/(K\cap G)\le \widebar{G}/G$, and 
hence $K/C'_G(P)$ has $p$--power order. Since $\widebar{P}\in\Syl_p(K)$, 
$\widebar{P}$ is an extension of $P$ by $K/(K\cap G)$, and 
$N_K(\widebar{P})$ is an extension of a subgroup of $(K\cap 
G)=(C'_G(P)\times{}P)$ by $K/(K\cap G)$.  Also, an element $x\in{}C'_G(P)$ 
normalizes $\widebar{P}$ if and only if 
$[x,\widebar{P}]\in\widebar{P}\cap{}C'_G(P)=1$.  Hence
	\beq N_K(\widebar{P}) = C_K(\widebar{P}){\cdot}\widebar{P} = 
	C'_G(\widebar{P})\times \widebar{P}, \tag{1} \eeq
where $C'_G(\widebar{P})=C'_G(P)\cap{}C_G(\widebar{P})$ has order prime to 
$p$ and is normal in $N_K(\widebar{P})$.  Since 
$C_{\widebar{G}}(\widebar{P})\le 
C_{\widebar{G}}(P)\le K$, (1) shows that 
$C_{\widebar{G}}(\widebar{P})\le{}C'_G(\widebar{P})\times\widebar{P}$, and 
hence that $\widebar{P}$ is $p$--centric in $\widebar{G}$.

It remains to show that $\widebar{P}$ is $\calf_p(\widebar{G})$--radical. Note 
first that $K\nsg N_{\widebar{G}}(P)$ by construction, so for any $x\in 
N_{\widebar{G}}(P)$, $x\widebar{P}x^{-1}\in\Syl_p(K)$. Since $K$ is an 
extension of $C_G'(P)\times{}P$ by the $p$--group $K/(K\cap G)$, and since 
$C_G'(P)\nsg K$, it follows that $K$ is a split extension of $C_G'(P)$ by
$\widebar{P}$. Hence for any $x\in N_{\widebar{G}}(P)$, $x\widebar{P}x^{-1} 
= y\widebar{P}y^{-1}$ for some $y\in C_G'(P)$. Consequently, the 
restriction map
	\beq
	N_{\widebar{G}}(\widebar{P})/C_{\widebar{G}}(\widebar{P}) \cong
	\Aut_{\widebar{G}}(\widebar{P})\Right5{} \Aut_{\widebar{G}}(P) 
	\cong N_{\widebar{G}}(P)/C_{\widebar{G}}(P)
	\tag{2} \eeq
is surjective. Also, if $x\in C_{\widebar{G}}(P)\le K$ normalizes
$\widebar{P}$, then $x\in N_K(\widebar{P})\cong \widebar{P}\times
C_G'(\widebar{P})$ by (1), and so $c_x\in\Inn(\widebar{P})$. Thus
the kernel of the map in (2) is contained in
$\Inn(\widebar{P})$. Consequently,
	\begin{small}  
	$$ \Out_{\widebar{G}}(\widebar{P}) =
	\Aut_{\widebar{G}}(\widebar{P}) / \Inn(\widebar{P}) \cong 
	\Aut_{\widebar{G}}(P) / \Aut_{\widebar{P}}(P)
	\cong \Out_{\widebar{G}}(P)/O_p(\Out_{\widebar{G}}(P)), $$
	\end{small}%
and it follows that $\widebar{P}$ is $\calf_p(\widebar{G})$--radical.
\end{proof}

This is now applied to show the following:

\begin{Prop} \label{Spin-cent-rad}
Fix an odd prime power $q$, and let $P\le\Spin_7(q)$ be any subgroup which 
is 2--centric and $\calf_2(\Spin_7(q))$--radical.  Then $P$ is centric in 
$\Spin_7(\fqbar)$; ie, $C_{\Spin_7(\fqbar)}(P)=Z(P)$.
\end{Prop}

\begin{proof}  Let $z$ be the central involution in $\Spin_7(q)$.  By 
Lemma \ref{str-p-rad}(e), $z\in{}P$, and $\widebar{P}\defeq{}P/\gen{z}$ is 
2--centric in $\Omega_7(q)$ and is $\calf_2(\Omega_7(q))$--radical. So by 
Lemma \ref{str-p-rad}(d), there is a 2--subgroup $\widehat{P}\le{}O_7(q)$ 
such that $\widehat{P}\cap\Omega_7(q)=\widebar{P}$, and such that 
$\widehat{P}$ is 2--centric in $O_7(q)$ and is $\calf_2(O_7(q))$--radical.

Let $V=\bigoplus_{i=1}^mV_i$ be a maximal decomposition of $V$ as an 
orthogonal direct sum of $\widehat{P}$--representations, and set 
$\bb_i=\bb|_{V_i}$.  We assume these are arranged so that for some $k$, 
$\dim(V_i)>1$ when $i\le{}k$ and $\dim(V_i)=1$ when $i>k$.  Let $V_+$ be 
the sum of those 1--dimensional components $V_i$ with square discriminant, 
and let $V_-$ be the sum of those 1--dimensional components $V_i$ with 
nonsquare discriminant.  We will be referring to the two decompositions
	$$ (V,\bb) = \bigoplus_{i=1}^m(V_i,\bb_i)
	= \bigoplus_{i=1}^k(V_i,\bb_i) \oplus 
	(V_+,\bb_+) \oplus (V_-,\bb_-) , $$
both of which are orthogonal direct sums.  We also write
	$$ V^{(\infty)}=\fqbar\otimes_{\F_q}V
	\qquad\textup{and}\qquad
	V_i^{(\infty)}=\fqbar\otimes_{\F_q}V_i, $$
and let $\bb^{(\infty)}$ and $\bb_i^{(\infty)}$ be the induced quadratic 
forms. 

\smallskip

\noindent\textbf{Step 1}\qua For each $i$, set 
	$$ D_i=\{\pm\Id_{V_i}\}\le O(V_i,\bb_i), $$
a subgroup of order 2; and write
	$$ D = \prod_{i=1}^mD_i\le O(V,\bb),\qquad\textup{and}\qquad
	D_\pm = \prod_{V_i\subseteq{}V_\pm}D_i \le O(V\pm,\bb_\pm). $$
Thus $D$ and $D_\pm$ are elementary abelian 2--groups of rank $m$ and 
$\dim(V_\pm)$, respectively.  We first claim that
	\beq \widehat{P} \ge D, \tag{1} \eeq
and that
	\beq \widehat{P} = \prod_{i=1}^m P_i \quad
	\textup{where $\forall\, i$, $P_i$ is 2--centric in $O(V_i,\bb_i)$ 
	and $\calf_2(O(V_i,\bb_i))$--radical.}
	\tag{2} \eeq
Clearly, $[D,\widehat{P}]=1$ (and $D$ is a 2--group), so 
$D\le{}\widehat{P}$ since $\widehat{P}$ is 2--centric.  This proves (1).  
The $V_i$ are thus distinct (pairwise nonisomorphic) as 
$\widehat{P}$--representations, since they are pairwise nonisomorphic as 
$D$--representations.  The decomposition as a sum of $V_i$'s is thus unique 
(not only up to isomorphism), since $\Hom_{\widehat{P}}(V_i,V_j)=0$ for 
$i\ne{}j$. 

Let $\widehat{C}$ be the group of elements of $O(V,\bb)$ which send each 
$V_i$ to itself, and let $\widehat{N}$ be the group of elements which permute 
the $V_i$.  By the uniqueness of the decomposition of $V$, 
	$$ \widehat{P} {\cdot} C_{O(V,\bb)}(\widehat{P}) \le 
	\widehat{C}=\prod_{i=1}^mO(V_i,\bb_i)
	\qquad\textup{and}\qquad
	N_{O(V,\bb)}(\widehat{P}) \le \widehat{N}. $$
Since $\widehat{P}$ is 2--centric in $O(V,\bb)$ and 
$\calf_2(O(V,\bb))$--radical, it 
is also 2--centric in $\widehat{N}$ and $\calf_2(\widehat{N})$--radical 
(this holds for any subgroup which contains $N_{O(V,\bb)}(\widehat{P})$).  So 
by 
Lemma \ref{str-p-rad}(b) (and since $\widehat{C}\nsg\widehat{N}$), 
$\widehat{P}$ 
is 2--centric in $\widehat{C}$ and $\calf_2(\widehat{C})$--radical.  Point 
(2) now follows from Lemma \ref{str-p-rad}(a).

\smallskip

\noindent\textbf{Step 2}\qua Whenever $\dim(V_i)>1$ (ie, $1\le{}i\le{}k$), 
then by Lemma \ref{2-irred}, $\dim(V_i)$ is even, and $\bb_i$ has square 
discriminant.  So by Lemma \ref{SO-invol}(a), 
$-\Id_{V_i}\in\Omega(V_i,\bb_i)$ for such $i$.  Together with (1), this 
shows that
	\beq \widebar{P} = \widehat{P} \cap \Omega_7(q) 
	\ge \prod_{i=1}^k D_i \times 
	\bigl(\Omega(V_+,\bb_+)\cap{}D_+\bigr) \times 
	\bigl(\Omega(V_-,\bb_-)\cap{}D_-\bigr) \,. \tag{3} \eeq
Also, by Lemma \ref{SO-invol}(a) again, 
	\beq \begin{split}  
	\bigl(\Omega(V_\pm,\bb_\pm)\cap{}D_\pm\bigr) 
	= \bigl( SO(V_\pm,&\bb_\pm)\cap{}D_\pm\bigr) \\
	&= \bigl\langle -\Id_{V_i\oplus{}V_j} \,\big|\, 
	k{+}1\le i<j \le m,\ V_i,V_j\subseteq V_\pm \bigr\rangle. 
	\end{split}
	\tag{4} \eeq

\smallskip

\noindent\textbf{Step 3}\qua  By (3) and (4), the $V_i$ are distinct as 
$\widebar{P}$--representations (not only as 
$\widehat{P}$--represen\-ta\-tions), except possibly when $\dim(V_\pm)=2$.  
We first check that this exceptional case cannot occur.  If $\dim(V_+)=2$ 
and its two irreducible summands are isomorphic as 
$\widebar{P}$--representations, then the image of $\widebar{P}$ under 
projection to $O(V_+,\bb_+)$ is just $\{\pm\Id_{V_+}\}$.  Hence we can 
write $V_+=W\oplus{}W'$, where $W{\perp}W'$ are 1--dimensional, 
$\widebar{P}$--invariant, and have nonsquare discriminant.  Also, 
$\dim(V_-)$ is odd, since $V_+$ and the $V_i$ for $i\le{}k$ are all even 
dimensional.  So $-\Id_{V_-\oplus{}W}$ lies in 
$C_{\Omega_7(q)}(\widebar{P})$ but not in $\widebar{P}$.  But this is 
impossible, since $\widebar{P}$ is 2--centric in $\Omega_7(q)$.  The 
argument when $\dim(V_-)=2$ is similar.

The $V_i$ are thus distinct as $\widebar{P}$--representations.  So for all 
$i\ne{}j$, $\Hom_{P}(V_i,V_j)=0$, and hence
	$$ \Hom_{\fqbar[P]}(V_i^{(\infty)},V_j^{(\infty)}) \cong
	\fqbar\otimes_{\F_q} \Hom_{\F_q[P]}(V_i,V_j) = 0. $$
Thus any element of $O(V^{(\infty)},\bb^{(\infty)})$ which centralizes 
$\widebar{P}$ sends each $V_i^{(\infty)}$ to itself.  In other words,
	\beq C_{\Spin_7(\fqbar)}(P)/\gen{z} \le 
	C_{\Omega_7(\fqbar)}(\widebar{P}) \le \prod_{i=1}^m 
	O(V_i^{(\infty)},\bb_i^{(\infty)}). \eeq

If $\dim(V_\pm)\ge2$, then since $\widebar{P}$ contains all involutions in 
$O(V_\pm,\bb_\pm)$ which are $\widehat{P}$--invariant and have even dimensional 
$(-1)$--eigenspace (see (3)), Lemma \ref{SO-invol}(c) shows that each 
element of $\Spin_7(\fqbar)$ which commutes with $P$ must act on $V_\pm$ 
via $\pm\Id$.  Also, for $1\le{}i\le{}k$, since 
$-\Id_{V_i}\in{}\widebar{P}$ by (3), each element in the centralizer of $P$ 
acts on $V_i$ with determinant 1 (Lemma \ref{SO-invol}(c) again).  We thus 
conclude that
	\beq C_{\Spin_7(\fqbar)}(P)/\gen{z} \le 
	\prod_{i=1}^k SO(V_i^{(\infty)},\bb_i^{(\infty)})
	\times \{\pm\Id_{V_+}\} \times \{\pm\Id_{V_-}\} . 
	\tag{5} \eeq

\smallskip

\noindent\textbf{Step 4}\qua We next show that 
	\beq C_{\Spin_7(\fqbar)}(P)/\gen{z} \le 
	\prod_{i=1}^k\{\pm\Id_{V_i}\} \times 
	\{\pm\Id_{V_+}\} \times \{\pm\Id_{V_-}\} .
	\tag{6} \eeq
Using (5), this means showing, for each $1\le{}i\le{}k$, that 
	\beq \pr_i\bigl(C_{\Spin_7(\fqbar)}(P)/\gen{z}\bigr) \le 
	\{\pm\Id_{V_i}\}; \tag{7} \eeq
where $\pr_i$ denotes the projection of 
$O_7(\fqbar)=O(V^{(\infty)},\bb^{(\infty)})$ to 
$O(V_i^{(\infty)},\bb_i^{(\infty)})$.  By Lemma \ref{2-irred}, 
$\dim(V_i)=2$ or $4$.  We consider these two cases separately.

\smallskip

\noindent\textbf{Case 4A}\qua If $\dim(V_i)=4$, then by (2) and 
Lemma \ref{str-p-rad}(c), $P'_i\defeq{}P_i\cap\Omega(V_i,\bb_i)$ is 
2--centric in $\Omega(V_i,\bb_i)$ and is 
$\calf_2(\Omega(V_i,\bb_i))$--radical.  Also, by Proposition \ref{SO3-4},
	$$ \Omega(V_i,\bb_i) \cong \Omega_4^+(q) \cong SL_2(q) \times_{C_2} 
	SL_2(q). $$
By Lemma \ref{str-p-rad}(a,f), under this identification, we have 
$P'_i=Q\times_{C_2}Q'$, where $Q$ and 
$Q'$ are 2--centric in $SL_2(q)$ and $\calf_2(SL_2(q))$--radical.  The Sylow 
2--subgroups of $SL_2(q)$ are quaternion groups of order $\ge8$, all 
subgroups of a quaternion 2--group are quaternion or cyclic, and cyclic 
2--subgroups of $SL_2(q)$ cannot be both 2--centric and 
$\calf_2(SL_2(q))$--radical.  So $Q$ and $Q'$ must be quaternion of order 
$\ge8$.  By \cite[3.6.3]{Suzuki}, any cyclic 2--subgroup of $SL_2(\fqbar)$ 
of order $\ge4$ is conjugate to a subgroup of diagonal matrices, whose 
centralizer is the group of all diagonal matrices in $SL_2(\fqbar)$.  
Knowing this, one easily checks that all nonabelian quaternion 2--subgroups 
of $SL_2(\fqbar)$ are centric in $SL_2(\fqbar)$.  It follows that $P'_i$ 
is centric in 
	$$ SO(V_i^{(\infty)},\bb_i^{(\infty)})\cong 
	SL_2(\fqbar)\times_{C_2}SL_2(\fqbar), $$
and hence that 
	$$ \pr_i\bigl(C_{\Spin_7(\fqbar)}(P)/\gen{z}\bigr) \le 
	C_{SO(V_i^{(\infty)},\bb_i^{(\infty)})}(P'_i) 
	= Z(P'_i) = \{\pm\Id_{V_i}\}. $$
Thus (7) holds in this case.

\smallskip

\noindent\textbf{Case 4B}\qua If $\dim(V_i)=2$, then 
$O(V_i,\bb_i)\cong{}O_2^\pm(q)$ is a dihedral group of order $2(q\mp1)$ 
\cite[Theorem 11.4]{Taylor}.  Hence $P_i\in\Syl_2(O(V_i,\bb_i))$, since 
the Sylow subgroups are the only radical 2--subgroups of a dihedral group.  
Fix $V_j$ for any $k<j\le{}m$, and choose $\alpha\in{}O(V_i,\bb_i)$ of 
determinant $(-1)$ whose $(-1)$--eigenspace has the same discriminant as 
$V_j$.  Since $P_i\in\Syl_2(O(V_i,\bb_i))$, we can assume (after 
conjugating if necessary) that $\alpha\in{}P_i$.  Then 
$(-\Id_{V_j})\oplus\alpha$ lies in $\widebar{P}=\widehat{P}\cap\Omega_7(q)$.  
Hence for any $g\in{}C_{\Spin_7(\fqbar)}(P)/\gen{z}$, 
$\pr_i(g)\in{}O(V_i^{(\infty)},\bb_i^{(\infty)})$ leaves both eigenspaces 
of $\alpha$ invariant, and has determinant 1 by (5).  Thus 
$\pr_i(g)=\pm\Id_{V_i}$; and so (7) holds in this case.

\smallskip

\noindent\textbf{Step 5}\qua Clearly, $-\Id_{V_\pm}$ lies in 
$SO(V_\pm,\bb_\pm)$ if and only if $\dim(V_\pm)$ is even (which is the 
case for exactly one of the two spaces $V_\pm$), and this holds if and only 
if $-\Id_{V_\pm}\in\Omega(V_\pm,\bb_\pm)$.  Also, since each $V_i$ for 
$1\le{}i\le{}k$ has square discriminant (Lemma \ref{2-irred} again),
$-\Id_{V_i}\in\Omega(V_i,\bb_i)$ for all such $i$.  Thus (6) and (1) imply 
that 
	$$ C_{\Spin_7(\fqbar)}(P)/\gen{z} \le 
	\widehat{P}\cap\Omega_7(q)=\widebar{P}, $$
and hence that $P$ is centric in $\Spin_7(\fqbar)$.  
\end{proof}

Proposition \ref{Spin-cent-rad} does \emph{not} hold in general if 
$\Spin_7(-)$ is replaced by an arbitrary algebraic group.  For example, 
assume $q$ is an odd prime power, and let $P\le{}SL_5(q)$ be the group of 
diagonal matrices of 2--power order.  Then $P$ is 2--centric in $SL_5(q)$ 
and $\calf_2(SL_5(q))$--radical, but is definitely not 2--centric in 
$SL_5(\fqbar)$.



\def\theshorttitle{Correction to: Construction of 2--local finite groups}

\volumenumber{9}
\volumeyear{2005}
\pagenumbers{3001}{21}
\published{9 February 2005}

\def\S{Section }

\newpage

\renewcommand{\thesection}{\arabic{section}}

\count0=\startpage


\gt\hfill      
\hbox to 77pt{\vbox to 0pt{\vglue -15pt\epsfbox{gtlogo.eps}\vss}\hss}
\break
{\small Volume 9 (2005) \startpage--\finishpage\ (temporary page numbers)\nl
Erratum 1\nl
Published:  \publishdate}

\section*{}
\addcontentsline{toc}{section}{Erratum}
\vglue -0.4in

\cl{\large\bf  Correction to:  Construction of 2--local finite groups}
\cl{\large\bf of a type 
studied by Solomon and Benson} 
\vglue 0.05truein 

\cl{\sc Ran Levi}
\cl{\sc Bob Oliver}
\vglue 0.03truein 

\cl{\small\sl Department of Mathematical Sciences, University of Aberdeen}
\cl{\small\sl  Meston Building 339, Aberdeen AB24 3UE, UK}
\cl{\small and}
\cl{\small\sl  LAGA, Institut Galil\'ee, Av. J-B Cl\'ement} 
\cl{\small\sl  93430 Villetaneuse, France}

\vglue 5pt 
\cl{\small Email: {\tt\mailto{ran@maths.abdn.ac.uk}}{\qua 
and \qua}{\tt\mailto{bob@math.univ-paris13.fr}}}

\vglue 10pt 

{\small\par\leftskip 25pt\rightskip 25pt
{\bf Abstract}\stdspace
A $p$--local finite group is an algebraic structure with a classifying
space which has many of the properties of $p$--completed classifying
spaces of finite groups. In our paper \cite{Solnew}, we constructed a
family of 2--local finite groups which are ``exotic'' in the following
sense: they are based on certain fusion systems over the Sylow
2--subgroup of $\Spin_7(q)$ ($q$ an odd prime power) shown by Solomon
not to occur as the 2--fusion in any actual finite group.  As predicted
by Benson, the classifying spaces of these 2--local finite groups are
very closely related to the Dwyer--Wilkerson space $BDI(4)$.  An error
in our paper \cite{Solnew} was pointed out to us by Andy Chermak, and we
correct that error here.
\vglue 5pt 

{\bf Keywords}\qua Classifying space, $p$--completion, finite groups, fusion 
\vglue 2pt 

{\bf AMS Classification}\qua 55R35; 55R37, 20D06, 20D20
\par}
\vglue 5pt

\setcounter{section}{1}
\setcounter{Thm}{1}

A \emph{saturated fusion system} over a finite $p$--group $S$ is a category 
whose objects are the subgroups of $S$, whose morphisms are all 
monomorphisms between the subgroups, and which satisfy certain axioms first 
formulated by Puig, and also described at the start of the first section in 
\cite{Solnew}.  The main result of \cite{Solnew} is the construction of 
saturated fusion systems over certain 2--groups, motivated by a theorem of 
Solomon \cite{Solomon}, which implies that these systems cannot be induced 
by fusion in any finite group.  Recently, Andy Chermak has pointed out to us 
that the fusion systems actually constructed in \cite{Solnew} are \emph{not}
saturated (do not satisfy all of Puig's axioms).  In this note, we describe 
how to modify that construction in a way so as to obtain saturated fusion 
systems of the desired type, and explain why all of the results in 
\cite{Solnew} (aside from \cite[Lemma A.10]{Solnew}) are true under this new 
construction.

The following is the main theorem in \cite{Solnew}:

\begin{Thm}{\rm\cite[Theorem 2.1]{Solnew}}\label{newSol(q)}\qua
Let $q$ be an odd prime power, and fix $S\in\Syl_2(\Spin_7(q))$.  Let 
$z\in{}Z(\Spin_7(q))$ be the central element of order 2.  Then there is a 
saturated fusion system $\calf=\calf_{\Sol}(q)$ which satisfies the 
following conditions:
\begin{enumerate}  
\item $C_{\calf}(z)=\calf_S(\Spin_7(q))$ as fusion systems over $S$.

\item All involutions of $S$ are $\calf$--conjugate.
\end{enumerate}
Furthermore, there is a unique centric linking system 
$\call=\call^c_{\Sol}(q)$ associated to $\calf$.
\end{Thm}

We have, in fact, found two errors in our proof of this theorem which we 
correct here.  The more serious one is in \cite[Lemma A.10]{Solnew}, which is 
not true as stated: the last sentence in its proof is wrong. This has 
several implications on the rest of our construction, all of which are 
systematically treated here.  There is also an error in the statement of 
\cite[Lemma 2.8(b)]{Solnew} which is corrected below (Lemma \ref{newE4-props}).

We first state and prove here a corrected version of \cite[Lemma 
A.10]{Solnew}, and then state a modified version of the main technical 
proposition, \cite[Proposition 1.2]{Solnew}, used to prove saturation.  
Afterwards, we describe the changes which are needed in \cite[\S2]{Solnew} to 
prove the main theorem.  In table \ref{newt1}, we list the correspondence 
between results and proofs in \cite[\S2]{Solnew} and those here.  This is 
intended as a guide to the reader who is not yet familiar with \cite{Solnew}, 
and who wants to read it simultaneously with this correction.

\def\strutt{\vrule width 0pt depth 6pt}
\def\struttt{\vrule width 0pt height 12pt}

\begin{table}[ht!]
\begin{tabular}{l|l|l}
\textbf{Reference in \cite{Solnew}} & \textbf{Reference here} & 
\textbf{Remarks}\strutt \\ \hline
Proposition 1.2\struttt 
& Proposition \ref{newconstruct-sFs} & more general statement \\
Theorem 2.1 & Theorem \ref{newSol(q)} & unchanged \\
Definition 2.2, Lemma 2.3 & --- & unchanged \\
Definition 2.4 & --- & omit definition of $\Gamma_n$ \\
Prp. 2.5, Def. 2.6, Lem. 2.7 & --- & unchanged \\
--- & Lem. \ref{newC(D8)}, Prp. \ref{newOmega-u} & added \\
--- & Definition \ref{newgamma-u-def} & new definition of $\Gamma_n$ \\
--- & Lemma \ref{newgamma-pres.fus.} & added \\
Lemma 2.8 & Lemma \ref{newE4-props} & (b) restated, \\
&& new proofs of (b), (e) \\
Proposition 2.9 & Proposition \ref{newrk3} & partly new proof \\
Lemma 2.10 & Lemma \ref{newpoint(e)} & partly new proof \\
Proposition 2.11 & Proposition \ref{newF_Sol(fqbar)} & partly new proof 
\end{tabular}\nocolon\caption{}\label{newt1}
\end{table}

The only difference between \cite[Lemma A.10]{Solnew} and the corrected 
version shown here is that in \cite{Solnew}, we claimed that the ``correction 
factor'' $Z$ must lie in a certain subgroup of order 2 in $SL_3(\Z/2^k)$, 
which is definitely not the case.  Also, for convenience, we state this 
lemma here for matrices over the 2--adic integers $\ztwo$, instead of for 
matrices over the finite rings $\Z/2^k$.

\begin{Lem}[{Modified \cite[Lemma A.10]{Solnew}}] \label{newnewA10}
Let $T_1$ and $T_2$ be the two maximal parabolic subgroups of $GL_3(2)$:
        $$ T_1=GL^1_2(\Z/2)=\bigl\{(a_{ij})\in{}GL_3(2)\,|\,
        a_{21}=a_{31}=0 \bigr\} $$
and
        $$ T_2=GL^2_1(\Z/2)=\bigl\{(a_{ij})\in{}GL_3(2)\,|\,
        a_{31}=a_{32}=0 \bigr\}. $$
Set $T_0=T_1\cap{}T_2$: the group of upper triangular matrices in $GL_3(2)$.  
Assume, for some $k\ge2$, that 
        $$ \mu_i\: T_i \Right4{} SL_3(\ztwo) $$
are lifts of the inclusions $T_i\Right1{}GL_3(2)=SL_3(2)$ (for $i=1,2$) 
such that $\mu_1|_{T_0}=\mu_2|_{T_0}$.  Then there is a homomorphism 
        $$ \mu\:GL_3(2)\Right2{}SL_3(\ztwo), $$ 
and an element $Z\in{}C_{SL_3(\ztwo)}(\mu_1(T_0))$,
such that $\mu|_{T_1}=\mu_1$, and $\mu|_{T_2}=c_Z\circ\mu_2$.  
\end{Lem}

\begin{proof}  By \cite[Lemma 4.4]{DW:DI4new}, there is a lifting 
$\mu'\:GL_3(2)\rTo SL_3(\ztwo)$ of the identity on $GL_3(2)$; and any 
two liftings to $SL_3(\ztwo)$ of the inclusion of $T_1$ or of $T_2$ into 
$GL_3(2)$ differ by conjugation by an element of $SL_3(\ztwo)$.  In 
particular, there are elements $Z_1,Z_2\in{}SL_3(\ztwo)$ such that 
$\mu_i=c_{Z_i}\circ\mu'|_{T_i}$ (for $i=1,2$).  
Set $\mu=c_{Z_1}\circ\mu'$, and $Z=Z_1Z_2^{-1}$.  Then $\mu_1=\mu|_{T_i}$, 
and $\mu|_{T_2}=c_Z\circ\mu_2$.  Since $\mu_1|_{T_0}=\mu_2|_{T_0}$, 
conjugation by $Z$ is the identity on $\mu_1(T_0)=\mu_2(T_0)$, and thus 
$Z\in{}C_{SL_3(\ztwo)}(\mu_1(T_0))$.  
\end{proof}

Whenever $G$ is a finite group, $S\in\sylp{G}$, $S_0\nsg{}S$, and 
$\Gamma\le\Aut(S_0)$, then 
        $$ \gen{\calf_S(G);\calf_{S_0}(\Gamma)} $$
denotes the smallest fusion system over $S$ which contains all $G$--fusion, 
and which also contains all restrictions of automorphisms in $\Gamma$.  
In other words, if $\calf$ denotes this fusion system, then for each 
$P,Q\le{}S$, $\Hom_{\calf}(P,Q)$ is the set of all composites 
        $$ P = P_0 \Right2{\varphi_1} P_1 \Right2{\varphi_2} P_2
        \Right2{} \cdots \Right2{}
        P_{k-1} \Right2{\varphi_k} P_k = Q, $$
where for each $i$, $P_i\le{}S$, and either 
$\varphi_i\in\Hom_{G}(P_{i-1},P_i)$, or (if $P_{i-1},P_i\le{}S_0$) 
$\varphi_i=\psi_i|_{P_{i-1}}$ for some $\psi_i\in\Gamma$ such that 
$\psi_i(P_{i-1})=P_i$.  

Whenever $G$ is a finite group and $S\in\sylp{G}$, an automorphism 
$\varphi\in\Aut(S)$ is said to \emph{preserve $G$--fusion} if for each 
$P,Q\le{}S$ and each $\alpha\in\Iso(P,Q)$, $\alpha\in\Iso_G(P,Q)$ if and 
only if $\varphi\alpha\varphi^{-1}\in\Iso_G(\varphi(P),\varphi(Q))$. The 
proof of Theorem \ref{newSol(q)} is based on the following proposition.  

\begin{Prop}[{Modified \cite[Proposition 1.2]{Solnew}}] \label{newconstruct-sFs}
Fix a finite group $G$, a prime $p$ dividing $|G|$, and a Sylow 
$p$--subgroup $S\in\sylp{G}$.  Fix a normal subgroup $Z\nsg{}G$ of order 
$p$, an elementary abelian subgroup $U\nsg{}S$ of rank two containing $Z$ 
such that $C_S(U)\in\sylp{C_G(U)}$, and a group $\Gamma\le\Aut(C_S(U))$ 
of automorphisms which preserve all $C_G(U)$--fusion, and such that 
$\gamma(U)=U$ for all $\gamma\in\Gamma$.  Set 
        $$ S_0=C_S(U) \qquad\textup{and}\qquad 
        \calf=\gen{\calf_S(G);\calf_{S_0}(\Gamma)}, $$  
and assume the following hold.
\begin{enumerate}  
\item All subgroups of order $p$ in $S$ different from $Z$ are $G$--conjugate.

\item $\Gamma$ permutes transitively the subgroups of order $p$ in $U$.

\item $\{\varphi\in\Gamma\,|\,\varphi(Z)=Z\}=\Aut_{N_G(U)}(C_S(U))$.

\item For each $E\le{}S$ which is elementary abelian of rank 
three, contains $U$, and is fully centralized in $\calf_S(G)$, 
        $$ \{\alpha\in\autf(C_S(E))\,|\,\alpha(Z)=Z\} = \Aut_G(C_S(E)). $$

\item For all $E,E'\le{}S$ which are elementary abelian of rank three and 
contain $U$, if $E$ and $E'$ are $\Gamma$--conjugate, then they are 
$G$--conjugate.
\end{enumerate}
Then $\calf$ is a saturated fusion system over $S$.  Also, for 
any $P\le{}S$ such that $Z\le{}P$,
        \beq \{\varphi\in\homf(P,S)\,|\,\varphi(Z)=Z\} = \Hom_G(P,S). 
        \tag{1} \eeq
\end{Prop}

\noindent This proposition is slightly more general than \cite[Proposition 
1.2]{Solnew}, in that $\Gamma$ is assumed only to be a group of automorphisms 
of $C_S(U)$ which preserves $C_G(U)$--fusion, and not a group of 
automorphisms of $C_G(U)$ itself.  This extra generality is necessary when 
proving that the fusion systems $\calf_{\Sol}(q)$, under our modified 
construction, are saturated.  The changes needed to prove this more 
general version of \cite[Proposition 1.2]{Solnew} are described in Section 2.

Proposition \ref{newconstruct-sFs} is applied with $G=\Spin_7(q)$, $Z=Z(G)$, 
and $S\in\Syl_2(G)$, $U\le{}S$, and $\Gamma\le\Aut(C_G(U))$ to be chosen 
shortly.  The error in \cite{Solnew} arose in the choice of $\Gamma$, as will 
be explained in detail below.

We first recall some of the definitions and notation used in \cite{Solnew}.  
Throughout, we fix an odd prime power $q$, let $\F_q$ be a field with $q$ 
elements, and let $\fqbar$ be its algebraic closure.  We write 
$SL_2(q^\infty)=SL_2(\fqbar)$, $\Spin_7(q^\infty)=\Spin_7(\fqbar)$, etc., 
for short.  For each $n$, $\psi^{q^n}$ denotes the automorphism of 
$\Spin_7(q^\infty)$ or of $SL_2(q^\infty)$ induced by the field 
isomorphism $(x\mapsto{}x^{q^n})$.  We then fix elements 
$z,z_1\in\Spin_7(q)$ of order 2, where $\gen{z}=Z(\Spin_7(q))$, set 
$U=\gen{z,z_1}$, and construct an explicit homomorphism
        $$ \omega \: SL_2(q^\infty)^3 \Right4{} \Spin_7(q^\infty) $$
such that 
        $$ \Im(\omega)=C_{\Spin_7(q^\infty)}(U) 
        \qquad\textup{and}\qquad \Ker(\omega)=\gen{(-I,-I,-I)}. $$
Write
        $$ H(q^\infty)\defeq\omega(SL_2(q^\infty)^3)
        =C_{\Spin_7(q^\infty)}(U)
        \quad\textup{and}\quad
        \trp[X_1,X_2,X_3]=\omega(X_1,X_2,X_3)  $$
for short.  In particular, 
        $$ z =\trp[I,I,-I] \qquad\textup{and}\qquad z_1=\trp[-I,I,I], $$
and thus 
        $$ U=\bigl\{\trp[\pm{}I,\pm{}I,\pm{}I]\bigr\} $$ 
(with all combinations of signs).  By \cite[Lemma 2.3 \& Proposition 
2.5]{Solnew}, there is an element $\tau\in{}N_{\Spin_7(q)}(U)$ of order 2 
such that
        \beq \tau{\cdot}\trp[X_1,X_2,X_3]{\cdot}\tau^{-1} = \trp[X_2,X_1,X_3] 
        \eeq
for all $X_1,X_2,X_3\in{}SL_2(q^\infty)$, and such that
        $$ N_{\Spin_7(q^\infty)}(U)=H(q^\infty){\cdot}\gen{\tau}. $$

We next fix elements $A,B\in{}SL_2(q)$ of order 4, such that 
$\gen{A,B}\cong{}Q_8$ (a quaternion group of order 8).  Most of the 
following notation is taken from \cite[Definition 2.6]{Solnew}.  We set
        $$ \AAA=\trp[A,A,A] \qquad\textup{and}\qquad \BBB=\trp[B,B,B]; $$
        $$ C(q^\infty) = \{ X\in{}C_{SL_2(q^\infty)}(A) \,|\,
        X^{2^k}=I,\textup{ some $k$}\} \cong \Z/2^\infty; $$
and 
        $$ Q(q^\infty)=\gen{C(q^\infty),B}. $$
Here, $\Z/2^\infty$ means a union of cyclic 2--groups $\Z/2^n$ for all 
$n$; equivalently, the group $\Z[\frac12]/\Z$. We then define
        \begin{align*}  
        A(q^\infty) &= \omega(C(q^\infty)^3) \cong (\Z/2^\infty)^3, \\
        S_0(q^\infty) &= \omega(Q(q^\infty)^3) \le H(q^\infty) \\
        S(q^\infty) &= S_0(q^\infty){\cdot}\gen{\tau} \le 
        H(q^\infty){\cdot}\gen{\tau} \le \Spin_7(q^\infty). 
        \end{align*}

In all cases, whenever a subgroup $\Theta(q^\infty)\le\Spin_7(q^\infty)$ has 
been defined, we set 
        $$ \Theta(q^n)=\Theta(q^\infty)\cap\Spin_7(q^n). $$  
Since $\Spin_7(q^n)$ is the fixed subgroup of $\psi^{q^n}$ acting on 
$\Spin_7(q^\infty)$ (cf. \cite[Lemma A.3]{Solnew}), $H(q^n)\le H(q^\infty)$ 
is the subgroup of all elements of the form $\trp[X_1,X_2,X_3]$, where either 
$X_i\in{}SL_2(q^n)$ for each $i$, or $\psi^{q^n}(X_i)=-X_i$ for each $i$.  
By \cite[Lemma 2.7]{Solnew}, for all $n$, 
        $$ S(q^n)\in\Syl_2(\Spin_7(q^n)). $$

The following lemma is what is needed to tell us how to choose a subgroup 
$\Gamma_n\le\Aut(S_0(q^n))$ so that the fusion system 
$\gen{\calf_{S(q^n)}(\Spin_7(q^n));\calf_{S_0(q^n)}(\Gamma_n)}$ is 
saturated.  Note that since each element of $C(q^\infty)$ has 2--power 
order, it makes sense to write $X^u\in{}C(q^\infty)$ for 
$X\in{}C(q^\infty)$ and $u\in\ztwo$.  

\begin{Lem} \label{newC(D8)}
Assume $\alpha\in\Aut(A(q^\infty))$ centralizes 
$\Aut_{S(q^\infty)}(A(q^\infty))$.  Then $\alpha$ has the form
        $$ \alpha(\trp[X_1,X_2,X_3])=\trp[X_1^v,X_2^v,X_3^u] $$
for some $u,v\in(\ztwo)^*$.
\end{Lem}

\begin{proof}  Set 
        $$ \Delta_0=\Aut_{S(q^\infty)}(A(q^\infty)) = 
        \gen{c_{\trp[B,I,I]}, c_{\trp[I,B,I]}, c_{\trp[I,I,B]}, c_\tau} $$
for short.  The second equality follows since $S(q^\infty)$ is by 
definition generated by $A(q^\infty)$ and the four elements listed.  Set
        $$ A_1 \defeq \gen{z,z_1,\AAA}\cong{}C_2^3, $$
the 2--torsion subgroup in $A(q^\infty)$.  The image of 
$\Delta_0\le\Aut(A(q^\infty))$ in the group
$\Aut(A_1)\cong{}GL_3(2)$ (the image under restriction) is the 
group of all automorphisms which leave $\gen{z}$ and 
$U=\gen{z,z_1}$ invariant (ie, the group of upper triangular matrices 
with respect to the ordered basis $\{z,z_1,\AAA\}$).  

By assumption, $[\alpha,\Delta_0]=1$, and in particular, 
$\alpha\Delta_0\alpha^{-1}=\Delta_0$.  Since each element of $\Delta_0$ 
sends $U$ to itself, this means that each element of $\Delta_0$ also sends 
$\alpha(U)$ to itself.  Also, $U$ is the only subgroup of rank $2$ left 
invariant by all elements of $\Delta_0$ (since $\Delta_0$ contains all 
automorphisms which leave $\gen{z}$ and $U$ invariant), and hence
$\alpha(U)=U$. 

It follows that $\alpha$ induces an automorphism $\alpha'$ of 
$A(q^\infty)/U=(C(q^\infty)/\gen{-I})^3$.  Also, $\alpha'$ commutes with 
the following automorphisms of $(C(q^\infty)/\gen{-I})^3$: 
        $$ (Y_1,Y_2,Y_3)\mapsto(Y_1^{\pm1},Y_2^{\pm1},Y_3^{\pm1})
        \qquad\textup{and}\qquad
        (Y_1,Y_2,Y_3)\mapsto(Y_2,Y_1,Y_3) $$
since these are induced by automorphisms in $\Delta_0$.  Hence $\alpha'$ has 
the form $\alpha'(Y_1,Y_2,Y_3)=(Y_1^v,Y_2^v,Y_3^u)$ for some 
$u,v\in(\ztwo)^*$.  Thus for all $X_1,X_2,X_3\in{}C(q^\infty)$, 
        $$ \alpha(\trp[X_1,X_2,X_3])=
        \trp[\pm{}X_1^v,\pm{}X_2^v,\pm{}X_3^u]. $$
Since all elements of $C(q^\infty)$ are squares, these signs must all be 
positive, and $\alpha$ has the form 
$\alpha(\trp[X_1,X_2,X_3]=\trp[X_1^v,X_2^v,X_3^u]$.
\end{proof}

For each $u\in(\ztwo)^*$, let $\delta_u\in\Aut(A(q^\infty))$ be the 
automorphism 
	$$\delta_u(\trp[X_1,X_2,X_3])=[X_1,X_2,X_3^u]. $$
Define $\gamma,\gamma_u\in\Aut(A(q^\infty))$ by setting
        $$ \gamma(\trp[X_1,X_2,X_3])=\trp[X_3,X_1,X_2] 
        \qquad\textup{and}\qquad
        \gamma_u=\delta_u\gamma\delta_u^{-1}. $$

\begin{Prop} \label{newOmega-u}
There is an element $u\in\ztwo$ such that $u\equiv1$ (mod 
$4$), and such that the subgroup $\Omega_u\le\Aut(A(q^\infty))$ given by 
        $$ \Omega_u\defeq
        \gen{\Aut_{\Spin_7(q^\infty)}(A(q^\infty)),\gamma_u} $$
is isomorphic to $C_2\times{}GL_3(2)$.  Furthermore, the following
hold:
\begin{enumerate}
\item  The subgroup of 
elements of $\Omega_u$ which act via the identity on all 2--torsion in 
$A(q^\infty)$ has order 2, and contains only the identity and the 
automorphism $(g\mapsto{}g^{-1})$.

\item For each $n\ge1$, 
        $$ \gen{\Aut_{\Spin_7(q^n)}(A(q^n)),\gamma_u} \cong
        C_2\times{}GL_3(2). $$
\end{enumerate}
\end{Prop}

\begin{proof}  For each $k\ge1$, let $A_k\le{}A(q^\infty)$ denote the 
$2^k$--torsion subgroup.  In particular,
        $$ A_1=\gen{z,z_1,\AAA}\cong{}C_2^3, 
        \qquad\textup{where}\qquad
        \AAA=\trp[A,A,A]. $$
Let $R_1,R_2,\dots\in C(q^\infty)\cong\Z/2^\infty$ be elements such that 
$R_1=-I$, $R_2=A$, and $(R_i)^2=R_{i-1}$ for all $k\ge2$.  Thus, 
$|R_i|=2^i$ for all $i$.  For each $k\ge1$, let 
$\{\mathbf{r}_1^{(k)},\mathbf{r}_2^{(k)},\mathbf{r}_3^{(k)}\}$ be 
the basis of $A_k$ defined by 
        $$ \mathbf{r}_1^{(k)}=\trp[I,I,R_k], \qquad 
        \mathbf{r}_2^{(k)}=\trp[R_k,I,I], 
        \qquad\textup{and}\qquad
        \mathbf{r}_3^{(k)}=\trp[R_{k+1},R_{k+1},R_{k+1}]. $$
In particular, $\mathbf{r}_1^{(1)}=z$, $\mathbf{r}_2^{(1)}=z_1$, and 
$\mathbf{r}_3^{(1)}=\AAA$.  Using these bases, we identify 
$\Aut(A_k)=GL_3(\Z/2^k)$ and $\Aut(A(q^\infty))=GL_3(\ztwo)$.  

Set 
        $$ \Delta_0=\Aut_{S(q^\infty)}(A(q^\infty)), \quad
        \Delta_1=\Aut_{\Spin_7(q^\infty)}(A(q^\infty)), 
        \quad\textup{and}\quad
        \Delta_2=\gen{\Delta_0,\gamma}. $$
In particular, $\Delta_2$ is the group of all signed permutations 
        $$ \trp[X_1,X_2,X_3]\mapsto \trp[X_{\sigma(1)}^{\pm1},
        X_{\sigma(2)}^{\pm1},X_{\sigma(3)}^{\pm1}] $$ 
for $\sigma\in\Sigma_3$.  

For each $i=0,1,2$ and each $k\ge1$, let $\Delta_i^{(k)}\le\Aut(A_k)$ be 
the image of $\Delta_i$ under restriction.  By \cite[Proposition 
A.8]{Solnew}, $\Delta_1^{(1)}=\Aut_{\Spin}(A_1)$ is the group of all 
elements of $\Aut(A_1)\cong{}GL_3(\Z/2)$ which send $z$ to itself.  Also,
$\Delta_0^{(1)}$ was seen in the proof of Lemma \ref{newC(D8)} to be 
the group of all automorphisms of $A_1$ which leaves both $z$ and 
$U=\gen{z,z_1}$ invariant; and a similar argument shows that 
$\Delta_2^{(1)}$ is the group of all automorphisms of $A_1$ which 
leaves $U$ invariant.  Hence, with respect to the ordered basis 
$\{z,z_1,\AAA\}$ of $A_1$, each group $\Delta_i^{(1)}\le\Aut(A_1)$ 
($i=0,1,2$) can be identified with the subgroup $T_i\le{}GL_3(\Z/2)$ of 
Lemma \ref{newnewA10}.  

By \cite[Proposition 2.5]{Solnew}, 
        $$ C_{\Spin_7(q^\infty)}(U)=H(q^\infty) \cong 
        (SL_2(q^\infty)^3)/\gen{(-I,-I,-I)}. $$
An element $\trp[X_1,X_2,X_3]\in{}H(q^\infty)$ ($X_i\in{}SL_2(q^\infty)$) 
centralizes $\AAA\defeq\trp[A,A,A]$ if and only if $[X_i,A]=1$ for each 
$i$, or $X_iAX_i^{-1}=-A=A^{-1}$ for each $i$.  
Set $\widebar{C}=C_{SL_2(q^\infty)}(A)$.  This is an abelian group (the 
union of the finite cyclic groups $C_{SL_2(q^n)}(A)$), and 
$N_{SL_2(q^\infty)}(A)=\widebar{C}{\cdot}\gen{B}$.  
Hence, since $A_1=\gen{U,\AAA}$ and $\BBB=\trp[B,B,B]$, we have
        \beq C_{\Spin_7(q^\infty)}(A_1) = C_{H(q^\infty)}(\AAA) = 
        \omega(\widebar{C}{}^3){\cdot}\gen{\BBB}. \tag{1} \eeq
Since $\omega(\widebar{C}{}^3)$ is abelian (thus centralizes 
$A(q^\infty)$), this shows that the kernel of each of the projection maps 
$\Delta_i\Onto2{}\Delta_i^{(1)}$ is generated by conjugation by $\BBB$; 
ie, by the automorphism $(g\mapsto{}g^{-1})$.  

Since each $\Delta_i$ is finite, their elements all have determinant of 
finite order in $(\ztwo)^*$, hence are $\pm1$ in all cases.  Also, for each 
$i$, $\Delta_i$ surjects onto $\Delta_i^{(1)}=T_i$ with kernel generated 
by the automorphism $(g\mapsto{}g^{-1})$ of determinant $(-1)$.  Hence the 
elements of determinant one in $\Delta_i$ are sent isomorphically to 
$T_i$, and define a lift
        $$ \mu_i\: T_i \Right5{} SL_3(\ztwo) $$
with respect to the given bases.  In particular, 
$\mu_1|_{T_0}=\mu_2|_{T_0}=\mu_0$.  

For all $i,j,k\in(\ztwo)^*$, we define $\psi_{i,j,k}\in\Aut(A(q^\infty))$ 
by setting 
        $$ \psi_{i,j,k}(\trp[X_1,X_2,X_3])=\trp[X_1^i,X_2^j,X_3^k] $$
for all $\trp[X_1,X_2,X_3]\in{}A(q^\infty)$.  Recall that 
$c_\tau(\trp[X_1,X_2,X_3])=\trp[X_2,X_1,X_3]$
(by choice of $\tau$).  Since 
$\Delta_0=\Aut_{S(q^\infty)}(A(q^\infty))$ is generated by 
$c_{\trp[B,I,I]}$, $c_{\trp[I,B,I]}$, $c_{\trp[I,I,B]}$, and $c_\tau$ 
(corresponding to generators of $S(q^\infty)/A(q^\infty)$), $\mu_0$ has 
image
        $$ \Im(\mu_0) = \mu_1(T_0) = 
        \gen{\psi_{-1,-1,1}, \psi_{1,-1,-1}, c_\tau}, $$
where $\psi_{-1,-1,1}=c_{\trp[B,B,I]}$ and 
$\psi_{1,-1,-1}=c_{\trp[I,B,B]}$.  

By Lemma \ref{newnewA10}, there is a homomorphism
        $$ \mu\: GL_3(2) \Right5{} SL_3(\ztwo), $$
and an element $Z\in\Aut(A(q^\infty))\cong{}GL_3(\ztwo)$ which commutes 
with all elements of $\mu_1(T_0)$, such that $\mu|_{T_1}=\mu_1$ and 
$\mu|_{T_2}=c_Z\circ\mu_2$.  By Lemma \ref{newC(D8)}, $Z=\psi_{v,v,u}$ (using 
the above notation) for some $u,v\in(\ztwo)^*$.  Since $\psi_{v,v,v}$ lies 
in the center of $\Aut(A(q^\infty))$ (it sends every element to its $v$-th 
power), we can assume that $v=1$ (without changing $c_Z$), and thus that 
$Z=\psi_{1,1,u}=\delta_u$.  Finally, since 
$\delta_{-1}=\psi_{1,1,-1}\in\Delta_0$, we can replace $\delta_u$ by 
$\delta_{-u}$ if necessary, and assume that $u\equiv1$ (mod $4$). 

Under the identification $GL_3(\ztwo)=\Aut(A(q^\infty))$, we now have
        $$ \Aut_{\Spin_7(q^\infty)}(A(q^\infty)) = (g\mapsto g^{-1}) 
        \times \mu(T_1) $$
and
        $$ \gamma_u = \delta_u\gamma\delta_u^{-1} = c_Z(\gamma) = 
        \mu\left(\left(\begin{smallmatrix}0&1&0\\1&1&0\\0&0&1
        \end{smallmatrix}\right)\right), $$
where the matrix is that of $\gamma|_{A_1}$ with respect to 
the basis $\{\mathbf{r}_1^{(1)},\mathbf{r}_2^{(1)},\mathbf{r}_3^{(1)}\}$.  
Also, $T_1$ is a maximal subgroup of $GL_3(2)$ --- the subgroup of 
invertible matrices which send $\mathbf{r}_1^{(1)}$ to itself --- and so 
$T_1$ together with this matrix generate $GL_3(2)$.  Thus 
        \beq \Omega_u = \gen{\Aut_{\Spin_7(q^\infty)}(A(q^\infty)),\gamma_u}
        = \gen{g\mapsto g^{-1}} \times \mu(GL_3(2)) \cong
        C_2\times{}GL_3(2). \tag{2} \eeq
This proves the first claim in the proposition.  Point (a) follows by 
construction, since each nonidentity element of $\mu(GL_3(2))$ acts 
nontrivially on $A_1$.  Point (b) follows from (2), once we know that 
each element of $\Spin_7(q^n)$ which normalizes $A_1$ (hence which 
normalizes $A(q^n)$) also normalizes $A(q^\infty)$ --- and this
follows
from (1).
\end{proof}

We are now in a position to define the fusion systems we want.  Roughly, 
they are generated by the fusion systems of $\Spin_7(q^n)$ together with 
one extra automorphism:  the cyclic permutation 
$\trp[X_1,X_2,X_3]\mapsto\trp[X_3,X_1,X_2]$ ``twisted'' by the automorphism 
$\delta_u$ of the last proposition.  By comparison, the construction in 
\cite{Solnew} was similar but without the twisting (ie, done with $u=1$), and 
the resulting fusion system is, in fact, not saturated.

We regard $Q(q^\infty)=C(q^\infty){\rtimes}\gen{B}$ as an infinite 
quaternion group:  $BA'B^{-1}=A'{}^{-1}$ for each 
$A'\in{}C(q^\infty)$, and each element of the coset $C(q^\infty){\cdot}B$ 
has order $4$.  Hence any automorphism of $C(q^\infty)$ extends to a unique 
automorphism of $Q(q^\infty)$ which sends $B$ to itself.

\begin{Defi} \label{newgamma-u-def}
Let $u\in(\ztwo)^*$ be as in Proposition \ref{newOmega-u}.  Let
        $$ \widehat{\gamma},\widehat{\delta}_u\in \Aut(S_0(q^\infty)) $$ 
be the automorphisms
        $$ \widehat\gamma(\trp[X_1,X_2,X_3])=\trp[X_3,X_1,X_2] $$
and
        $$ \widehat\delta_u(\trp[X_1,X_2,A'B^j])=\trp[X_1,X_2,(A')^uB^j] $$
for all $X_i\in{}Q(q^\infty)$, $A'\in{}C(q^\infty)$, and $i,j\in\Z$; and 
set $\widehat{\gamma}_u= \widehat{\delta}_u\widehat{\gamma} 
\widehat{\delta}_u^{-1}\in \Aut(S_0(q^\infty))$.  For each $n\ge1$, set 
        $$ \Gamma_n = \gen{\Inn(S_0(q^n)),c_\tau,\widehat{\gamma}_u}
        \le \Aut(S_0(q^n)); $$
and set 
        $$ \calf_n = \calf\sol(q^n) = 
        \gen{\calf_{S(q^n)}(\Spin_7(q^n)),\calf_{S_0(q^n)}(\Gamma_n)}. $$
\end{Defi}

In order to be able to apply Proposition \ref{newconstruct-sFs}, it is 
important to know that $\Gamma_n$ is fusion preserving.  This follows 
immediately from the following lemma.

\begin{Lem} \label{newgamma-pres.fus.}
For each $n\ge1$, the automorphisms $c_\tau$, $\widehat{\gamma}$, 
$\widehat{\delta}_u$, and $\widehat{\gamma}_u$ all preserve 
$H(q^n)$--fusion in $S_0(q^n)$.  
\end{Lem}

\begin{proof}  
The fusion in $SL_2(q^n)$ is generated by inner automorphisms of its Sylow 
subgroup $Q(q^n)$, together with the groups $\Aut(P)$ for subgroups 
$P\le{}Q(q^n)$ isomorphic to $Q_8$.  This follows from Alperin's fusion 
theorem, since these are the only subgroups whose automorphism group is 
not a 2--group.  Thus any automorphism of $Q(q^n)$ --- in particular, the 
automorphism $\varphi_u$ defined by $\varphi_u(A')=A'{}^u$, $\varphi_u(B)=B$ 
for $A'\in{}C(q^n)$ --- preserves fusion.  

Set $H_0(q^n)=\{\trp[X_1,X_2,X_3]\,|\,X_i\in{}SL_2(q^n),\ \forall i\}$.  
Fix a generator $Y$ of $C(q^{2n})$; then $C(q^n)=\gen{Y^2}$, and 
$\psi^{q^n}(Y)=-Y$.  For each $g=\trp[X_1,X_2,X_3]\in{}H(q^n)$, 
$\trp[\psi^{q^n}(X_1),\psi^{q^n}(X_2),\psi^{q^n}(X_3)]=\trp[X_1,X_2,X_3]$; and 
hence there is some fixed $\epsilon=\pm1$ for which 
$\psi^{q^n}(X_i)=\epsilon{\cdot}X_i$ for $i=1,2,3$.  When $\epsilon=1$, this 
means that $X_i\in{}SL_2(q^n)$ for each $i$; while if $\epsilon=-1$ it 
means that $X_i\in{}SL_2(q^n){\cdot}Y$ for each $i$.  So every element of 
$H(q^n)$ either lies in $H_0(q^n)$, or has the form $g{\cdot}\trp[Y,Y,Y]$ 
for some $g\in{}H_0(q^n)$.  

Assume $g\in{}H(q^n)$ and $P,Q\le{}S_0(q^n)$ are such that $gPg^{-1}=Q$.  
Clearly, $P\le{}H_0(q^n)$ if and only if $Q\le{}H_0(q^n)$ ($H_0(q^n)$ is 
normal in $H(q^n)$).  
We claim that there is $h\in{}H(q^n)$ such that 
        \beq h\widehat{\delta}_u(P)h^{-1}=\widehat{\delta}_u(Q)
        \qquad\textup{and}\qquad
        c_h\circ\widehat{\delta}_u|_P=\widehat{\delta}_u\circ{}c_g|_P \,.
        \tag{1} \eeq
Let $P_i,Q_i\le{}SL_2(q^{2n})$ be the projections of $P$ and $Q$ to the 
$i$-th factor ($i=1,2,3$), and write $g=\trp[X_1,X_2,X_3]$ (thus 
$X_iP_iX_i^{-1}=Q_i$).  Consider the following cases.  
\begin{enumerate} 
\item Assume $g\in{}H_0(q^n)$ and $P,Q\le{}H_0(q^n)$.  Since $\varphi_u$ 
preserves fusion in $SL_2(q^n)$, there is $Y_3\in{}SL_2(q^n)$ such 
that $Y_3\varphi_u(P_3)Y_3^{-1}=\varphi_u(Q_3)$ and
$c_{Y_3}\circ\varphi_u|_{P_3}=\varphi_u\circ{}c_{X_3}|_{P_3}$.  Then 
$h\defeq\trp[X_1,X_2,Y_3]$ satisfies (1).  

\item Assume $g\notin{}H_0(q^n)$ and $P,Q\le{}H_0(q^n)$.  Write 
$g=g'{\cdot}\trp[Y,Y,Y]$, 
where $g'\in{}H_0(q^n)$.  Choose $h'$ as in (a), so that (1) is satisfied 
with $g,h$ replaced by $g',h'$.  Then the element $h=h'{\cdot}\trp[Y,Y,Y^u]$ 
satisfies (1). 

\item Finally, assume that $P,Q\nleq{}H_0(q^n)$.  Then none of 
the subgroups $P_i,Q_i\le{}SL_2(q^{2n})$ is contained in $SL_2(q^n)$.
By the same procedure as was used in (a), we can find
$h\in{}H_0(q^{2n})$ which satisfies (1); the problem is to do this so
that $h\in{}H(q^n)$. 

As noted above, fusion in $SL_2(q^{2n})$ is generated by inner
automorphisms of $Q(q^{2n})$ and automorphisms of subgroups isomorphic
to $Q_8$.  Hence if $P_i$ is not isomorphic to $Q_8$ or one of its
subgroups, then there is $X'_i\in{}Q(q^{2n})$ such that 
$c_{X_i}|_{P_i}=c_{X'_i}|_{P_i}$.  If, on the other hand, 
$P_i$ (and hence $Q_i$) is isomorphic to $Q_8$ or $C_4$, then 
$P_i\le\gen{A,Y^rB}$ and $Q_i\le\gen{A,Y^sB}$ for some odd $r,s\in\Z$, 
and we can choose $k\in\Z$ such that $X'_i\defeq{}Y^k\in{}Q(q^{2n})$ 
has the same conjugation action as $X_i$.  Thus
in all cases, we can write $X_i=X'_iX''_i$ for some $X'_i\in{}Q(q^{2n})$ 
and some $X''_i\in{}C_{SL_2(q^{2n}){\cdot}\gen{Y}}(P_i)$.  

The subgroups of $Q(q^{2n})$ which are centralized by elements in the 
coset $SL_2(q^n){\cdot}Y$ are precisely the cyclic subgroups.  (The 
quaternion subgroups of order $\ge8$ are all centric in $SL_2(q^{2n})$.)  
Hence we can choose elements $Y''_i$ as follows:  $Y''_i=1$ if 
$X''_i\in{}SL_2(q^n)$, and $Y''_i\in{}SL_2(q^n){\cdot}Y$ and centralizes 
$\varphi_u(P_i)$ if $X''_i\in{}SL_2(q^n){\cdot}Y$.  We now define
	$$h=\trp[X'_1,X'_2,\varphi_u(X'_3)]{\cdot}\trp[Y''_1,Y''_2,Y''_3];$$
then $h\in{}H(q^n)$ and satisfies (1).

\end{enumerate}

This shows that $\widehat{\delta}_u$ preserves $H(q^n)$--fusion as an 
automorphism of $S_0(q^n)$.  Also, $\widehat{\gamma}$ and $c_\tau$ preserve 
fusion, since both extend to automorphisms of $H(q^n)$; and hence 
$\widehat{\gamma}_u= \widehat{\delta}_u\widehat{\gamma} 
\widehat{\delta}_u^{-1}$ also preserves fusion.  
\end{proof}

Let $\psi=\psi^{q^n}\in\Aut(\Spin_7(q^\infty))$ be induced by the field 
automorphism $x\mapsto{}x^{q^n}$.  
By \cite[Proposition A.9(a)]{Solnew}, if $E\le\Spin_7(q^n)$ is an arbitrary 
elementary abelian 2--subgroup of rank 4, then there is an element 
$a\in\Spin_7(q^\infty)$ such that $aEa^{-1}=E_*$, and we define
        $$ x_\calc(E) = a^{-1}\psi(a). $$
Then $x_\calc(E)\in{}E$, and is independent of the choice of $a$.  

In the following lemma, we correct the statement and proof of points
(b) and (e).  The proof of (e) is affected by both the changes
in the statement of (b) and those in the definition of $\Gamma_n$.

\begin{Lem}{\rm\cite[Lemma 2.8]{Solnew}}\label{newE4-props}\qua
Fix $n\ge1$, set $E_*=\gen{z,z_1,\AAA,\BBB}\le{}S(q^n)$, and let $\calc$ 
be the $\Spin_7(q^n)$--conjugacy class of $E_*$.  Let $\EUfour$ be the set 
of all elementary abelian subgroups $E\le{}S(q^n)$ of rank $4$ which contain 
$U=\gen{z,z_1}$.  Fix a generator $X\in{}C(q^n)$ (the 2--power torsion in 
$C_{SL_2(q^n)}(A)$), and choose $Y\in{}C(q^{2n})$ such that $Y^2=X$.  Then 
the following hold.
\begin{enumerate}  
\item $E_*$ has type I.

\item Each subgroup in $\EUfour$ which contains $\AAA$ is of the form
	\begin{multline*}  
        E_{ijk}=\gen{z,z_1,\AAA,\trp[X^iB,X^jB,X^kB]}
        \quad\textup{or}\\
        E'_{ijk}=\gen{z,z_1,\AAA,\trp[X^iYB,X^jYB,X^kYB]}. 
	\end{multline*}
Each subgroup in $\EUfour$ is $H(q^n)$--conjugate to one of these subgroups  
$E_{ijk}$ or $E'_{ijk}$ for some $i,j,k\in\Z$.

\item $x_{\calc}(E_{ijk})=\trp[(-I)^i,(-I)^j,(-I)^k]$ and
$x_{\calc}(E'_{ijk})=\trp[(-I)^i,(-I)^j,(-I)^k]{\cdot}\AAA$.

\item All of the subgroups $E'_{ijk}$ have type II.  The subgroup 
$E_{ijk}$ has type I if and only if $i\equiv{}j$ (mod $2$), and lies in 
$\calc$ (is conjugate to $E_*$) if and only if $i\equiv{}j\equiv{}k$ (mod 
$2$).  The subgroups $E_{000}$, $E_{001}$, and $E_{100}$ thus represent 
the three conjugacy classes of rank four elementary abelian subgroups of 
$\Spin_7(q^n)$ (and $E_*=E_{000}$).

\item For any $\varphi\in\Gamma_n\le\Aut(S_0(q^n))$ and any $E\in\EUfour$,
$\varphi(x_{\calc}(E))=x_{\calc}(\varphi(E))$.
\end{enumerate}
\end{Lem}

\begin{proof}  We prove only points (b) and (e) here, and refer to 
\cite{Solnew} for the proofs of the other points.

\smallskip

\noindent\textbf{(b) } Assume first that $\AAA\in{}E$; 
ie, that $E\ge A_1=\gen{z,z_1,\AAA}$.  By definition,  
$S(q^n)$ is generated by $A(q^n)$, whose elements clearly
centralize $A_1$; and elements $\trp[B^i,B^j,B^k]$ for $i,j,k\in\{0,1\}$.
Since an element of this form centralizes $\AAA$ only if $i=j=k$, this
shows that
        $$ E\le{}C_{S(q^n)}(\gen{z,z_1,\AAA})=A(q^n){\cdot}\gen{\BBB}. $$
Since $A(q^n)$ is a finite abelian 2--group of rank $3$, we have
$E=\gen{z,z_1,\AAA,g\BBB}$ for some $g\in{}A(q^n)$.  Also,
	\begin{small}  
        $$ A(q^n) = \omega(C(q^\infty)^3)\cap{}\Spin_7(q^n) =
        \bigl\{\trp[X^i,X^j,X^k],\trp[X^iY,X^jY,X^kY]\,\big|\,
        i,j,k\in\Z\bigr\}, $$
	\end{small}%
and hence $E=E_{ijk}$ or $E'_{ijk}$ for some $i,j,k\in\Z$.

Now let $E\in\EUfour$ be arbitrary.  Each element of $E$ has the form 
$\trp[X_1,X_2,X_3]$, where either $X_i\in{}SL_2(q^n)$ for all $i$, or 
$X_i\in{}SL_2(q^n){\cdot}Y$ for all $i$ --- and the elements of the first 
type ($X_i\in{}SL_2(q^n)$) form a subgroup of index at most 2.  Since $U$ 
has index $4$ in $E$, this means that there is some 
$g=\trp[X_1,X_2,X_3]\in{}E{\sminus}U$ for which $X_i\in{}SL_2(q^n)$ for all 
$i$.  Also, $|X_i|=4$ for all $i$, since 
$g\notin{}U=\{\trp[\pm{}I,\pm{}I,\pm{}I]\}$.  Since all elements of order 
$4$ in $SL_2(q^n)$ are conjugate (cf. \cite[3.6.23]{Suzukinew}), this implies 
that $g$ is $H(q^n)$--conjugate to $\AAA$, and hence that $E$ is 
$H(q^n)$--conjugate to one of the above subgroups $E_{ijk}$ or $E'_{ijk}$.

\smallskip

\noindent\textbf{(e) } By construction, $x_{\calc}(-)$ is preserved under 
conjugation by elements of the group $\Spin_7(q^n)$.  Since $\Gamma_n$ is 
generated by 
$\widehat{\gamma}_u$ and conjugation by elements of $\Spin_7(q^n)$, it 
suffices to prove the result when $\varphi=\widehat{\gamma}_u$.  Since 
$u\equiv1$ (mod $4$) by Proposition \ref{newOmega-u}, 
$\widehat{\gamma}_u(\AAA)=\AAA$, 
$\widehat{\gamma}_u(E_{ijk})=E_{k',i',j'}$, and 
$\widehat{\gamma}_u(E'_{ijk})=E'_{k'',i'',j''}$ for some 
$i'\equiv{}i''\equiv{}i$ (mod $2$), and similarly for the other indices.  
Hence by (c), $\widehat{\gamma}_u(x_{\calc}(E))=
x_{\calc}(\widehat{\gamma}_u(E))$ whenever $E=E_{ijk}$ or $E=E'_{ijk}$ for 
some $i,j,k\in\Z$.

Now assume $E\in\EUfour$ is not one of the subgroups $E_{ijk}$ or $E'_{ijk}$.  
By (b), there is $g\in{}H(q^n)$ such that $E'\defeq{}gEg^{-1}$ is of this 
form.  Since $\widehat{\gamma}_u$ preserves $H(q^n)$--fusion by Lemma 
\ref{newgamma-pres.fus.}, there is $h\in{}H(q^n)$ such that the following 
square commutes:
        \begin{diagram}  
        E & \rTo^{\widehat{\gamma}_u} & \widehat{\gamma}_u(E) \\
        \dTo<{c_g} && \dTo<{c_h} \\
        E' & \rTo^{\widehat{\gamma}_u} & \widehat{\gamma}_u(E') \rlap{\,.}
        \end{diagram}
We have seen that 
$\widehat{\gamma}_u(x_\calc(E'))=x_{\calc}(\widehat{\gamma}_u(E'))$; and 
also that $c_g$ and $c_h$ preserve $x_{\calc}(-)$.  Hence 
$\widehat{\gamma}_u(x_\calc(E))=x_{\calc}(\widehat{\gamma}_u(E))$ by the 
commutativity of the square.
\end{proof}

The following is the crucial result needed to apply Proposition 
\ref{newconstruct-sFs}.  The statement is exactly the same as that in 
\cite{Solnew}, but the proof has to be modified slightly due to the changed 
definition of $\Gamma_n$ (hence of $\calf_n$).  

\begin{Prop}{\rm\cite[Proposition 2.9]{Solnew}}\label{newrk3}\qua
Fix $n\ge1$.  Let $E\le{}S(q^n)$ be an elementary abelian
subgroup of rank $3$ which contains $U$, and such that
$C_{S(q^n)}(E)$ $\in\Syl_2(C_{\Spin_7(q^n)}(E))$.  Then \beq
\{\varphi\in\Aut_{\calf_n}(C_{S(q^n)}(E))\,|\,\varphi(z)=z\} =
\Aut_{\Spin_7(q^n)}(C_{S(q^n)}(E)). \tag{1} \eeq
\end{Prop}

\begin{proof}  Set 
        $$ \Spin=\Spin_7(q^n), \qquad S=S(q^n), \qquad \Gamma=\Gamma_n,
        \qquad\textup{and}\qquad \calf=\calf_n $$ 
for short, and consider the subgroups
        $$ R_0=R_0(q^n)\defeq A(q^n)
        \qquad\textup{and}\qquad
        R_1=R_1(q^n)\defeq C_S(\gen{U,\AAA})=\gen{R_0,\BBB}. $$
Then
        $$ R_0 \cong (C_{2^k})^3 \qquad\textup{and}\qquad
        R_1 = R_0 \rtimes \gen{\BBB}, $$
where $2^k$ is the largest power which divides $q^n\pm1$, and where 
$\BBB=\trp[B,B,B]$ has order 2 and acts on $R_0$ via $(g\mapsto{}g^{-1})$. 
Also,
        $$ \gen{U,\AAA} = \gen{\trp[\pm{}I,\pm{}I,\pm{}I],\trp[A,A,A]} 
        \cong C_2^3 $$ 
is the 2--torsion subgroup of $R_0$.  It was shown in the proof of 
\cite[Proposition 2.9]{Solnew} that
        \beq \textup{$R_0$ is the only subgroup of $S$ isomorphic to 
        $(C_{2^k})^3$.} \tag{2} \eeq

Let $E\le{}S$ be an elementary abelian subgroup of rank $3$ which
contains $U$, and such that $C_{S}(E)\in\Syl_2(C_{\Spin}(E))$.  
There are two cases to
consider:  that where $E\le{}R_0$ and that where $E\nleq{}R_0$.  

\smallskip

\noindent\textbf{Case 1}\qua 
Assume $E\le{}R_0$.  Since $R_0$ is abelian of rank 3, we must
have $E=\gen{U,\AAA}$, the 2--torsion subgroup of $R_0$, and $C_S(E)=R_1$.  
Also, by (2), neither $R_0$ nor $R_1$ is isomorphic to any
other subgroup of $S$; and hence
        \beq \Aut_{\calf}(R_i) = \bigl\langle 
        \Aut_{\Spin}(R_i),\Aut_{\Gamma}(R_i) \bigr\rangle
        = \bigl\langle \Aut_{\Spin}(R_i),\gamma_u|_{R_i} \bigr\rangle
        \qquad \textup{for $i=0,1$.} \tag{3} \eeq
(Recall that $\Aut_{\Gamma}(R_i)$ is generated by
$\Aut_{S_0(q^n)}(R_i)$ and the restrictions of $c_\tau$ and $\gamma_u$.)
Hence by Proposition \ref{newOmega-u}(b), 
        \beq \Aut_{\calf}(R_0) = \bigl\langle
        \Aut_{\Spin}(R_0),\gamma_u|_{R_0} \bigr\rangle \cong{}
        C_2 \times GL_3(2). \eeq
In other words, if we let $A_1=\gen{z,z_1,\AAA}$ denote the 2--torsion 
subgroup of $R_0$, then restriction to $A_1$ sends $\Aut_\calf(R_0)$ onto 
$\Aut(A_1)$ with kernel $\gen{c_{\BBB}}$ of order 2.  Since 
$\Aut_{\Spin}(A_1)$ is the group of all automorphisms of $A_1$ which send 
$z$ to itself \cite[Proposition A.8]{Solnew}, this shows that
        \beq \Aut_{\Spin}(R_0) = 
        \bigl\{\varphi\in\Aut_\calf(R_0) \,\big|\, \varphi(z)=z \bigr\}. 
        \tag{4} \eeq

In the proof of \cite[Proposition 2.9]{Solnew} (see formula (5) in that 
proof), we show the first of the following two equalities:
        \begin{align*}  
        \Aut_{\Spin}(R_1) 
        &= \Inn(R_1){\cdot} 
        \bigl\{\varphi\in\Aut_{\Spin}(R_1)\,\big|\,
        \varphi(\BBB)=\BBB \bigr\} \\
        &= \Inn(R_1){\cdot} 
        \bigl\{\varphi\in\Aut(R_1)\,\big|\,\varphi(\BBB)=\BBB,\
        \varphi|_{R_0} \in \Aut_{\Spin}(R_0) \bigr\} . 
        \tag{5}
        \end{align*}
The second equality holds since $R_1=\gen{R_0,\BBB}$, and since $R_0$ is 
the unique abelian subgroup of $R_1$ of index 2.  Since 
$\widehat{\gamma}_u(R_0)=R_0$ and $\widehat{\gamma}_u(\BBB)=\BBB$, this, 
together with (3), shows that 
        $$ \Aut_{\calf}(R_1) = \Inn(R_1) {\cdot}
        \bigl\{\varphi\in\Aut(R_1) \,\big|\, \varphi(\BBB)=\BBB,\ 
        \varphi|_{R_0} \in \Aut_{\calf}(R_0) \bigr\}, $$
and hence that
        \begin{align*}  
        \bigl\{\varphi\in&\Aut_{\calf}(R_1) \,\big|\, \varphi(z)=z \bigr\} 
        \\
        &= \Inn(R_1){\cdot} 
        \bigl\{\varphi\in\Aut(R_1) \,\big|\, \varphi(\BBB)=\BBB,\ 
        \varphi(z)=z,\ 
        \varphi|_{R_0} \in \Aut_{\calf}(R_0) \bigr\} \\
        &= \Inn(R_1){\cdot} 
        \bigl\{\varphi\in\Aut(R_1) \,\big|\, \varphi(\BBB)=\BBB,\ 
        \varphi|_{R_0} \in \Aut_{\Spin}(R_0) \bigr\} , \tag{6}
        \end{align*}
where the second equality follows from (4).  If $\varphi\in\Aut(R_1)$ is 
such that $\varphi(\BBB)=\BBB$ and $\varphi|_{R_0}=c_x|_{R_0}$ for some 
$x\in{}N_{\Spin}(R_0)$, then $x$ normalizes $R_1$ (the centralizer of the 
2--torsion in $R_0$), $c_x(\BBB)=c_y(\BBB)$ for some 
$y\in{}R_1=R_0{\cdot}\gen{\BBB}$ by (5), so we can assume $y\in{}R_0$.  
Since $R_0=A(q^\infty)$ is abelian (so $[y,R_0]=1$), this implies 
$c_{y^{-1}x}|_{R_1}=\varphi$, and hence $\varphi\in\Aut_{\Spin}(R_1)$.  So 
$\{\varphi\in \Aut_{\calf}(R_1)\,|\, \varphi(z)=z\} \le\Aut_{\Spin}(R_1)$ 
by (6), and the opposite inclusion is clear.

\smallskip

\noindent\textbf{Case 2}\qua Now assume that $E\nleq{}R_0$.  By assumption, 
$U\le{}E$ (hence $E\le{}C_S(E)\le{}C_S(U)$), and $C_S(E)$ is a Sylow 
subgroup of $C_{\Spin}(E)$.  Also, $E$ contains an element of the form 
$g{\cdot}\trp[B^i,B^j,B^k]$ for $g\in{}R_0=A(q^n)$ and some $i,j,k$ not
all even, and hence $A(q^n)\nleq{}C_S(E)$.  Hence by (2), 
$C_S(E)$ is not isomorphic to 
$R_1=C_S(\gen{z,z_1,\AAA})$, and this shows that $E$ is not 
$\Spin$--conjugate to $\gen{z,z_1,\AAA}$.  By \cite[Proposition A.8]{Solnew}, 
$\Spin$ contains exactly two conjugacy classes of rank 3 subgroups 
containing $z$, and thus $E$ must have type II.  So by \cite[Proposition 
A.8(d)]{Solnew}, $C_S(E)$ is elementary abelian of rank 4, and also 
has type II.  

Let $\calc$ be the $\Spin_7(q^n)$--conjugacy class of the subgroup 
$E_*=\gen{U,\AAA,\BBB}\cong{}C_2^4$, which by Lemma \ref{newE4-props}(a) has 
type I.  Let $\cale'$ be the set of all subgroups of $S$ which are 
elementary abelian of rank 4, contain $U$, and are not in $\calc$.  By Lemma 
\ref{newE4-props}(e), for any $\varphi\in\Iso_{\Gamma}(E',E'')$ and any 
$E'\in\cale'$, $E''\defeq\varphi(E')\in\cale'$, and $\varphi$ sends 
$x_{\calc}(E')$ to $x_{\calc}(E'')$.  The same holds for 
$\varphi\in\Iso_{\Spin}(E',E'')$ by definition of the elements 
$x_{\calc}(-)$ (\cite[Proposition A.9]{Solnew}).  Since $C_S(E)\in\cale'$, this 
shows that all elements of $\Aut_{\calf}(C_S(E))$ send the element 
$x_{\calc}(C_S(E))$ to itself.  By \cite[Proposition A.9(c)]{Solnew}, 
$\Aut_{\Spin}(C_S(E))$ is the group of automorphisms which are the 
identity on the rank two subgroup $\gen{x_{\calc}(C_S(E)),z}$; and (1) now 
follows.
\end{proof}

The proof of the following lemma is essentially unchanged.  

\begin{Lem}{\rm\cite[Lemma 2.10]{Solnew}}\label{newpoint(e)}\qua
Fix $n\ge1$, and let $E,E'\le{}S(q^n)$ be two elementary abelian subgroups 
of rank three which contain $U$, and which are $\Gamma_n$--conjugate.  Then 
$E$ and $E'$ are $\Spin_7(q^n)$--conjugate.
\end{Lem}

\begin{proof}  Consider the sets
        $$ \calj_1 = \bigl\{X\in SL_2(q^n) \,\big|\, X^2=-I \bigr\} $$
and
        $$ \calj_2 = \bigl\{X\in SL_2(q^{2n}) \,\big|\, \psi^{q^n}(X)=-X,\ 
        X^2=-I \bigr\}. $$
Here, as usual, $\psi^{q^n}$ is induced by the field automorphism 
$(x\mapsto{}x^{q^n})$.  It was shown in the proof of \cite[Lemma 2.10]{Solnew} 
that all elements in each of these sets are $SL_2(q)$--conjugate to each 
other.

Since $E$ and $E$ contain $U$, $E,E'\le{}C_{\Spin_7(q^n)}(U)$.  
By \cite[Proposition 2.5(a)]{Solnew},
        $$ C_{\Spin_7(q^n)}(U)=H(q^n)
        \defeq \omega(SL_2(q^\infty)^3)\cap\Spin_7(q^n). $$
Thus 
        $$ E=\gen{z,z_1,\trp[X_1,X_2,X_3]} \qquad\textup{and}\qquad
        E'=\gen{z,z_1,\trp[X'_1,X'_2,X'_3]}, $$
where the $X_i$ are all in $\calj_1$ or all in $\calj_2$, and similarly 
for the $X'_i$.  Also, $E$ and $E'$ are $\Gamma_n$--conjugate, and each 
element of $\Gamma_n$ leaves $U=\gen{z,z_1}$ and $\omega(SL_2(q^n)^3)$ 
invariant.  Hence either $E$ and $E'$ are both contained in 
$\omega(SL_2(q^n)^3)$, in which case the $X_i$ and $X'_i$ are all in 
$\calj_1$; or neither is contained in $\omega(SL_2(q^n)^3)$, in which case 
the $X_i$ and $X'_i$ are all in $\calj_2$.  This shows that the $X_i$ and 
$X'_i$ are all $SL_2(q^n)$--conjugate, and so $E$ and $E'$ are 
$\Spin_7(q^n)$--conjugate.
\end{proof}

We are now ready to prove:

\begin{Prop}{\rm\cite[Proposition 2.11]{Solnew}}\label{newF_Sol(fqbar)}\qua
For a fixed odd prime power $q$, let 
$S(q^n)\le{}S(q^\infty)\le\Spin_7(q^\infty)$ be as defined above.  Let 
$z\in{}Z(\Spin_7(q^\infty))$ be the central element of order 2.  Then for 
each $n$, $\calf_n=\calf_{\Sol}(q^n)$ is saturated as a fusion system over 
$S(q^n)$, and satisfies the following conditions:
\begin{enumerate}  
\item For all $P,Q\le{}S(q^n)$ which contain $z$, if $\alpha\in\Hom(P,Q)$ is 
such that $\alpha(z)=z$, then $\alpha\in\Hom_{\calf_n}(P,Q)$ if and only if 
$\alpha\in\Hom_{\Spin_7(q^n)}(P,Q)$.
\item $C_{\calf_n}(z)=\calf_{S(q^n)}(\Spin_7(q^n))$ as fusion systems over 
$S(q^n)$.
\item All involutions of $S(q^n)$ are $\calf_n$--conjugate.
\end{enumerate}
Furthermore, $\calf_m\subseteq\calf_n$ for $m|n$.  The union of 
the $\calf_n$ is thus a category $\calf_{\Sol}(q^\infty)$ 
whose objects are the finite subgroups of $S(q^\infty)$.
\end{Prop}

\begin{proof}  We apply Proposition \ref{newconstruct-sFs}, where $p=2$, 
$G=\Spin_7(q^n)$, $S=S(q^n)$, $Z=\gen{z}=Z(G)$, $U=\gen{z,z_1}$, 
$C_G(U)=H(q^n)$, $S_0=C_S(U)=S_0(q^n)$, and $\Gamma=\Gamma_n\le\Aut(S_0)$. 
By Lemma \ref{newgamma-pres.fus.}, $\widehat{\gamma}_u$ preserves 
$H(q^n)$--fusion in $S_0$.  Since $\Gamma$ is generated by 
$\widehat{\gamma}_u$ and certain automorphisms of $H(q^n)$, this shows 
that all automorphisms in $\Gamma$ preserve $H(q^n)$--fusion.

Condition (a) in Proposition \ref{newconstruct-sFs} (all noncentral involutions 
in $G$ are conjugate) holds since all subgroups in $\cale_2$ are conjugate 
(\cite[Proposition A.8]{Solnew}), and condition (b) holds by definition of 
$\Gamma$.  Condition (c) holds since
        $$ \{\gamma\in\Gamma\,|\,\gamma(z)=z\} = 
        \Inn(S_0(q^n)){\cdot}\gen{c_\tau} = \Aut_{N_G(U)}(S_0(q^n)) $$
by \cite[Proposition 2.5(b)]{Solnew}.  Condition (d) was shown in Proposition 
\ref{newrk3}, and condition (e) in Lemma \ref{newpoint(e)}.  So by Proposition 
\ref{newconstruct-sFs}, $\calf_n$ is a saturated fusion system, and 
$C_{\calf_n}(Z)=\calf_{S(q^n)}(\Spin_7(q^n))$.  

The proofs of the other statements remain unchanged.
\end{proof}

In Section 3 of \cite{Solnew}, these corrections affect only the proof of 
Lemma 3.1.  In that proof, the groups $E_{100}$ and $E_{001}$ 
are not $\Gamma_1$--conjugate under the new definitions; instead $E_{100}$ 
is $\Gamma_1$--conjugate to a subgroup which is $\Spin_7(q)$--conjugate to 
$E_{001}$ (and hence the two are $\calf$--conjugate).  Also, when showing 
that $\autf(E_{001})$ is the group of all automorphisms which fix 
$z=x_{\calc}(E_{001})$, it is important to know that all 
$\Gamma_1$--isomorphisms between subgroups in that conjugacy class preserve 
the elements $x_{\calc}(-)$ (as shown in Lemma \ref{newE4-props}(e)), and not 
just that $\Gamma_1$--automorphisms of $E_{001}$ do so.  

The changes do not affect the later sections.  Just note that we are able 
to consider $\calf_{\Sol}(q^n)$ as a subcategory of $\calf_{\Sol}(q^{kn})$ 
for $k>1$, because they were both chosen using the same ``correction 
factor'' $u\in(\ztwo)^*$ in the definitions of $\Gamma_n$ and 
$\Gamma_{kn}$.

\section{Proof of Proposition \ref{newconstruct-sFs}}

Proposition \ref{newconstruct-sFs} follows from Lemmas 1.3, 1.4, and 1.5 in 
\cite{Solnew}, once they are restated to assume the hypotheses of this new 
proposition, and not those of \cite[Proposition 1.2]{Solnew}.  The 
only one of these lemmas whose proof is affected by the change in 
hypotheses is Lemma 1.4, and so we restate and reprove it here.  

\begin{Lem}\label{new9.2.2}
Assume the hypotheses of Proposition \ref{newconstruct-sFs}, and let 
	$$\calf=\gen{\calf_S(G);\calf_{S_0}(\Gamma)}$$ 
be the fusion system 
generated by $G$ and $\Gamma$.  Then for all $P,P'\le{}S$ which contain 
$Z$,
        $$ \{\varphi\in\homf(P,P') \,|\, \varphi(Z)=Z \} = \Hom_G(P,P'). $$
\end{Lem}

\begin{proof} Upon replacing $P'$ by $\varphi(P)\le{}P'$, we can assume 
that $\varphi$ is an isomorphism, and thus that it factors as a composite 
of isomorphisms 
        $$ P=P_0 \RIGHT2{\varphi_1}{\cong} P_1 \RIGHT2{\varphi_2}{\cong} 
        P_2 \RIGHT2{\varphi_3}{\cong} \cdots \RIGHT2{\varphi_{k-1}}{\cong} 
        P_{k-1} \RIGHT2{\varphi_k}{\cong} P_k=P', $$
where for each $i$, $\varphi_i\in\Hom_G(P_{i-1},P_i)$ or 
$\varphi_i\in\Hom_\Gamma(P_{i-1},P_i)$.  Let $Z_i\le{}Z(P_i)$ be the 
subgroups of order $p$ such that $Z_0=Z_k=Z$ and $Z_i=\varphi_i(Z_{i-1})$. 

To simplify the discussion, we say that a morphism in $\calf$ is of type 
$(G)$ if it is given by conjugation by an element of $G$, and of type 
$(\Gamma)$ if it is the restriction of an automorphism in $\Gamma$.  More 
generally, we say that a morphism is of type $(G,\Gamma)$ if it is the 
composite of a morphism of type $(G)$ followed by one of type $(\Gamma)$, 
etc.  We regard $\Id_P$, for all $P\le{}S$, to be of both types, even
if $P\nleq{}S_0$.  By definition, if any nonidentity isomorphism is of 
type $(\Gamma)$, then its source and image are both contained in 
$S_0=C_S(U)$.  In particular, if $P_i\nleq{}S_0$ for any $0<i<k$, then 
$\varphi_i$ and $\varphi_{i+1}$ are both of type $(G)$, and we can remove 
it and replace the two morphisms by their composite.  We can thus assume 
that $[P_i,U]=1$ for all $0<i<k$, and hence that $Z_i\le{}Z(P_iU)$ for all 
$i$.

For each $i$, using \cite[Lemma 1.3]{Solnew}, choose some 
$\psi_i\in\Hom_{\calf}(P_iU,S)$ such that $\psi_i(Z_i)=Z$.  More 
precisely, using points (1) and (2) in \cite[Lemma 1.3]{Solnew}, we can choose 
$\psi_i$ to be of type $(\Gamma)$ if $Z_i\le{}U$ (the inclusion if 
$Z_i=Z$), and to be of type $(G,\Gamma)$ if $Z_i\nleq{}U$.  Set 
$P'_i=\psi_i(P_i)$. To keep track of the effect of morphisms on the 
subgroups $Z_i$, we write them as morphisms between pairs, as shown below. 
Thus, $\varphi$ factors as a composite of isomorphisms
        $$ (P'_{i-1},Z) \Right4{\psi_{i-1}^{-1}} (P_{i-1},Z_{i-1}) 
        \Right4{\varphi_i} (P_i,Z_i) \Right4{\psi_i} (P'_i,Z).  $$
If $\varphi_i$ is of type $(G)$, then this composite (after replacing 
adjacent morphisms of the same type by their composite) is of type 
$(\Gamma,G,\Gamma)$.  If $\varphi_i$ is of type $(\Gamma)$, then the 
composite is again of type $(\Gamma,G,\Gamma)$ if either $Z_{i-1}\le{}U$ 
or $Z_i\le{}U$, and is of type $(\Gamma,G,\Gamma,G,\Gamma)$ if neither 
$Z_{i-1}$ nor $Z_i$ is contained in $U$.  So we are reduced to assuming
that $\varphi$ is of one of these two forms.

\smallskip

\noindent\textbf{Case 1}\qua  Assume first that $\varphi$ is of type 
$(\Gamma,G,\Gamma)$; ie, a composite of isomorphisms of the form
        $$ (P_0,Z) \RIGHT3{\varphi_1}{(\Gamma)} (P_1,Z_1) 
        \RIGHT3{\varphi_2}{(G)} (P_2,Z_2) 
        \RIGHT3{\varphi_3}{(\Gamma)} (P_3,Z). $$
Then $Z_1=Z$ if and only if $Z_2=Z$ because $\varphi_2$ is of type $(G)$. 
If $Z_1=Z_2=Z$, then $\varphi_1$ and $\varphi_3$ are of type $(G)$ by 
Proposition \ref{newconstruct-sFs}(c), 
and the result follows.  

If $Z_1\ne{}Z\ne{}Z_2$, then $U=ZZ_1=ZZ_2$, and thus $\varphi_2(U)=U$.  
Neither $\varphi_1$ nor $\varphi_3$ can be the identity, so 
$P_i\le{}S_0=C_S(U)$ for all $i$ by definition of $\Hom_\Gamma(-,-)$.  Let 
$\gamma_1,\gamma_3\in\Gamma$ be such that $\varphi_i=\gamma_i|_{P_{i-1}}$, 
and let $g\in{}N_G(U)$ be such that $\varphi_2=c_g$.  Then 
$gC_S(U)g^{-1}\in\sylp{C_G(U)}$, so there is $h\in{}C_G(U)$ such that 
$hg\in{}N(C_S(U))$.  Then $c_{hg}\in\Gamma$ by Proposition 
\ref{newconstruct-sFs}(c).  

Since $\gamma_3$ preserves $C_G(U)$--fusion among subgroups of $S_0$, there 
is $g'\in{}C_G(U)$ such that 
$c_{g'}=\gamma_3\circ{}c_h^{-1}\circ\gamma_3^{-1}$.  Thus
$\varphi$ is the composite
        $$ P_0 \RIGHT3{\varphi_1}{(\Gamma)} P_1 
        \RIGHT3{c_{hg}}{(\Gamma)} hP_2h^{-1}
        \RIGHT3{\gamma_3}{(\Gamma)} \gamma_3(hP_2h^{-1})
        \RIGHT3{c_{g'}}{(G)} P_3. $$
The composite of the first three isomorphisms is of type $(\Gamma)$ and 
sends $Z$ to $g'{}^{-1}Zg'=Z$, hence is of type $(G)$ by Proposition 
\ref{newconstruct-sFs}(c) again, and so $\varphi$ is also of type $(G)$.  

\smallskip

\noindent\textbf{Case 2}\qua  Assume now that $\varphi$ is of type 
$(\Gamma,G,\Gamma,G,\Gamma)$; more precisely, that it is a composite of the 
form
	{ \begin{multline*}  
         \def\2#1#2{\!\RIGHT2{#1}{#2}}
        (P_0,Z) \2{\varphi_1}{(\Gamma)} (P_1,Z_1) 
        \2{\varphi_2}{(G)} (P_2,Z_2) 
        \2{\varphi_3}{(\Gamma)} (P_3,Z_3) \\
        \def\2#1#2{\!\RIGHT2{#1}{#2}}
        \2{\varphi_4}{(G)} (P_4,Z_4) 
        \2{\varphi_5}{(\Gamma)} (P_5,Z), 
	\end{multline*} }
where $Z_2,Z_3\nleq{}U$.  Then $Z_1,Z_4\le{}U$ and are distinct from $Z$, 
and the groups $P_0,P_1,P_4,P_5$ all contain $U$ since $\varphi_1$ and 
$\varphi_5$ (being of type $(\Gamma)$) leave $U$ invariant.  In 
particular, $P_2$ and $P_3$ contain $Z$, since $P_1$ and $P_4$ do and 
$\varphi_2,\varphi_4$ are of type $(G)$.  We can also assume that 
$U\le{}P_2,P_3$, since otherwise $P_2\cap{}U=Z$ or $P_3\cap{}U=Z$, 
$\varphi_3(Z)=Z$, and hence $\varphi_3$ is of type $(G)$ by Proposition 
\ref{newconstruct-sFs}(c) again.  Finally, we assume that 
$P_2,P_3\le{}S_0=C_S(U)$, since otherwise $\varphi_3=\Id$.  

Let $E_i\le{}P_i$ be the rank three elementary abelian subgroups defined 
by the requirements that $E_2=UZ_2$, $E_3=UZ_3$, and 
$\varphi_i(E_{i-1})=E_i$.  In particular, $E_i\le{}Z(P_i)$ for $i=2,3$ 
(since $Z_i\le{}Z(P_i)$, and $U\le{}Z(P_i)$ by the above remarks); and 
hence $E_i\le{}Z(P_i)$ for all $i$.  Also, $U=ZZ_4\le\varphi_4(E_3)=E_4$ 
since $\varphi_4(Z)=Z$, and thus $U=\varphi_5(U)\le{}E_5$.  Via similar 
considerations for $E_0$ and $E_1$, we see that $U\le{}E_i$ for all $i$.

Set $H=C_G(U)$ for short.  Let $\cale_3$ be the set of all elementary 
abelian subgroups $E\le{}S$ of rank three which contain $U$, and let 
$\cale^*_3\subseteq\cale_3$ be the set of all $E\in\cale_3$ such that
$C_S(E)\in\sylp{C_G(E)}$.  Then for all $E\in\cale_3$,
        $$ C_S(E)=C_{S_0}(E) \qquad\textup{and}\qquad C_G(E)=C_H(E) $$ 
since $E\ge{}U$ (and $S_0=C_S(U)$).  Thus $\cale^*_3$ is the set of all 
subgroups $E\in\cale_3$ which are fully centralized in the fusion system 
$\calf_{S_0}(H)$, and so each subgroup in $\cale_3$ is $H$--conjugate to a 
subgroup in $\cale^*_3$.  

Let $\gamma_1,\gamma_3,\gamma_5\in\Gamma$ be such that 
$\varphi_i=\gamma_i|_{P_{i-1}}$ for odd $i$.  Let $g_2,g_4\in G$ be such 
that $\varphi_i$ is conjugation by $g_i$ for $i=2,4$.  We will construct a 
commutative diagram of the following form
        \begin{diagram}[w=30pt]
        (P_0,E_0) & \rTo^{\gamma_1} & (P_1,E_1) & \rTo^{c_{g_2}} & 
        (P_2,E_2) & \rTo^{\gamma_3} & (P_3,E_3) & \rTo^{c_{g_4}} & 
        (P_4,E_4) & \rTo^{\gamma_5} & (P_5,E_5) \\
        \dTo<{c_{a_0}} && \dTo<{c_{a_1}} && \dTo<{c_{a_2}} && 
        \dTo<{c_{a_3}} && \dTo<{c_{a_4}} && \dTo<{c_{a_5}} \\
        (C'_0,E'_0) & \rTo^{\gamma_1} & (C'_1,E'_1) & \rTo^{c_{h_2}} & 
        (C'_2,E'_2) & \rTo^{\gamma_3} & (C'_3,E'_3) & \rTo^{c_{h_4}} & 
        (C'_4,E'_4) & \rTo^{\gamma_5} & (C'_5,E'_5) 
        \end{diagram}
where $E'_i\in\cale^*_3$, $C'_i=C_S(E'_i)$, $h_2,h_4\in{}G$, and 
$a_i\in{}H=C_G(U)$.  To do this, first choose $a'_0,a'_2,a'_4\in{}H$ such 
that $E'_i\defeq{}a'_iE_ia'_i{}^{-1}\in\cale^*_3$ for $i=0,2,4$, and set 
$E'_i=\gamma_i(E'_{i-1})$ for $i=1,3,5$.  Then 
$C_S(E'_i)=\gamma_i(C_S(E'_{i-1}))$ for $i=1,3,5$, so 
$C_S(E'_i)\in\sylp{C_H(E'_i)}$ and $E'_i\in\cale^*_3$ for all $i$.  So we 
can choose $x_0\in{}C_H(E'_0)$ such that 
$x_0(a'_0P_0a'_0{}^{-1})x_0^{-1}\le{}C_S(E'_0)$, and set 
$a_0=x_0a'_0\in{}H$.  Since $\gamma_1\in\Aut(S_0)$ preserves $H$--fusion, 
there is $a_1\in{}H$ which makes the first square in the above diagram 
commute.  Now choose $x_2\in{}C_H(E'_2)$ such that
        $$ x_2\bigl((a'_2g_2a_1^{-1})C_S(E'_1)(a'_2g_2a_1^{-1})^{-1}
        \bigr)x_2^{-1} \le C_S(E'_2), $$
and set $a_2=x_2a'_2$ and $h_2=a_2g_2a_1^{-1}$.  Upon continuing this 
procedure we obtain the above diagram.

Let $\widehat{\varphi}\in\Iso_\calf(C_S(E'_0),C_S(E'_5))$ be the composite 
of the morphisms in the bottom row of the above diagram.  Then 
$\widehat{\varphi}(Z)=Z$, since $\varphi$ and the $c_{a_i}$ all send $Z$ 
to itself.  By Proposition \ref{newconstruct-sFs}(e), the rank three 
subgroups $E'_i$ are all $G$--conjugate to each other.  Choose $g\in{}G$ 
such that $gE'_5g^{-1}=E'_0$.  Then $g{\cdot}C_S(E'_5){\cdot}g^{-1}$ and 
$C_S(E'_0)$ are both Sylow $p$--subgroups of $C_G(E'_0)$, so there is 
$h\in{}C_G(E'_0)$ such that $(hg)C_S(E'_5)(hg)^{-1}=C_S(E'_0)$.  By 
Proposition \ref{newconstruct-sFs}(d), 
$c_{hg}\circ\widehat{\varphi}\in\Aut_\calf(C_S(E'_0))$ is of type $(G)$, 
so $\widehat{\varphi}$ is of type $(G)$, and hence 
$\varphi=c_{a_5}^{-1}\circ\widehat{\varphi}\circ{}c_{a_0}$ is also of type 
$(G)$.
\end{proof}

\end{document}

\end{document}